\documentstyle[11pt,amstex, amscd, graphicx]{article}

\begin{document}

\newtheorem{theorem}{Theorem}[section]
\newtheorem{prop}[theorem]{Proposition}
\newtheorem{lemma}[theorem]{Lemma}
\newtheorem{cor}[theorem]{Corollary}
\newtheorem{definition}[theorem]{Definition}
\newtheorem{conj}[theorem]{Conjecture}
\newtheorem{claim}[theorem]{Claim}

\newcommand{\boundary}{\partial}
\newcommand{\C}{{\mathbb C}}
\newcommand{\integers}{{\mathbb Z}}
\newcommand{\natls}{{\mathbb N}}
\newcommand{\ratls}{{\mathbb Q}}
\newcommand{\reals}{{\mathbb R}}
\newcommand{\proj}{{\mathbb P}}
\newcommand{\lhp}{{\mathbb L}}
\newcommand{\tube}{{\mathbb T}}
\newcommand{\cusp}{{\mathbb P}}
\newcommand\AAA{{\mathcal A}}
\newcommand\BB{{\mathcal B}}
\newcommand\CC{{\mathcal C}}
\newcommand\DD{{\mathcal D}}
\newcommand\EE{{\mathcal E}}
\newcommand\FF{{\mathcal F}}
\newcommand\GG{{\mathcal G}}
\newcommand\HH{{\mathcal H}}
\newcommand\II{{\mathcal I}}
\newcommand\JJ{{\mathcal J}}
\newcommand\KK{{\mathcal K}}
\newcommand\LL{{\mathcal L}}
\newcommand\MM{{\mathcal M}}
\newcommand\NN{{\mathcal N}}
\newcommand\OO{{\mathcal O}}
\newcommand\PP{{\mathcal P}}
\newcommand\QQ{{\mathcal Q}}
\newcommand\RR{{\mathcal R}}
\newcommand\SSS{{\mathcal S}}
\newcommand\TT{{\mathcal T}}
\newcommand\UU{{\mathcal U}}
\newcommand\VV{{\mathcal V}}
\newcommand\WW{{\mathcal W}}
\newcommand\XX{{\mathcal X}}
\newcommand\YY{{\mathcal Y}}
\newcommand\ZZ{{\mathcal Z}}
\newcommand\CH{{\CC\HH}}
\newcommand\MF{{\MM\FF}}
\newcommand\PMF{{\PP\kern-2pt\MM\FF}}
\newcommand\ML{{\MM\LL}}
\newcommand\PML{{\PP\kern-2pt\MM\LL}}
\newcommand\GL{{\GG\LL}}
\newcommand\Pol{{\mathcal P}}
\newcommand\half{{\textstyle{\frac12}}}
\newcommand\Half{{\frac12}}
\newcommand\Mod{\operatorname{Mod}}
\newcommand\Area{\operatorname{Area}}
\newcommand\ep{\epsilon}
\newcommand\hhat{\widehat}
\newcommand\Proj{{\mathbf P}}
\newcommand\U{{\mathbf U}}
 \newcommand\Hyp{{\mathbf H}}
\newcommand\D{{\mathbf D}}
\newcommand\Z{{\mathbb Z}}
\newcommand\R{{\mathbb R}}
\newcommand\Q{{\mathbb Q}}
\newcommand\E{{\mathbb E}}
\newcommand\til{\widetilde}
\newcommand\length{\operatorname{length}}
\newcommand\tr{\operatorname{tr}}
\newcommand\gesim{\succ}
\newcommand\lesim{\prec}
\newcommand\simle{\lesim}
\newcommand\simge{\gesim}
\newcommand{\simmult}{\asymp}
\newcommand{\simadd}{\mathrel{\overset{\text{\tiny $+$}}{\sim}}}
\newcommand{\ssm}{\setminus}
\newcommand{\diam}{\operatorname{diam}}
\newcommand{\pair}[1]{\langle #1\rangle}
\newcommand{\T}{{\mathbf T}}
\newcommand{\inj}{\operatorname{inj}}
\newcommand{\pleat}{\operatorname{\mathbf{pleat}}}
\newcommand{\short}{\operatorname{\mathbf{short}}}
\newcommand{\vertices}{\operatorname{vert}}
\newcommand{\collar}{\operatorname{\mathbf{collar}}}
\newcommand{\bcollar}{\operatorname{\overline{\mathbf{collar}}}}
\newcommand{\I}{{\mathbf I}}
\newcommand{\tprec}{\prec_t}
\newcommand{\fprec}{\prec_f}
\newcommand{\bprec}{\prec_b}
\newcommand{\pprec}{\prec_p}
\newcommand{\ppreceq}{\preceq_p}
\newcommand{\sprec}{\prec_s}
\newcommand{\cpreceq}{\preceq_c}
\newcommand{\cprec}{\prec_c}
\newcommand{\topprec}{\prec_{\rm top}}
\newcommand{\Topprec}{\prec_{\rm TOP}}
\newcommand{\fsub}{\mathrel{\scriptstyle\searrow}}
\newcommand{\bsub}{\mathrel{\scriptstyle\swarrow}}
\newcommand{\fsubd}{\mathrel{{\scriptstyle\searrow}\kern-1ex^d\kern0.5ex}}
\newcommand{\bsubd}{\mathrel{{\scriptstyle\swarrow}\kern-1.6ex^d\kern0.8ex}}
\newcommand{\fsubeq}{\mathrel{\raise-.7ex\hbox{$\overset{\searrow}{=}$}}}
\newcommand{\bsubeq}{\mathrel{\raise-.7ex\hbox{$\overset{\swarrow}{=}$}}}
\newcommand{\tw}{\operatorname{tw}}
\newcommand{\base}{\operatorname{base}}
\newcommand{\trans}{\operatorname{trans}}
\newcommand{\rest}{|_}
\newcommand{\bbar}{\overline}
\newcommand{\UML}{\operatorname{\UU\MM\LL}}
\newcommand{\EL}{\mathcal{EL}}
\newcommand{\tsum}{\sideset{}{'}\sum}
\newcommand{\tsh}[1]{\left\{\kern-.9ex\left\{#1\right\}\kern-.9ex\right\}}
\newcommand{\Tsh}[2]{\tsh{#2}_{#1}}
\newcommand{\qeq}{\mathrel{\approx}}
\newcommand{\Qeq}[1]{\mathrel{\approx_{#1}}}
\newcommand{\qle}{\lesssim}
\newcommand{\Qle}[1]{\mathrel{\lesssim_{#1}}}
\newcommand{\simp}{\operatorname{simp}}
\newcommand{\vsucc}{\operatorname{succ}}
\newcommand{\vpred}{\operatorname{pred}}
\newcommand\fhalf[1]{\overrightarrow {#1}}
\newcommand\bhalf[1]{\overleftarrow {#1}}
\newcommand\sleft{_{\text{left}}}
\newcommand\sright{_{\text{right}}}
\newcommand\sbtop{_{\text{top}}}
\newcommand\sbot{_{\text{bot}}}
\newcommand\sll{_{\mathbf l}}
\newcommand\srr{_{\mathbf r}}
\newcommand\geod{\operatorname{\mathbf g}}
\newcommand\mtorus[1]{\boundary U(#1)}
\newcommand\A{\mathbf A}
\newcommand\Aleft[1]{\A\sleft(#1)}
\newcommand\Aright[1]{\A\sright(#1)}
\newcommand\Atop[1]{\A\sbtop(#1)}
\newcommand\Abot[1]{\A\sbot(#1)}
\newcommand\boundvert{{\boundary_{||}}}
\newcommand\storus[1]{U(#1)}
\newcommand\Momega{\omega_M}
\newcommand\nomega{\omega_\nu}
\newcommand\twist{\operatorname{tw}}
\newcommand\modl{M_\nu}
\newcommand\MT{{\mathbb T}}
\newcommand\Teich{{\mathcal T}}
\renewcommand{\Re}{\operatorname{Re}}
\renewcommand{\Im}{\operatorname{Im}}

\title{Cannon-Thurston Maps for Surface Groups I: Amalgamation
  Geometry and Split Geometry}

\author{Mahan Mj.\\ \date{} 
 RKM Vidyamandira and RKMVERI,  \\
Belur Math, WB-711202, India \\
  email: brmahan@@gmail.com}

\maketitle

\begin{abstract}

We introduce the notion of manifolds of
amalgamation geometry  and its
generalisation,  split geometry.
We show that the limit
set of any surface  group of  split geometry
 is locally connected, by constructing a natural Cannon-Thurston map.

\smallskip

\begin{center}


\end{center}

\end{abstract}

\tableofcontents

\section{Introduction}

\subsection{Statement of Results}
In this paper and its successor \cite{mahan-split}, 
we continue our study of Cannon-Thurston maps and limit sets of
Kleinian groups
initiated in \cite{mitra-trees}, 
\cite{brahma-pared} and \cite{brahma-ibdd}.
Several questions and conjectures have been made in this context by
different authors:

\noindent $\bullet 1$ In Section 6 of \cite{CT}, Cannon and Thurston
raise the following problem: \\
{\bf Question:}
Suppose a closed surface group $\pi_1 (S)$ acts freely and properly
discontinuously on ${\Bbb{H}}^3$ by isometries. Does the inclusion
$\tilde{i} : \widetilde{S} \rightarrow {\Bbb{H}}^3$ extend
continuously to the boundary?

The authors of \cite{CT} point out that for a simply degenerate group,
this is equivalent to asking if the limit set is locally connected.

\smallskip

\noindent $\bullet 2$ In \cite{ctm-locconn}, McMullen makes the
following more general conjecture: \\
{\bf Conjecture:}
For any hyperbolic 3-manifold $N$ with finitely generated fundamental
group, there exists a continuous, $\pi_1(N)$-equivariant map \\
\begin{center}
$F: \partial \pi_1 (N) \rightarrow \Lambda \subset S^2_{\infty}$

\end{center}

where the boundary $\partial \pi_1(N)$ is constructed by scaling the
metric on the Cayley graph of $\pi_1 (N)$ by the conformal factor of
$d(e,x)^{-2}$, then taking the metric completion. (cf. Floyd
\cite{Floyd})

\smallskip

\noindent $\bullet 3$ The  author raised the following question
in his thesis \cite{mitra-thesis} (see also \cite{bestvinahp}): \\
{\bf Question:}
Let $G$ be a hyperbolic group in the sense of Gromov acting freely and
properly discontinuously by isometries on a hyperbolic metric space
$X$. Does the inclusion of the Cayley graph $i: \Gamma_G \rightarrow
X$ extend continuously to the (Gromov) compactifications? \\
A similar question may be asked for relatively hyperbolic groups (in
the sense of Gromov \cite{gromov-hypgps} and 
Farb \cite{farb-relhyp}). 

The question for relatively hyperbolic groups unifies
 all the above questions and
conjectures.

\smallskip

In this paper we introduce the notion of what we call {\it
  amalgamation geometry} which is, in a way, a considerable
  generalisation of the notion of {\it i-bounded geometry} introduced
  in \cite{brahma-ibdd}. We then generalise it by weakening the
 hypothesis to the notion of {\it split geometry}.
A crucial step in this paper is to prove:

\smallskip

\noindent {\bf Theorems \ref{crucial-split} and  \ref{crucial-punct-split}:} Let
$\rho : \pi_1(S) \rightarrow 
PSL_2(C)$ be a faithful representation of a 
surface group with or without punctures,
and without accidental parabolics. Let $M =
{{\Bbb{H}}^3}/{\rho (\pi_1 (S))}$ be of split geometry. 
Let $i$ be an embedding of $S$
in $M$ that induces a homotopy equivalence. Then the embedding
$\tilde{i} : \widetilde{S} \rightarrow \widetilde{M} = {\Bbb{H}}^3$
extends continuously to a map $\hat{i}: {\Bbb{D}}^2 \rightarrow
{\Bbb{D}}^3$. Further, the limit set of ${\rho (\pi_1 (S))}$ is
locally connected.

\smallskip

In fact our methods prove the following considerably stronger result by
combining the techniques of this paper with those of \cite{mitra-trees} and
\cite{brahma-pared}. This is a partial  affirmation of McMullen's conjecture
above. 

\smallskip

\noindent {\bf Theorem \ref{main3} :}
Suppose that $N^h \in H(M,P)$ is a hyperbolic structure of {\em
split geometry} 
on a pared manifold $(M,P)$ with incompressible boundary $\partial_0 M$. Let
$M_{gf}$ denote a geometrically finite hyperbolic structure adapted
to $(M,P)$. Then the map  $i: \widetilde{M_{gf}}
\rightarrow \widetilde{N^h}$ extends continuously to the boundary
$\hat{i}: \widehat{M_{gf}}
\rightarrow \widehat{N^h}$. If $\Lambda$ denotes the limit set of
$\widetilde{M}$, then $\Lambda$ is locally connected.

\smallskip

In \cite{mahan-split}, we shall
 show that the Minsky model is of split geometry. 
Combining this with Theorems \ref{crucial-split} and
 \ref{crucial-punct-split}, 
we shall get

\smallskip

\noindent {\bf Theorem \cite{mahan-split}:}
Let $\rho$ be a representation of a surface group $H$ (corresponding
to the surface $S$) into
$PSl_2(C)$ without accidental parabolics. Let $M$ denote the (convex
core of) ${\Bbb{H}}^3 / \rho 
(H)$.  Further suppose that $i: S \rightarrow M$, taking
 parabolic to parabolics, induces a homotopy
equivalence.
  Then the inclusion
  $\tilde{i} : \widetilde{S} \rightarrow \widetilde{M}$ extends continuously
  to a map 
  $\hat{i} : \widehat{S} \rightarrow \widehat{M}$. Hence the limit set
  of $\widetilde{S}$ is locally connected.

\smallskip

Again, combining the Minsky model with Theorem \ref{main3}, we shall 
get\\

\noindent {\bf Theorem \cite{mahan-split}:}
Suppose that $N^h \in H(M,P)$ is a hyperbolic structure 
on a pared manifold $(M,P)$ with incompressible boundary $\partial_0 M$. Let
$M_{gf}$ denote a geometrically finite hyperbolic structure adapted
to $(M,P)$. Then the map  $\tilde{i}: \widetilde{M_{gf}}
\rightarrow \widetilde{N^h}$ extends continuously to the boundary
$\hat{i}: \widehat{M_{gf}}
\rightarrow \widehat{N^h}$. If $\Lambda$ denotes the limit set of
$\widetilde{M}$, then $\Lambda$ is locally connected.

\subsection{History and Present State of the Problem}

The first major result that started this entire program was Cannon and
Thurston's result  \cite{CT} for hyperbolic 3-manifolds fibering over the
circle with fiber a closed surface group. 

This was generalised by Minsky who proved the Cannon-Thurston result
for bounded 
geometry Kleinian closed surface groups \cite{minsky-jams}. 

An
alternate approach (purely in terms of coarse geometry ignoring all
local information) was given by the author in \cite{mitra-trees}
generalising the results of both  Cannon-Thurston and Minsky. We proved the
Cannon-Thurston result for hyperbolic 3-manifolds of bounded geometry
without parabolics and with freely indecomposable fundamental group. A
different approach based on Minsky's work was
given by Klarreich \cite{klarreich}.

Bowditch \cite{bowditch-ct}
\cite{bowditch-stacks} proved the Cannon-Thurston result
for punctured surface Kleinian groups of
bounded geometry.
In  \cite{brahma-pared} we gave an
alternate proof of Bowditch's results and simultaneously generalised the
results of Cannon-Thurston, Minsky, Bowditch, and those of 
\cite{mitra-trees} to
all 3 manifolds of bounded geometry whose cores are incompressible
away from cusps.  The proof has the
advantage that it reduces to a proof for manifolds without parabolics
when the 3 manifold in question has freely indecomposable fundamental
group and no accidental parabolics.

McMullen \cite{ctm-locconn} proved the Cannon-Thurston result
for punctured torus groups, using
 Minsky's model for these groups \cite{minsky-torus}. 
 In \cite{brahma-ibdd} we identified a large-scale 
coarse geometric structure involved in the Minsky model for punctured
torus groups (and called it {\bf
  i-bounded geometry}). {\em i-bounded geometry} can roughly be
regarded as that geometry of ends where the boundary torii of Margulis
tubes have uniformly bounded diameter. We gave a  proof for models of
{\it i-bounded geometry}. In combination with the methods of
\cite{brahma-pared} this was enough to bring under the same umbrella
all known results on Cannon-Thurston maps for 3 manifolds whose cores are
incompressible away from cusps. In particular,
when $(M,P)$ is the 
pair $S \times I, \delta S \times I$, for $S$ a punctured torus or
four-holed sphere, we gave an alternate proof of
McMullen's
result \cite{ctm-locconn}. 

In this paper, we define  {\it amalgamation geometry}
and prove the Cannon-Thurston result for models of
{\it amalgamation geometry}. We then weaken this assumption
to what we call {\it split geometry} and prove the Cannon-Thurston
property
for such geometries. In \cite{mahan-split}  we shall
show that the Minsky model for
general
simply or totally degenerate surface groups \cite{minsky-elc1}
\cite{minsky-elc2} gives rise to a model of {\em split geometry}. 
 This will allow us to
conclude that all surface groups have the Cannon-Thurston property and
hence have locally connected limit sets. In the sequel to this paper
\cite{mahan-split}, we  show that the Minsky model for surface groups
has split geometry. This proves that surface groups (and more
generally Kleinian groups corresponding to manifolds whose cores are
incompressible away from cusps) have locally-connected limit sets. 

\subsection{Scheme and Outline of the Paper}

We first describe in brief, the philosophy of the proof. Given a
simply degenerate surface ($S$) 
group (without accidental parabolics),
Thurston \cite{Thurstonnotes} proves that a unique ending lamination
$\lambda$ exists. Let $M = {\Hyp}^3/\rho(\pi_1 (S))$.
In this situation, it follows from
\cite{Thurstonnotes} that any sequence of simple closed curves
$\sigma_i$, whose geodesic realizations exit the end, converges to
$\lambda$. Dual to $\lambda$, there exists an $\Bbb{R}$-tree $\TT$ and
a free action of $\pi_1(S)$ on $\TT$. Now, each $\sigma_i$ gives rise
to a splitting of $\pi_1(S)$, and hence an action of $\pi_1(S)$ on a
simplicial tree $\TT_i$. The sequence of these actions converges to
the action of $\pi_1(S)$ on $\TT$ dual to $\lambda$ (see for instance,
\cite{otal-book}). 

The guiding motif of this paper is to find geometric realizations of
this sequence of splittings in terms of contiguous blocks $B_i$ (each
homeomorphic to $S \times I$). By a geometric realization of a
splitting we mean the following:\\
Margulis tubes $T_i$ are chosen, exiting the end of $M$. Let
$\sigma_i$ denote the core geodesic of $T_i$. We require that $T_i$
splits some block $B_i$, i.e. $B_i \setminus T_i$ is homeomorphic to
$(S \setminus A(\sigma_i )) \times I$, where $A(\sigma_i )$ is an
annular neighborhood of a geodesic representative of $\sigma_i$ {\em
  on $S$}. We require further control on the geometry of the
complementary pieces 
$(S \setminus A(\sigma_i )) \times I$. 

For conceptual simplicity, assume the $T_i$'s are
separating. Different degrees of control on the geometry of the pieces
$(S \setminus A(\sigma_i )) \times I$ give rise to different
geometries. Fix a piece of
$(S \setminus A(\sigma_i )) \times I$ and call it $K$.
 It is better to look at the universal cover $\til{B_i}$ and a lift
 $\til{K} \subset \til{B_i}$. We adjoin the lifts of $T_i$ bounding
 $\til{K}$ to $\til{K}$ and call it $K_1$. $K_1$ shall be referred to
 as a {\bf component} of the relevant geometry. \\
{\bf 1) Amalgamation Geometry:} The simplest geometry arising from
this situation is the case where all $K_1$'s are uniformly quasiconvex
{\em in the hyperbolic metric on $\til{M}$}. This is called
amalgamation geometry, and can in brief e described as the geometry in
which all components are uniformly (hyperbolically) quasiconvex. \\
{\bf 2) Graph Amalgamation Geometry:} Amalgamation geometry is too
restrictive. As a first step towards relaxing this hypothesis, we do
not demand that the convex hulls $CH(K_1)$'s be contained in uniformly
bounded neighborhoods of the respective $K_1$'s in the hyperbolic
metric. Instead we construct an auxiliary metric called the{\bf
  graph-metric}. Roughly speaking, the {\em graph-metric} is the
natural simplicial metric on the nerve of the covering of $\til{M}$ by
the components $K_1$. {\em Graph Amalgamation Geometry} is the
condition that the convex hulls $CH(K_1)$'s lie in  uniformly bounded
neighborhood of $K_1$'s in the {\em graph metric}. \\
{\bf 3) Split  Geometry:} So far, we have assumed that each Margulis
tube $T_i$ is contained wholly in a block $B_i$, splitting
it. However, 
as was pointed out to the author by Yair Minsky and Dick Canary, this
is not the most general situation. The $T_i$'s may {\em interlock}. To
take care of this situation, we allow each tube $T_i$ to cut through
(partly or wholly) a uniformly bounded number of blocks. The notions
of complementary components and the graph metric still make sense. The
rsulting geometry is termed {\em split geometry}. 

We shall take one step at a time in this paper, relaxing the
hypothesis in the order above. The additional arguments to be
introduced as we proceed from one geometry to the next (more general)
one will be described as modifications of the core argument relevant
to amalgamation geometry. 

In the sequel \cite{mahan-split}, we shall show that simply degenerate
ends of hyperbolic 3-manifolds enjoy split geometry. \medskip \medskip

\noindent {\bf Outline:}
A brief outline of the paper follows. Section 2 deals with
preliminaries.
We also define {\it amalgamation geometry} via the construction of a model
manifold.

 Section 3 deals with relative hyperbolicity {\it a la} Gromov
 \cite{gromov-hypgps}, Farb \cite{farb-relhyp} and Bowditch
 \cite{bowditch-relhyp}. 

As in \cite{mitra-ct}, \cite{mitra-trees}, \cite{brahma-pared},
\cite{brahma-ibdd}, a crucial part of our
proof proceeds  by constructing a {\it ladder-like} set $B_\lambda \subset
\widetilde{M}$  from a geodesic segment $\lambda \subset
\widetilde{S}$ and then a retraction $\Pi_\lambda$ of $\widetilde{M}$
onto $B_\lambda$. 

In
Section 4, we construct a model geometry for the universal covers of
building blocks  and the relevant geometries (electric and graph
models) that will concern us.

We also construct the paths that go to build up the
{\it ladder-like set}
$B_\lambda$. We further construct the
restriction of the retraction $\Pi_\lambda$ to  blocks and show
that the retraction does not increase distances much. 

In Section 5, we put the blocks and retractions together (by adding them one
on top of another) to  build the {\it ladder-like} $B_\lambda$
and  prove the main
technical  theorem - the existence of 
of a retract $\Pi_\lambda$ of $\widetilde{M}$ onto $B_\lambda$. This
shows that $B_\lambda$ is quasiconvex in $\widetilde{M}$ equipped with
a model pseudometric.

In Section 6, we put together the ingredients from Sections 2, 3, 4
and 5 to prove the existence of a Cannon-Thurston map for simply or
doubly degenerate Kleinian groups corresponding to representations of
{\em closed} surface groups that have {\em
  amalgamation geometry}.

In Section 7, we extend these results to include surface groups with
punctures.

In Section 8, we weaken the hypothesis of {\em amalgamation geometry}
to what we have called {\em graph amalgamation geometry} and describe
the modifications necessary to extend our results to such geometries.

In Section 9, we weaken the hypothesis further to {\em split
  geometry} which allows for Margulis tubes to cut across the blocks.

In Section 10, we further generalise
these result to include 
hyperbolic manifolds whose cores are incompressible away from
cusps. (We had termed such manifolds {\em pared manifolds with
  incompressible boundary} in \cite{brahma-pared}.)

In Section 11, we give a scheme for proving that the Minsky model for
surface groups \cite{minsky-elc1} has split geometry. Details will
appear in the second part of this paper \cite{mahan-split}.

In Section 12, we propose an extension of the Sullivan-McMullen
dictionary between Kleinian groups and complex dynamics, and suggest
an analogue of Yoccoz puzzles in the 3 dimensional setting. 

\medskip

\noindent {\bf Acknowledgements:} Its a pleasure to thank Jeff Brock,
Dick Canary and Yair
Minsky
for their support, both personal and mathematical,
during the course of this work. In particular, the generalisations of
amalgamation geometry to graph amalgamation geometry and split
geometry
 were made to fill a gap in a 
previous
 version of this paper. The gap was brought to my notice by Minsky and
 Canary.   I, nevertheless, claim credit for
any errors and gaps that might still persist.

\section{Preliminaries and Amalgamation Geometry}

\subsection{Hyperbolic Metric Spaces}

We start off with some preliminaries about hyperbolic metric
spaces  in the sense
of Gromov \cite{gromov-hypgps}. For details, see \cite{CDP}, \cite{GhH}. Let $(X,d)$
be a hyperbolic metric space. The 
{\bf Gromov boundary} of 
 $X$, denoted by $\partial{X}$,
is the collection of equivalence classes of geodesic rays $r:[0,\infty)
\rightarrow\Gamma$ with $r(0)=x_0$ for some fixed ${x_0}\in{X}$,
where rays $r_1$
and $r_2$ are equivalent if $sup\{ d(r_1(t),r_2(t))\}<\infty$.
Let $\widehat{X}$=$X\cup\partial{X}$ denote the natural
 compactification of $X$ topologized the usual way(cf.\cite{GhH} pg. 124).

{\bf Definitions:}  A subset $Z$ of $X$ is said to be 
{\bf $k$-quasiconvex}
 if any geodesic joining points of  $ Z$ lies in a $k$-neighborhood of $Z$.
A subset $Z$ is {\bf quasiconvex} if it is $k$-quasiconvex for some
$k$. 
(For  simply connected real hyperbolic
manifolds this is equivalent to saying that the convex hull of the set
$Z$ lies in a  bounded neighborhood of $Z$. We shall have occasion to
use this alternate characterisation.)
A map $f$ from one metric space $(Y,{d_Y})$ into another metric space 
$(Z,{d_Z})$ is said to be
 a {\bf $(K,\epsilon)$-quasi-isometric embedding} if
 
\begin{center}
${\frac{1}{K}}({d_Y}({y_1},{y_2}))-\epsilon\leq{d_Z}(f({y_1}),f({y_2}))\leq{K}{d_Y}({y_1},{y_2})+\epsilon$
\end{center}
If  $f$ is a quasi-isometric embedding, 
 and every point of $Z$ lies at a uniformly bounded distance
from some $f(y)$ then $f$ is said to be a {\bf quasi-isometry}.
A $(K,{\epsilon})$-quasi-isometric embedding that is a quasi-isometry
will be called a $(K,{\epsilon})$-quasi-isometry.

A {\bf $(K,\epsilon)$-quasigeodesic}
 is a $(K,\epsilon)$-quasi-isometric embedding
of
a closed interval in $\Bbb{R}$. A $(K,K)$-quasigeodesic will also be called
a $K$-quasigeodesic.

Let $(X,{d_X})$ be a hyperbolic metric space and $Y$ be a subspace that is
hyperbolic with the inherited path metric $d_Y$.
By 
adjoining the Gromov boundaries $\partial{X}$ and $\partial{Y}$
 to $X$ and $Y$, one obtains their compactifications
$\widehat{X}$ and $\widehat{Y}$ respectively.

Let $ i :Y \rightarrow X$ denote inclusion.

{\bf Definition:}   Let $X$ and $Y$ be hyperbolic metric spaces and
$i : Y \rightarrow X$ be an embedding. 
 A {\bf Cannon-Thurston map} $\hat{i}$  from $\widehat{Y}$ to
 $\widehat{X}$ is a continuous extension of $i$.

The following  lemma (Lemma 2.1 of \cite{mitra-ct})
 says that a Cannon-Thurston map exists
if for all $M > 0$ and $y \in Y$, there exists $N > 0$ such that if $\lambda$
lies outside an $N$ ball around $y$ in $Y$ then
any geodesic in $X$ joining the end-points of $\lambda$ lies
outside the $M$ ball around $i(y)$ in $X$.
For convenience of use later on, we state this somewhat
differently.

\begin{lemma}
A Cannon-Thurston map from $\widehat{Y}$ to $\widehat{X}$
 exists if  the following condition is satisfied:

Given ${y_0}\in{Y}$, there exists a non-negative function  $M(N)$, such that 
 $M(N)\rightarrow\infty$ as $N\rightarrow\infty$ and for all geodesic segments
 $\lambda$  lying outside an $N$-ball
around ${y_0}\in{Y}$  any geodesic segment in $\Gamma_G$ joining
the end-points of $i(\lambda)$ lies outside the $M(N)$-ball around 
$i({y_0})\in{X}$.

\label{contlemma}
\end{lemma}

\smallskip

The above result can be interpreted as saying that a Cannon-Thurston map 
exists if the space of geodesic segments in $Y$ embeds properly in the
space of geodesic segments in $X$.

\subsection{Amalgamation Geometry}

We start with a hyperbolic surface $S$  without punctures. The
hyperbolic structure is arbitrary, but it is important that a choice
be made. 

\smallskip

{\bf The Amalgamated Building Block}

\smallskip

\noindent

For the construction of an amalgamated block $B$, $I$ will denote the closed
interval $[0,3]$. We will describe a geometry on $S \times I$. 
$B$ has a {\bf geometric core} $K$ with bounded geometry boundary
and a preferred geodesic $\gamma (= \gamma_B )$ of bounded length.

 There
will exist $\epsilon_0 , \epsilon_1 , D$ (independent of the block $B$)
such that the following hold:

\begin{enumerate}

\item $B$ is identified with $S \times I$\\
\item $B$ has a {\bf geometric core} $K$ identified with $S \times
  [1,2]$. ( $K$, in its intrinsic path metric, may be thought of, for
  convenience,  as a
  convex hyperbolic manifold with boundary consisting of pleated
  surfaces. But we will have occasion to use geometries that are only
  quasi-isometric to such geometries when lifted to universal
  covers. As of now, we do not impose any further restriction on the
  geometry of $K$. ) \\
\item $\gamma$ is homotopic to a simple closed curve on $S \times \{ i
  \} $ for any $ i\in I$ \\
\item $\gamma$ is small, i.e. the length of $\gamma$ is bounded above
  by $\epsilon_0$ \\
\item The intrinsic metric on $S \times i$ (for $i = 1,2$) has bounded
  geometry, i.e. any closed geodesic on $S \times \{ i \}$ has length
  bounded below by $\epsilon_1$. Further, the diameter of $S \times \{
  i \}$ is bounded above by $D$. (The latter restriction would have
  followed naturally had we assumed that  the curvature of $S \times
  \{ i \}$ is hyperbolic or at least pinched negative.) \\
\item There exists a regular neighborhood $N_k ( \gamma ) \subset K$
  of $\gamma$ which is homeomorphic to a solid torus, such that $N_k
  (\gamma ) \cap S \times \{ i \}$ is homeomorphic to an open annulus
  for $i = 1, 2$. We shall have occasion to denote $N_k (\gamma )$ by
  $T_\gamma$ and call it the Margulis tube corresponding to $\gamma$.  \\
\item $S \times [0,1]$ and $S \times [1,2]$ are given the product
  structures corresponding to the bounded geometry structures on $S
  \times \{ i \}$, for $i = 1,2 $ respectively. \\
\end{enumerate}

We next describe the geometry of the {\it geometric core}
$K$. $K - T_\gamma$ has one or two components according as
$\gamma$ does not or does separate $S$. These components shall be
called {\bf amalgamation components} of $K$. Let $K_1$ denote 
such an {\it amalgamation component. } Then a lift $\widetilde{K_1}$
of $K_1$ to $\widetilde{K}$ is bounded by lifts $\widetilde{T_\gamma
}$ of $T_\gamma$. The union of such a lift $\widetilde{K_1}$ along
with the lifts $\widetilde{T_\gamma }$ that bound it will be called an
{\bf amalgamation component} of $\widetilde{K}$. 

Note that two amalgamation components of $\widetilde{K}$, if they
intersect, shall do so along a lift $\widetilde{T_\gamma }$ of
$T_\gamma$. In this case, they shall be referred to as {\bf adjacent
  amalgamation components}. 

In addition to the above structure of $B$, we require in addition that
there exists $C > 0$ (independent of $B$) such that

\smallskip
\noindent
$\bullet$ Each amalgamation component  of $\widetilde{K}$ is
$C$-quasiconvex in the intrinsic metric on $\widetilde{K}$.

\smallskip

\noindent {\bf Note 1:} Quasiconvexity of an amalgamation component
follows from the 
fact that any geometric subgroup of infinite index in a surface group
is quasiconvex in the latter. The restriction above is therefore to
ensure uniform quasiconvexity. We shall strengthen this restriction
further when we describe the geometry of $\widetilde{M}$, where $M$ is
a 3-manifold built up of blocks of {\it amalgamation geometry} and
those of bounded geometry by gluing them end to end. We shall require
that each amalgamation component is {\bf uniformly quasiconvex in
  $\widetilde{M}$} rather than just in $\widetilde{K}$. \\
\noindent {\bf Note 2:} 
So far, the restrictions on $K$ are quite mild. There are really two
restrictions. One is the existence of a bounded length simple closed
geodesic whose regular neighborhood intersects the bounding surfaces
of $K$ in annulii. The second restriction is that the two boundary
surfaces of $K$ have {\it bounded geometry}. 

\smallskip

The copy of $S \times I$ thus obtained, with the restrictions above,
will be called a {\bf building block of amalgamated geometry} or an
{\bf amalgamation geometry building block}, or simply an {\bf
  amalgamation block}.

\smallskip

{\bf Thick Block}

Fix constants $D, \epsilon$ and let $\mu = [p,q]$ be  an $\epsilon$-thick
Teichmuller geodesic of length less than $D$. $\mu$ is
$\epsilon$-thick means that for any $x \in \mu$ and any closed
geodesic $\eta$ in the hyperbolic
surface $S_x$ over $x$, the length of $\eta$ is greater than
$\epsilon$.
Now let $B$ denote the universal curve over $\mu$ reparametrized
such that the length of $\mu$ is covered in unit time. 
 Thus $B = S \times [0,1]$
topologically.

 $B$ is given the path
metric and is 
called a {\bf thick building block}.

Note that after acting by an element of the mapping class group, we
might as well assume that $\mu$ lies in some given compact region of
Teichmuller space. This is because the marking on $S \times \{ 0 \}$
is not important, but rather its position relative to $S \times \{ 1 \}$
Further, since we shall be constructing models only upto
quasi-isometry, we might as well assume that $S \times \{ 0 \}$ and $S
\times \{ 1 \}$ {\em lie in the orbit} under the mapping class group
of some fixed base surface. Hence $\mu$ can be further simplified to
be a Teichmuller geodesic joining a pair $(p, q)$ amongst a finite set of
points in the orbit of a fixed hyperbolic surface $S$.

\smallskip 
 
{\bf The Model Manifold}

\smallskip

\noindent Note that the boundary of an amalgamation block $B_i$ consists of $S \times \{
0,3 \}$ and the intrinsic path metric on each such $S \times \{ 0 \}$
or
 $S \times \{ 3 \}$ is of bounded geometry.
Also,  the boundary of a thick block $B$ consists of $S \times \{
0,1 \}$, where $S_0, S_1$ lie in some given bounded region of
Teichmuller space. The intrinsic path metrics on each such $S \times \{ 0 \}$
or
 $S \times \{ 1 \}$ is  the path metric on $S$. 

The model
manifold of {\bf amalgamation geometry} is obtained from  $S \times J$
(where $J$ is a sub-interval of $\Bbb{R}$, 
which may be  semi-infinite or bi-infinite. In the former case, we
choose the usual normalisation $J = [ 0, {\infty })$ ) by first choosing a
  sequence
 of blocks $B_i$ (thick or amalgamated) and corresponding intervals $I_i =
[0,1]$ or $ [0,3]$  according as $B_i$ is thick or amalgamated. The
  metric on $S \times I_i$ is then declared to be that on the  building block
  $B_i$. Implicitly, we are requiring that the surfaces along which
  gluing occurs have the same metric. Thus we have,

\smallskip

\noindent {\bf Definition:} A manifold $M$ homeormorphic to $S \times
  J$, where $ J = [0,
  {\infty })$ or $J = ( - \infty , \infty )$, is said to be a model of {\bf
  amalgamation geometry} if  \\

\begin{enumerate}

\item there is a fiber preserving homeomorphism from $M$ to
  $\widetilde{S} \times J$ 
that lifts to  a quasi-isometry of universal covers \\ 

\item there exists a sequence $I_i$ of intervals (with disjoint
 interiors)
 and blocks $B_i$
where the metric on $S \times I_i$ is the same as
 that on some  building block $B_i$ \\

\item $\bigcup_i I_i = J$ \\

\item There exists $ C > 0$ such that for all amalgamated blocks $B$
  and geometric cores $K \subset B$, all amalgamation components of
  $\widetilde{K}$ are $C$-quasiconvex in $\widetilde{M}$

\end{enumerate}

{\bf Note:} The last restriction (4) above is a global restriction
on the geometry of amalgamation components, not just a local one
(i.e. quasiconvexity in $\widetilde{M}$ rather than $\widetilde{B}$ is
required.)

\smallskip

The figure below illustrates schematically what the model looks
like. Filled squares correspond to solid torii along which amalgamation
occurs. The adjoining piece(s) denote amalgamation blocks of $K$.
 The blocks which have no filled squares are the
{\it thick blocks} and those with filled squares are the {\it amalgamated
  blocks}

\smallskip

\begin{center}

\includegraphics[height=4cm]{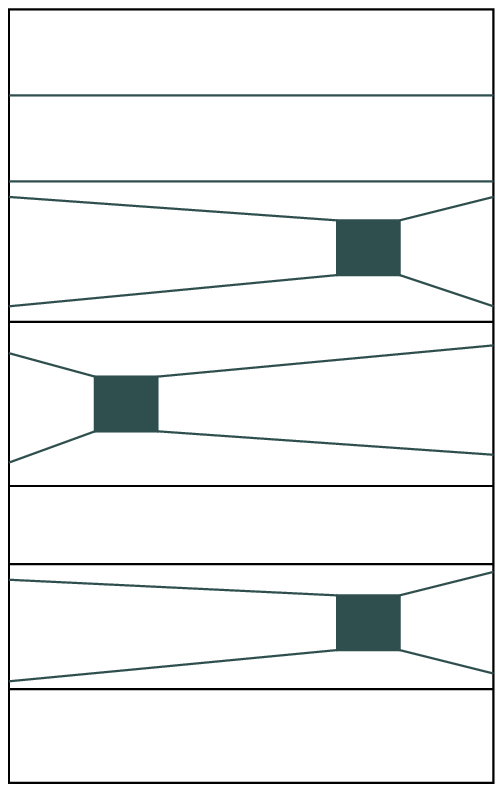}

\smallskip

\underline{Figure 1: {\it Model of amalgamated geometry (schematic)} }

\end{center}

\smallskip

\noindent {\bf Definition:} A manifold $M$ homeomorphic to $S \times
  J$, where $ J = [0,
  {\infty })$ or $J = ( - \infty , \infty )$, is said to have {\bf
  amalgamated geometry} if there exists $K, \epsilon > 0$  and a model
  manifold $M_1$ of {\em amalgamation geometry} such that \\
\begin{enumerate}
\item there exists
  a homeomorphism $\phi$ from $M$ to $M_1$. This induces from the block
  decomposition of $M_1$ a block decomposition of $M$. \\
\item We require in addition that  the induced homeomorphism
  $\tilde{\phi}$ between universal covers of blocks is a $(K,
  \epsilon)$
 quasi-isometry.
\end{enumerate}

We shall usually suppress the homeomorphism $\phi$ and take $M$ itself
to be a model manifold of {\em amalgamation geometry}.

A geometrically tame hyperbolic 3-manifold is said to have
{\bf amalgamated geometry} if each end has amalgamated geometry.

\smallskip

{\bf Note:} We shall later have occasion to introduce a different
model, called the {\bf graph model}

\section{Relative Hyperbolicity}

In this section, we shall recall first certain notions of relative
hyperbolicity due to Farb \cite{farb-relhyp}, Klarreich
\cite{klarreich} and the author \cite{brahma-ibdd}. Using these, we shall
derive certain Lemmas that will be useful in studying the geometry of
the universal covers of building blocks.

\subsection{Electric Geometry}

We start with a surface $S$ (assumed hyperbolic for the time being) of $(K,
\epsilon )$
bounded
geometry, i.e. $S$ has diameter bounded by $K$ and injectivity radius
bounded below by $\epsilon$. Let $\sigma$ be a simple closed geodesic
on $S$. 
Replace $\sigma$ by a copy of $\sigma \times [0,1]$, by cutting open along
$\sigma$ and gluing in a copy of $\sigma \times [0,1] = A_\sigma$. 
(This is like `grafting' but we shall not have much use for this
similarity in this paper.) Let $S_G$ denote the grafted surface. $S_G
- A_\sigma$ has one or two components according as $\sigma$ does
not or does separate $S$. Call these {\bf amalgamation component(s)}
of $S$ We shall denote amalgamation components as $S_A$. We construct
a pseudometric  on $S_G$, by declaring the metric 
on each amalgamation component to be zero and to be the product metric on
$A_\sigma$. Thus we define: 

\smallskip

\noindent 
$\bullet$ the length of any path that lies in the interior of
 an amalgamation component to be zero \\
$\bullet$  the length of any path  that lies in $A_\sigma$
to be  its (Euclidean) length in the path metric on $A_\sigma$ \\
$\bullet$ the length of any other path to be the sum of lengths of
pieces of the above two kinds. \\

\smallskip

This allows us to define distances by taking the infimum of lengths of
  paths joining pairs of points and gives us a path pseudometric,
  which we call the {\bf electric metric} on $S_G$. The electric metric  also
  allows us to define
geodesics. Let us
call $S_G$ equipped with the  above pseudometric $(S_{Gel} , d_{Gel})$ (to be
  distinguished from a `dual' construction of an electric metric
  $S_{el}$ used in \cite{brahma-ibdd}, where  the geodesic $\sigma$,
  rather than its complementary component(s) is electrocuted.)

\smallskip

{\bf Important Note:} We may and shall regard $S$ as a graph of groups with
vertex group(s) the subgroup(s) corresponding to amalgamation
component(s) and edge group $Z$, the fundamental group of
$A_\sigma$. Then $\widetilde{S}$ equipped  with the lift of the above
pseudometric is quasi-isometric to  the tree corresponding to the
splitting
on which $\pi_1 (S)$ acts.

\smallskip

We shall be interested in the universal cover $\widetilde{S_{Gel}}$ of
$S_{Gel}$. Paths in $S_{Gel}$ and $\widetilde{S_{Gel}}$ will be called
electric paths 
(following Farb \cite{farb-relhyp}). Geodesics and quasigeodesics in 
the electric metric will be called electric geodesics and electric
quasigeodesics respectively.

\smallskip

{\bf Definitions:} \\
$\bullet$
 A path $\gamma : I \rightarrow Y$ in a path metric space $Y$ is a
K-quasigeodesic if we have
\begin{center}
$L({\beta}) \leq K L(A) + K$
\end{center}
for any subsegment $\beta = \gamma |[a,b]$ and any rectifiable
path $A : [a,b] \rightarrow Y$ with the
same endpoints. \\
$\bullet$  $\gamma$ is said to be an electric
$K, \epsilon$-quasigeodesic in $\widetilde{S_{Gel}}$
{ \bf without backtracking } if
 $\gamma$ is an electric $K$-quasigeodesic in $\widetilde{S_{Gel}}$ and
 $\gamma$ does not return to any
any lift  $\widetilde{S_A} \subset
\widetilde{S_{Gel}}$ (of an amalgamation component
$S_A \subset S$) after leaving it.\\

\smallskip

We collect together certain facts about the electric metric that Farb
proves in \cite{farb-relhyp}. $N_R(Z)$ will denote the
$R$-neighborhood about the subset $Z$ in the hyperbolic metric.
 $N_R^e(Z)$ will denote the
$R$-neighborhood about the subset $Z$ in the electric metric.

\begin{lemma} (Lemma 4.5 and Proposition 4.6 of \cite{farb-relhyp})
\begin{enumerate}
\item {\it Electric quasi-geodesics electrically track hyperbolic
  geodesics:} Given $P > 0$, there exists $K > 0$ with the following
  property: For some $\widetilde{S_{Gel}}$, 
let $\beta$ be any electric $P$-quasigeodesic without backtracking
from $x$ to
  $y$, and let $\gamma$ be the hyperbolic geodesic from $x$ to $y$. 
Then $\beta \subset N_K^e ( \gamma )$. \\
\item {\it Hyperbolicity:} There exists $\delta$ such that each
  $\widetilde{S_{Gel}}$  is $\delta$-hyperbolic, independent of the
  curve $\sigma$ whose lifts are electrocuted. \\
\end{enumerate}
\label{farb1}
\end{lemma}

{\bf Note:} As pointed out before, $S_{Gel}$ is quasi-isometric to a
tree and is therefore hyperbolic. The above assertion holds in far
greater generality than stated. We discuss this below.

\smallskip

We consider
 a hyperbolic metric space $X$ and a collection $\mathcal{H}$
of
{\em (uniformly) $ C$-quasiconvex uniformly separated subsets}, i.e.
there exists $D > 0$ such that for $H_1, H_2 \in \mathcal{H}$, $d_X (H_1,
H_2) \geq D$. In this situation $X$ is hyperbolic relative to the
collection $\mathcal{H}$. The result in this form is due to Klarreich
\cite{klarreich}. We give the general version of Farb's theorem below and
refer to \cite{farb-relhyp} and Klarreich \cite{klarreich} for proofs.

\begin{lemma} (See Lemma 4.5 and Proposition 4.6 of \cite{farb-relhyp}
 and Theorem 5.3 of Klarreich \cite{klarreich})
Given $\delta , C, D$ there exists $\Delta$ such that
if $X$ is a $\delta$-hyperbolic metric space with a collection
$\mathcal{H}$ of $C$-quasiconvex $D$-separated sets.
then,

\begin{enumerate}
\item {\it Electric quasi-geodesics electrically track hyperbolic
  geodesics:} Given $P > 0$, there exists $K > 0$ with the following
  property: Let $\beta$ be any electric $P$-quasigeodesic from $x$ to
  $y$, and let $\gamma$ be the hyperbolic geodesic from $x$ to $y$. 
Then $\beta \subset N_K^e ( \gamma )$. \\
\item $\gamma$ lies in a {\em hyperbolic} $K$-neighborhood of $N_0 ( \beta
  )$, where $N_0 ( \beta )$ denotes the zero neighborhood of $\beta$
  in the {\em electric metric}. \\
\item {\it Hyperbolicity:} 
  $X$  is $\Delta$-hyperbolic. \\
\end{enumerate}
\label{farb1A}
\end{lemma}

A special kind of {\em geodesic without backtracking} will be
necessary for universal covers $\widetilde{S_{Gel}}$
of surfaces with some electric metric. Let $\sigma$, $A_\sigma$ be as before.

\smallskip

Let $\lambda_e$ be an electric geodesic in some 
$(\widetilde{S_{Gel}},d_{Gel})$. 
Then, each segment of $\lambda_e$ between two lifts
$\widetilde{A_\sigma}$ of $A_\sigma$ (i.e. lying inside a lift of an
amalgamation component) is required to be perpendicular
to the bounding geodesics. We shall refer to these segments of $\lambda_e$ as
{\bf amalgamation segments} because they lie inside lifts of
the amalgamation components.  

Let $a , b$ be the points at which
$\lambda_e$ enters and leaves a lift $\widetilde{A_\sigma}$ of
$A_\sigma$. If $a, b$ lie on the same side, i.e. on a lift of either
$\sigma \times \{ 0 \}$ or $\sigma \times \{ 1 \}$, then we join $a,
b$ by the geodesic joining them. If they lie on opposite sides of 
 $\widetilde{A_\sigma}$, then 
 assume, for convenience,  that $a$ lies  on a lift of $\sigma \times
 \{ 0 \}$ and $b$ lies on a lift of $\sigma \times
 \{ 1 \}$. Then we join $a$ to $b$ by a union of 2 geodesic segments
 $[a,c]$ and $[d,b]$
 lying along $\widetilde{\sigma} \times \{ 0 \}$ and
 $\widetilde{\sigma} \times \{ 1 \}$ respectively (for some lift
 $\widetilde{A_\sigma}$), along with a `horizontal' segment $[c,d]$,
 where $[c,d] \subset \widetilde{A_\sigma}$ projects to a segment of
 the form $\{ x \} \times [0,1] \subset \sigma \times [0,1]$. We
 further require that the sum of the lengths $d(a,c)$ and $d(d,b)$
 is the minimum possible. The union of the  three segments $[a,c],
 [c,d], [d,b]$ shall be denoted by  $[a,b]_{int}$ and shall be referred to
as an {\bf interpolating segment}.See figure below.

\smallskip

\begin{center}

\smallskip

\includegraphics[height=4cm]{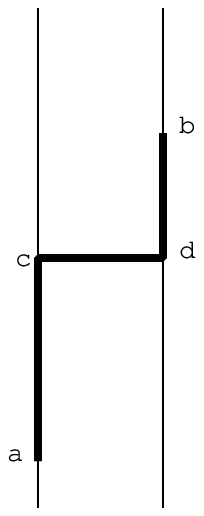}

\underline{Figure 2: {\it Interpolating segment} }

\end{center}

\smallskip

 The union of the {\em amalgamation
  segments}
along with the {\em interpolating segments} gives rise to a preferred
representative of a quasigeodesic without backtracking
 joining the end-points of $\lambda_{Gel}$.
 Such a representative of the class of $\lambda_{Gel}$
shall be called the {\bf canonical representative} of
$\lambda_{Gel}$. Further, the underlying set of the canonical
representative in the {\em hyperbolic metric} shall be called the {\bf
  electro-ambient representative} $\lambda_q$ of $\lambda_e$. 
Since $\lambda_q$  turns out to be a hyperbolic quasigeodesic (Lemma
  \ref{Gea} below), we
shall
also call it an {\bf electro-ambient quasigeodesic}. See Figure 3
below:

\smallskip

\begin{center}

\includegraphics[height=4cm]{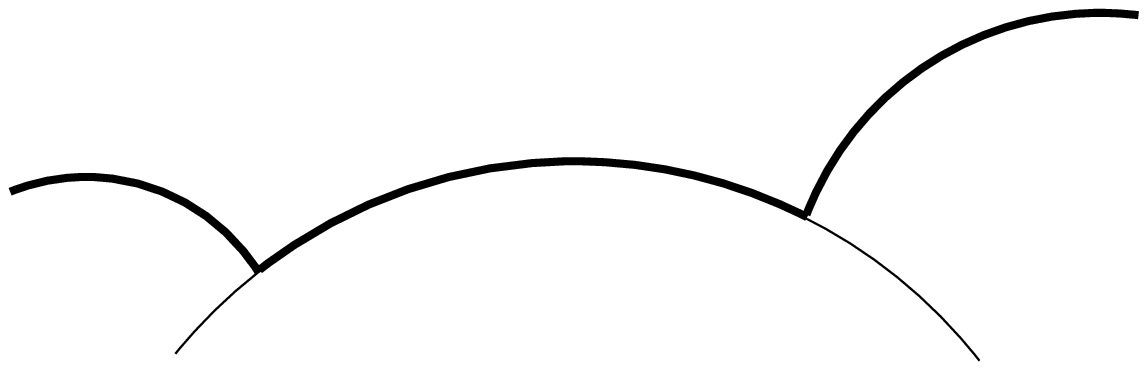}

\smallskip

\underline{Figure 3:{\it Electro-ambient quasigeodesic} }

\end{center}

\smallskip

{\bf Remark:}
We note first that if we collapse each lift of $A_\sigma$ along the $I
( = [0,1])$-fibres, (and thus obtain a geodesic that is a lift of
$\sigma$), then $\lambda_{Gel}$ becomes an electric geodesic
$\lambda_{el}$ in the universal cover  $\widetilde{S_{el}}$ of
$S_{el}$. Here $S_{el}$ denotes the space obtained by electrocuting
the geodesic $\sigma$ (See Section 3.1 of \cite{brahma-ibdd}. 

\smallskip

Let $c : S_G \rightarrow S$ be the map that {\it collapses} $I$-fibres,
i.e. it maps the annulus $A_\sigma = \sigma \times I$ to the geodesic
$\sigma$ by taking $(x,t)$ to $x$. The lift $\tilde{c} :
\widetilde{S_G} \rightarrow \widetilde{S}$ collapses
each lift of $A_\sigma$ along the $I
( = [0,1])$-fibres to  a geodesic that is a lift of
$\sigma$). Also it takes  $\lambda_{Gel}$ to an electric geodesic
$\lambda_{el}$ in the universal cover  $\widetilde{S_{el}}$ of
$S_{el}$ (that $\lambda_{el}$ is an electric geodesic in
$\widetilde{S_{el}}$ follows easily, say from normal forms). These
were precisely the electro-ambient quasigeodesics in the space
$\widetilde{S_el}$ (See Section 3.1 of \cite{brahma-ibdd} for definitions). 

\smallskip

{\bf Remark:} The electro-ambient geodesics in the sense of
\cite{brahma-ibdd} and those in the present paper differ slightly. The
difference is due to the grafting annulus $A_\sigma$ that we use here
in place of $\sigma$. What is interesting is that whether we
electrocute $\sigma$ (to obtain $S_{el}$) or its complementary
components  (to obtain $S_{Gel}$), we obtain very nearly the same
electro-ambient geodesics. In fact modulo $c$, they {\it are} the
same. 

\smallskip

We now recall a Lemma from \cite{brahma-ibdd}:

\begin{lemma}
(See Lemma 3.7 of \cite{brahma-ibdd} )
There exists $(K, \epsilon)$ such that each electro-ambient
representative $\lambda_{el}$ of an electric geodesic in
$\widetilde{S_{el}}$
is a $(K,
\epsilon)$ hyperbolic quasigeodesic.
\label{ea}
\end{lemma}

Since $\tilde{c}$ is clearly a quasi-isometry, it follows easily that:

\begin{lemma}
There exists $(K, \epsilon)$ such that each electro-ambient
representative $\lambda_{Gel}$ of an electric geodesic in
$\widetilde{S_{Gel}}$
is a $(K,
\epsilon)$ hyperbolic quasigeodesic.
\label{Gea}
\end{lemma}
 
In the above form, {\em electro-ambient quasigeodesics} are considered
only 
in the context of surfaces, closed geodesics on them and their
complementary ({\it amalgamation}) components. A considerable
generalisation of this was obtained in \cite{brahma-ibdd}, which will
be necessary while considering the global geometry of $\widetilde{M}$
(rather than the geometry of $\widetilde{B}$, for an amalgamated
building block $B$).

We recall a  definition from \cite{brahma-ibdd}:

\smallskip

\noindent {\bf Definitions:} Given a collection $\mathcal{H}$
of $C$-quasiconvex, $D$-separated sets and a number $\epsilon$ we
shall say that a geodesic (resp. quasigeodesic) $\gamma$ is a geodesic
(resp. quasigeodesic) {\bf without backtracking} with respect to
$\epsilon$ neighborhoods if $\gamma$ does not return to $N_\epsilon
(H)$ after leaving it, for any $H \in \mathcal{H}$. 
A geodesic (resp. quasigeodesic) $\gamma$ is a geodesic
(resp. quasigeodesic) {\bf without backtracking} if it is a geodesic
(resp. quasigeodesic) without backtracking with respect to
$\epsilon$ neighborhoods for some $\epsilon \geq 0$.

\smallskip

{\bf Note:} For strictly convex sets, $\epsilon = 0$ suffices, whereas
for convex sets any $\epsilon > 0$ is enough.

\smallskip

Let $X$ be a $\delta$-hyperbolic metric
space, and $\mathcal{H}$ a family of $C$-quasiconvex, $D$-separated,
 collection of subsets. Then by Lemma \ref{farb1A},
$X_{el}$ obtained by electrocuting the subsets in $\mathcal{H}$ is
a $\Delta = \Delta ( \delta , C, D)$ -hyperbolic metric space. Now,
let $\alpha = [a,b]$ be a hyperbolic geodesic in $X$ and $\beta $ be
an electric 
$P$-quasigeodesic without backtracking
 joining $a, b$. Replace each maximal subsegment, (with end-points $p,
 q$, say)
starting from the left
 of
$\beta$ lying within some $H \in \mathcal{H}$
by a hyperbolic  geodesic $[p,q]$. The resulting
{\bf connected}
path $\beta_q$ is called an {\em electro-ambient representative} in
$X$. 

\smallskip

In  \cite{brahma-ibdd} we noted that $\beta_q$ {\it need
  not be a hyperbolic quasigeodesic}. However, we did adapt Proposition
  4.3 of Klarreich \cite{klarreich} to obtain the following:

\begin{lemma} (See Proposition 4.3 of \cite{klarreich}, also see Lemma
  3.10 of \cite{brahma-ibdd}) 
Given $\delta$, $C, D, P$ there exists $C_3$ such that the following
holds: \\
Let $(X,d)$ be a $\delta$-hyperbolic metric space and $\mathcal{H}$ a
family of $C$-quasiconvex, $D$-separated collection of quasiconvex
subsets. Let $(X,d_e)$ denote the electric space obtained by
electrocuting elements of $\mathcal{H}$.  Then, if $\alpha , \beta_q$
denote respectively a hyperbolic geodesic and an electro-ambient
$P$-quasigeodesic with the same end-points, then $\alpha$ lies in a
(hyperbolic) 
$C_3$ neighborhood of $\beta_q$.
\label{ea-strong}
\end{lemma}

{\bf Note:} The above Lemma will be needed while considering geodesics
in $\widetilde{M}$.

\subsection{Electric isometries}

Recall that $S_G$ is a grafted surface  obtained
from a (fixed) hyperbolic metric by grafting an annulus $A_\sigma$
in place of a geodesic
$\sigma$. 

Now let $\phi$ be any diffeomorphism of $S_G$ that fixes $A_\sigma$
pointwise and
(in case $( S_G - A_\sigma )$ has two components) preserves each
amalgamation component as a set, i.e. $\phi$ sends each amalgamation
component to itself. Such a $\phi$ will be called a {\bf component
  preserving diffeomorphism}. Then in the electrocuted surface
$S_{Gel}$, any electric geodesic has length equal to the number of
times it crosses $A_\sigma$. It follows that $\phi$ is an isometry of
$S_{Gel}$. (See Lemma 3.12 of \cite{brahma-ibdd} for an analogous
result in $S_{el}$.) We state this below.

\begin{lemma}
Let $\phi$ denote a component preserving diffeomorphism of $S_G$. 
 Then $\phi $ induces an isometry of $(S_{Gel},d_{Gel})$.
\label{phi-isom1}
\end{lemma}

Everything in the above can be lifted to the universal cover
$\widetilde{S_{Gel}}$. We  let $\widetilde{\phi}$ denote the lift of $\phi$ to
$\widetilde{S_{Gel}}$.
This gives \\

\begin{lemma}
Let $\widetilde{\phi}$ denote a lift of a component preserving
 diffeomorphism $\phi$ to 
$(\widetilde{S_{Gel}},d_{Gel})$. 
 Then $\widetilde{\phi} $ induces an isometry of
$ ( \widetilde{S_{Gel}},d_{Gel})$.  
\label{phi-isom2}
\end{lemma}

\subsection{Nearest-point Projections}

We need the following basic lemmas from \cite{mitra-trees} and
\cite{brahma-ibdd}. 

The following Lemma  says nearest point projections in a $\delta$-hyperbolic
metric space do not increase distances much.

\begin{lemma}
(Lemma 3.1 of \cite{mitra-trees})
Let $(Y,d)$ be a $\delta$-hyperbolic metric space
 and  let $\mu\subset{Y}$ be a $C$-quasiconvex subset, 
e.g. a geodesic segment.
Let ${\pi}:Y\rightarrow\mu$ map $y\in{Y}$ to a point on
$\mu$ nearest to $y$. Then $d{(\pi{(x)},\pi{(y)})}\leq{C_3}d{(x,y)}$ for
all $x, y\in{Y}$ where $C_3$ depends only on $\delta, C$.
\label{easyprojnlemma}
\end{lemma}

The next lemma  says that quasi-isometries and nearest-point projections on
hyperbolic metric spaces `almost commute'.

\begin{lemma}
(Lemma 3.5 of \cite{mitra-trees})Suppose $(Y_1,d_1)$ and $(Y_2,d_2)$
are $\delta$-hyperbolic.
Let $\mu_1$ be some geodesic segment in $Y_1$ joining $a, b$ and let $p$
be any vertex of $Y_1$. Also let $q$ be a vertex on $\mu_1$ such that
${d_1}(p,q)\leq{d_2}(p,x)$ for $x\in\mu_1$. 
Let $\phi$ be a $(K,{\epsilon})$ - quasiisometric embedding
 from $Y_1$ to $Y_2$.
Let $\mu_2$ be a geodesic segment 
in $Y_2$ joining ${\phi}(a)$ to ${\phi}(b)$ . Let
$r$ be a point on $\mu_2$ such that ${d_2}({\phi}(p),r)\leq{d_2}({\phi(p)},x)$ for $x\in\mu_2$.
Then ${d_2}(r,{\phi}(q))\leq{C_4}$ for some constant $C_4$ depending   only on
$K, \epsilon $ and $\delta$. 
\label{almost-commute}
\end{lemma}

For our purposes we shall need the above Lemma for quasi-isometries
from $\widetilde{S_a}$ to  $\widetilde{S_b}$ for two different
hyperbolic structures on the same surface. We shall also need it for
 electrocuted surfaces.

\smallskip

Yet another property that we shall require for nearest point
projections is that nearest point projections in the electric metric
and in the `almost hyperbolic' metric (coming as a lift of the metric
on $S_G$) almost agree. Let $\widetilde{S_G} = Y$ be
the universal cover of a surface with the grafted metric. Equip $Y$
with the path metric $d$ 
as usual. Then
$Y$ is  quasi-isometric to the hyperbolic plane. 
Recall that $d_{Gel}$ denotes the electric
metric on $Y$ obtained by electrocuting the lifts of complementary
components. Now, 
let $\mu = [a,b]$ be a  geodesic on $(Y,d)$ and let
$\mu_q$ denote the electro-ambient quasigeodesic 
joining $a, b$ (See Lemma \ref{ea}). Let $\pi$ denote the nearest point projection in
$(Y,d)$. Tentatively, let $\pi_e$ denote the nearest point projection in
$(Y,d_{Gel})$. Note that $\pi_e$ is not well-defined. It is defined upto a
bounded amount of discrepancy in the electric metric $d_e$. But we
would like to make $\pi_e$ well-defined upto a bounded amount of
discrepancy in the  metric $d$.

\smallskip

{\bf Definition:} Let $y \in Y$ and let $\mu_q$ be an electro-ambient
representative of an electric geodesic $\mu_{Gel}$ in
$(Y,d_{Gel})$. Then $\pi_e(y) = z \in \mu_q$ if the ordered pair $\{
d_{Gel}(y,\pi_e(y)), d(y, \pi_e(y) ) \}$ is minimised at $z$. 

The proof of the following Lemma shows
 that this gives us a definition of $\pi_e$ which is ambiguous by
a finite amount of discrepancy not only in the electric metric but
also
in the hyperbolic metric.

\begin{lemma} There exists $C > 0$ such that the following holds.
Let $\mu$ be a hyperbolic geodesic joining $a, b$. Let
  $\mu_{Gel}$ be an electric geodesic
  joining $a, b$. Also let $\mu_q$ be the electro-ambient
  representative of $\mu_{Gel}$. Let $\pi_h$ denote the nearest point
  projection of $Y$ onto $\mu$. 
$d(\pi_h(y) , \pi_e(y))$ is uniformly bounded.
\label{hyp=elproj}
\end{lemma}

{\bf Proof:} 
This Lemma is similar to   Lemma
3.16 of \cite{brahma-ibdd}, but its proof is somewhat different. For
the purposes of this lemma we shall refer to the metric on
$\widetilde{S_G}$ as the hyperbolic metric whereas it is in fact only
quasi-isometric to it.

$[u, v]$ and $[u,v]_q$ will denote respectively the hyperbolic
geodesic and the electro-ambient quasigeodesic
joining $u, v$. Since $[u,v]_q$ is a quasigeodesic by Lemma \ref{ea},
it suffices to 
show that for any $y$, its hyperbolic and electric  projections $\pi_h
(y), \pi_e (y)$ almost
agree. 

First note that  any  hyperbolic geodesic $\eta$ in $\widetilde{S_G}$
is also an 
electric geodesic. This follows from the fact that $(\widetilde{S_G} ,
d_{Gel} )$ maps to a tree $T$ (arising from the splitting along $\sigma$)
with the pullback of every vertex a set of diameter zero in the
pseudometric $d_{Gel}$. Now if a path in $\widetilde{S_G}$ projects to
a path in $T$ that is not a geodesic, then it must backtrack. Hence,
it must leave an amalgamating component and return to it. Such a path
can clearly not be a hyperbolic geodesic in $\widetilde{S_G}$ (since
each amalgamating component is convex). 

Next, it follows that  hyperbolic projections
automatically minimise electric distances. Else as in the preceding
paragraph, $[y,\pi_h (y)]$ would have to cut a lift of
$\widetilde{\sigma} = \widetilde{\sigma_1}$
 that separates $[u,v]_q$. Further,  $[y,\pi_h (y)]$ cannot return to
 $\widetilde{\sigma_1 }$ after leaving it. 
Let $z$ be the first point at which
 $[y,\pi_h (y)]$ meets $\widetilde{\sigma_1}$. Also let $w$ be the
point on  $[u,v]_q \cap \widetilde{\sigma_1}$ that is nearest to
$z$. Since  amalgamation segments of $[u,v]_q$ meeting
$\widetilde{\sigma_1}$ are perpendicular to the latter, it follows
that $d(w,z) < d(w,\pi_h (y) )$ and therefore  $d(y,z) < d(y,\pi_h (y)
)$ contradicting the definition of $\pi_h (y)$. Hence hyperbolic projections
automatically minimise electric distances.

Further, it follows by repeating the argument in the first paragraph
 that $[y,\pi_h (y)]$ and $[y, \pi_e (y)]$ pass through the same set
 of amalgamation components in the same order; in particular they
 cut across the same set of lifts of $\widetilde{\sigma}$. Let
 $\widetilde{\sigma_2}$ be the last such lift. Then
 $\widetilde{\sigma_2}$ forms the boundary of an amalgamation
 component $\tilde{S_A}$ whose intersection with $[u,v]_q$ is of the
 form $[a,b] \cup [b,c] \cup [c,d]$, where $[a,b] \subset
 \widetilde{\sigma_3}$ and $[c,d] \subset \widetilde{\sigma_4}$ are
 subsegments of two lifts of $\sigma$ and $[b,c]$ is perpendicular to
 these two. Then the nearest-point 
projection of $\widetilde{\sigma_2}$ onto each of $[a,b], [b,c],
 [c,d]$ has uniformly bounded diameter. Hence the nearest point
 projection of $\tilde{\sigma_2}$ onto the hyperbolic geodesic $[a,d]
 \subset \widetilde{S_A}$ has uniformly bounded diameter.
 The result follows. $\Box$

\subsection{Coboundedness and Consequences}

In this Section, we collect together a few more results that
strengthen Lemmas \ref{farb1} and  \ref{farb1A}.

{\bf Definition:} A collection $\mathcal{H}$ of uniformly
$C$-quasiconvex sets in a $\delta$-hyperbolic metric space $X$
is said to be {\bf mutually D-cobounded} if 
 for all $H_i, H_j \in \mathcal{H}$, $\pi_i
(H_j)$ has diameter less than $D$, where $\pi_i$ denotes a nearest
point projection of $X$ onto $H_i$. A collection is {\bf mutually
  cobounded } if it is mutually D-cobounded for some $D$. 

\smallskip

\begin{lemma} 
Suppose $X$ is a $\delta$-hyperbolic metric space with a collection
$\mathcal{H}$ of $C$-quasiconvex $K$-separated $D$-mutually cobounded
subsets. There exists $\epsilon_0 = \epsilon_0 (C, K, D, \delta )$ such that
the following holds:

Let $\beta$ 
  be an electric $P$-quasigeodesic without backtracking
and $\gamma$ a hyperbolic geodesic,
  both joining $x, y$. Then, given $\epsilon \geq \epsilon_0$
 there exists $D = D(P, \epsilon )$ such that \\
\begin{enumerate}
\item {\it Similar Intersection Patterns 1:}  if
  precisely one of $\{ \beta , \gamma \}$ meets an
  $\epsilon$-neighborhood $N_\epsilon (H_1)$
of an electrocuted quasiconvex set
  $H_1 \in \mathcal{H}$, then the length (measured in the intrinsic path-metric
  on  $N_\epsilon (H_1)$ ) from the entry point
  to the 
  exit point is at most $D$. \\
\item {\it Similar Intersection Patterns 2:}  if
 both $\{ \beta , \gamma \}$ meet some  $N_\epsilon (H_1)$
 then the length (measured in the intrinsic path-metric
  on  $N_\epsilon (H_1)$ ) from the entry point of
 $\beta$ to that of $\gamma$ is at most $D$; similarly for exit points. \\
\end{enumerate}
\label{farb2A}
\end{lemma}
 
\noindent Summarizing, we have: \\
$\bullet$ If $X$ is a hyperbolic metric space and
$\mathcal{H}$ a collection of uniformly quasiconvex mutually cobounded
separated subsets,
then $X$ is hyperbolic relative to the collection $\mathcal{H}$ and
satisfies {\em Bounded Penetration}, i.e. hyperbolic geodesics and
electric quasigeodesics have similar intersection patterns in the
sense of Lemma \ref{farb2A}. \\

The relevance of co-boundedness comes from the following Lemma which
is essentially due to Farb \cite{farb-relhyp}.

\begin{lemma}
Let $M^h$ be a hyperbolic manifold of {\em i-bounded geometry}, with
Margulis tubes $T_i \in \mathcal{T}$ and horoballs $H_j \in
\mathcal{H}$. Then the lifts $\widetilde{T_i}$ and $\widetilde{H_j}$
are mutually co-bounded.
\label{coboundedHor&T}
\end{lemma}

The proof given in \cite{farb-relhyp} is for a collection of separated
horospheres, but the same proof works for neighborhoods of geodesics
and horospheres as well.

A closely related  theorem was proved by  McMullen
(Theorem 8.1 of \cite{ctm-locconn}).

As usual, $N_R (Z)$ will denote the $R$-neighborhood of the set $Z$. \\
Let $\cal{H}$ be a locally finite collection of horoballs in a convex
subset $X$ of ${\Bbb{H}}^n$ 
(where the intersection of a horoball, which meets $\partial X$ in a point, 
 with $X$ is
called a horoball in $X$).

\smallskip

{\bf Definition:} The $\epsilon$-neighborhood of a bi-infinite
geodesic in ${\Bbb{H}}^n$ will be called a {\bf thickened geodesic}. 

\smallskip

\begin{theorem} \cite{ctm-locconn}
Let $\gamma: I \rightarrow X \setminus \bigcup \cal{H}$ be an ambient
$K$-quasigeodesic (for $X$ a convex subset  of ${\Bbb{H}}^n$) and let
$\mathcal{H}$  denote a uniformly separated
collection of horoballs and thickened geodesics.
Let $\eta$ be the hyperbolic geodesic with the same endpoints as
$\gamma$. Let $\cal{H}({\eta})$  
be the union of all the horoballs and thickened geodesics
 in $\cal{H}$ meeting $\eta$. Then
$\eta\cup\mathcal{H}{({\eta})}$ is (uniformly) quasiconvex and $\gamma
(I) \subset  
B_R (\eta \cup \cal{H} ({\eta}))$, where $R$ depends only on
$K$. 
\label{ctm}
\end{theorem}

\section{Universal Covers of Building Blocks and Electric Geometry}

\subsection{Graph Model of Building Blocks}

{\bf Amalgamation Blocks}

\smallskip

Given a geodesic segment $\lambda \subset \widetilde{S}$ and a basic 
{\it amalgamation building block } $B$, let $\lambda = [a,b] \subset
\widetilde{S} \times \{ 0
\}$ be a geodesic segment, where $\widetilde{S} \times \{ 0
\} \subset \widetilde{B}$.

We shall now build a graph model for $\widetilde{B}$ which will be
quasi-isometric to an electrocuted version of the original model,
where amalgamation components of the geometric core $K$
are electrocuted.

 $\widetilde{S} \times \{ 0 \} $ and $\widetilde{S} \times \{ 3 \}$
are equipped with hyperbolic metrics.
  $\widetilde{S} \times \{ 1 \}$ and $\widetilde{S}
\times  \{ 2 \}$ are grafted surfaces with electric metric obtained by
electrocuting the 
amalgamation components. 
This constructs $4$ `sheets' of $\widetilde{S}$
comprising the `horizontal skeleton' of the `graph model' of
$\widetilde{B}$. Now for the vertical strands. On each vertical
element of the form $x \times [0,1]$ and $x \times [2,3]$ put the
Euclidean metric. 

To do this precisely, one needs to take a bit more care and perform
the construction in the universal cover. For each
amalgamation component of $\widetilde{K}$ (recall that such
a component is a  lift of an amalgamation component of $K$ to the 
universal cover along with bounding lifts $\widetilde{T_\sigma}$
 of the Margulis tubes). For
each such component $\widetilde{K_i}$ we construct $\widetilde{K_i}
\times [0,1/2]$, so that any two copies $\widetilde{K_i}
\times [0,1/2]$ and $\widetilde{K_j}
\times [0,1/2]$ intersect (if at all they do)
 only along the original bounding 
lifts $\widetilde{T_\sigma}$
 of the Margulis tubes. In particular the copies $\widetilde{K_i}
\times [0,1/2]$ intersect $\widetilde{K}$ along $\widetilde{K_i}
\times  \{ 0 \}$. 
Next put the zero metric on each copy of $\widetilde{K_i}
\times \{ 1/2 \}$. 

This construction is very closely related to the `coning'
construction introduced by Farb in \cite{farb-relhyp}. 

The resulting copy of $\widetilde{B}$ will be called the {\bf graph
  model of an amalgamation block}.

Next, we give an $I$-bundle structure to $K$ that preserves the
grafting annulus. Thus $A_\sigma \times [1,2]$ has a structure of a
Margulis tube. Let $\phi$ denote a map from  $S\times \{ 1 \}$ to
$S\times \{ 2 \}$ mapping $(x,1)$ to $(x,2)$. Clearly there is a bound
$l_B$ on the length in $K$ of $x \times [1,2]$ as $x$ ranges over $S
\times \{ 1 \}$. That is to say that the core $K$ has a bounded
thickness. This bound depends on the block $B$ we are considering.

Let $\tilde{\phi}$ denote the lift of $\phi$ to $\widetilde{K}$ Then
$\tilde{\phi}$ is a $(k, \epsilon )$-quasi-isometry where $k,
\epsilon$ depend on the block $B$. 

\smallskip

{\bf Thick Block}

\smallskip

For a thick block $B = \widetilde{S} \times [0,1]$, recall that $B$ is
the universal curve over a `thick' Teichmuller geodesic $\lambda_{Teich} =
[a,b]$
of length less than
some fixed $D > 0$. Each $S \times \{ x \}$ is identified with the
hyperbolic surface over $(a + \frac{x}{b-a} )$ (assuming that the
Teichmuller geodesic is parametrized by arc-length). 

Here $S \times \{ 0 \} $ is identified with the hyperbolic surface
corresponding to $a$, $S \times \{ 1 \}$ is identified with the
hyperbolic surface corresponding to $b$ and each $(x,a)$ is joined to
$(x,b)$ by a segment of length $1$. 

The resulting model of $\widetilde{B}$ is called a {\bf graph model of
  a thick block}.

Metrics on graph models are called {\bf graph metrics}.

\smallskip

{\bf Admissible Paths}

\smallskip

Admissible paths consist of the following :

\begin{enumerate}

\item Horizontal segments along some $\widetilde{S} \times \{ i \}$
  for $ i = \{ 0,1,2,3 \}$ (amalgamated blocks) or $i = \{ 0, 1 \}$ (thick
  blocks).

\item Vertical segments $x \times [0,1]$ or $x \times [2,3]$ for amalgamated
  blocks
or $x \times [0,1]$ for thick blocks.

\item Vertical segments of length $\leq l_B$ joining $x  \times \{ 1 \}$ to
$x  \times \{ 2 \}$ for amalgamated blocks.

\end{enumerate}

\subsection{Construction of Quasiconvex Sets for Building Blocks}

In the next section, we will construct a set $B_\lambda$ containing
$\lambda$ and a retraction $\Pi_\lambda$ of $\widetilde{M}$ onto
it. $\Pi_\lambda$ will have the property that it does not stretch
distances much. This will show that $B_\lambda$ is quasi-isometrically
embedded in $\widetilde{M}$. 

In this subsection, we describe the construction of $B_\lambda$
restricted to a building block $B$.

\smallskip

{\bf Construction of $B_\lambda (B)$ - Thick Block}

\smallskip

Let the thick block be the universal curve over a Teichmuller geodesic
$[\alpha , \beta ]$. Let $S_\alpha$ denote the hyperbolic surface over
$\alpha$ and $S_\beta$ denote the hyperbolic surface over $\beta$.

First, let $\lambda = [a,b]$ be a geodesic segment in
$\widetilde{S}$. Let $\lambda_{B0}$ denote $\lambda \times \{ 0
\}$. 

Next, let  $\psi$ be  the lift of the 'identity' map from 
$\widetilde{S_\alpha}$ to
$\widetilde{S_\beta}$. 
. Let $\Psi$ denote
the induced map on geodesics and let $\Psi ( \lambda )$ denote the
hyperbolic geodesic 
joining $\psi (a), \psi (b)$. Let $\lambda_{B1}$ denote $\Psi (\lambda
) \times \{ 1 \}$.

For the universal cover $\widetilde{B}$ of the thick block $B$, define:

\begin{center}

$B_\lambda (B) = \bigcup_{i=0,1} \lambda_{Bi}$

\end{center}

{\bf Definition:} Each $\widetilde{S} \times i$ for $i = 0, 1$
 will be called a {\bf horizontal sheet} of
$\widetilde{B}$ when $B$ is a thick block.

\smallskip

{\bf Construction of $B_\lambda (B)$ - Amalgamation Block}

\smallskip

First, recall that $\lambda = [a,b]$ is a geodesic segment in
$\widetilde{S}$. Let $\lambda_{B0}$ denote $\lambda \times \{ 0
\}$. 

Next, let $\lambda_{Gel}$ denote the electric geodesic joining $a,
b$ in the electric pseudo-metric on $\widetilde{S}$ obtained by
electrocuting lifts of $\sigma$. Let $\lambda_{B1}$ denote
$\lambda_{Gel} \times \{ 1 \}$. 

Third, recall that $\tilde{\phi }$ is the lift of a component
preserving diffeomorphism  $\phi$ to
$\widetilde{S}$ equipped with  the electric metric $d_{Gel}$. Let
$\tilde{\Phi }$ denote
the induced map on geodesics, i.e. if $\mu = [x,y] \subset (
\widetilde{S} , d_{Gel} )$, then $\tilde\Phi ( \mu ) = [ \phi (x), \phi (y)
]$ is the geodesic joining $\phi (x), \phi (y)$. Let $\lambda_{B2}$
denote $\Phi ( \lambda_{Gel} ) \times \{ 2 \}$.

Fourthly, let $\Phi ( \lambda )$ denote the hyperbolic geodesic
joining $\phi (a), \phi (b)$. Let $\lambda_{B3}$ denote $\Phi (\lambda
) \times \{ 3 \}$.

For the universal cover $\widetilde{B}$ of the thin block $B$, define:

\begin{center}

$B_\lambda (B) = \bigcup_{i=0,\cdots , 3} \lambda_{Bi}$

\end{center}

{\bf Definition:} Each $\widetilde{S} \times i$ for $i = 0
\cdots 3$ will be called a {\bf horizontal sheet} of
$\widetilde{B}$ when $B$ is a thick block.

\smallskip

{\bf Construction of $\Pi_{\lambda ,B}$ - Thick Block}

\smallskip

On
$\widetilde{S} \times \{ 0 \}$, let $\Pi_{B0}$ denote nearest point
projection
onto $\lambda_{B0}$ in the path metric on $\widetilde{S} \times \{ 0 \}$. 

On
$\widetilde{S} \times \{ 1 \}$, let $\Pi_{B1}$ denote nearest point
projection
onto $\lambda_{B1}$ in the path metric on $\widetilde{S} \times \{ 1 \}$.

For the universal cover $\widetilde{B}$ of the thick block $B$, define:

\begin{center}

$\Pi_{\lambda ,B}(x) = \Pi_{Bi}(x) , x \in \widetilde{S} \times \{ i
  \} , i=0,1$

\end{center}

\smallskip

{\bf Construction of $\Pi_{\lambda ,B}$ - Amalgamation Block}

\smallskip

On 
$\widetilde{S} \times \{ 0 \}$, let $\Pi_{B0}$ denote nearest point
projection onto $\lambda_{B0}$. Here the nearest point projection is
taken in the path metric on $\widetilde{S} \times \{ 0 \}$ which is a
hyperbolic 
metric space.

On $\widetilde{S} \times \{ 1 \}$, let $\Pi_{B1}$ 
denote  the nearest point projection onto $\lambda_{B1}$. Here the
nearest point projection is taken in the sense of the definition
preceding
Lemma \ref{hyp=elproj}, i.e.  minimising the ordered pair $(d_{Gel},
d_{hyp})$ (where $d_{Gel}, d_{hyp}$ refer to electric and hyperbolic
metrics respectively.)

On $\widetilde{S} \times \{ 2 \}$, let $\Pi_{B2}$ 
denote  the nearest point projection onto $\lambda_{B2}$. Here, again the
nearest point projection is taken in the sense of the definition
preceding
Lemma \ref{hyp=elproj}.

Again, on 
$\widetilde{S} \times \{ 3 \}$, let $\Pi_{B3}$ denote nearest point
projection onto $\lambda_{B3}$. Here the nearest point projection is
taken in the path metric on $\widetilde{S} \times \{ 3 \}$ 
which is a hyperbolic
metric space.

For the universal cover $\widetilde{B}$ of the thin block $B$, define:

\begin{center}

$\Pi_{\lambda ,B}(x) = \Pi_{Bi}(x) , x \in \widetilde{S} \times \{ i
  \} , i=0,\cdots , 3$

\end{center}

{\bf $\Pi_{\lambda , B}$ is a retract - Thick Block}

The proof for a thick block is exactly as in \cite{mitra-trees} and
\cite{brahma-ibdd}. We omit it here.

\begin{lemma} (Lemma 4.1 of \cite{brahma-ibdd}
There exists $C > 0$ such that the following holds: \\
Let $x, y \in \widetilde{S} \times \{ 0, 1\} \subset \widetilde{B}$
for some thick block $B$. 
Then $d( \Pi_{\lambda , B} (x), \Pi_{\lambda , B} (y)) \leq C d(x,y)$.
\label{retract-thick}
\end{lemma}

{\bf $\Pi_{\lambda , B}$ is a retract - Amalgamation Block}

The  main ingredient in this case is Lemma \ref{hyp=elproj}. 

\smallskip

\begin{lemma}
There exists $C > 0$ such that the following holds: \\
Let $x, y \in \widetilde{S} \times \{ 0,1,2,3 \} \subset \widetilde{B}$
for some amalgamated block $B$. 
Then $d_{Gel}( \Pi_{\lambda , B} (x), \Pi_{\lambda , B} (y)) \leq C d_{Gel}(x,y)$.
\label{retract-thin}
\end{lemma}

{\bf Proof:} It is enough to show this for the  following cases: \\

\begin{enumerate}

\item $x, y \in \widetilde{S} \times \{ 0 \} $ OR 
 $x, y \in \widetilde{S} \times \{ 3 \} $. \\

\item $x = (p,0)$ and $y = (p,1)$ for some $p$ \\

\item $x, y$ both lie in the geometric core $K$ \\

\item $x = (p,2)$ and $y = (p,3)$ for some $p$. \\

\end{enumerate}

{\bf Case 1:} This follows from Lemma \ref{easyprojnlemma}

\smallskip

{\bf Case 2 and Case 4:} These follow from Lemma \ref{hyp=elproj} which says
that the hyperbolic and electric projections of $p$ onto the
hyperbolic geodesic $[a,b]$ and the electro-ambient geodesic
$[a,b]_{ea}$ respectively `almost agree'.  If $\pi_h$ and $\pi_e$
denote the hyperbolic and electric projections, then
there exists $C_1 > 0$ such that

\begin{center}

$d_{Gel}( \pi_h(p), \pi_e(p)) \leq C_1$

\end{center}

Hence

\begin{center}

 $d_{Gel}( \Pi_{\lambda , B} ((p,i)), \Pi_{\lambda , B} ((p,i+1)))
 \leq C_1 + 1$, for $i = 0, 2$.

\end{center}

\smallskip

{\bf Case 3:} This follows from the fact that $K$ in the graph model
with the electric metric is essentially the tree coming from the
splitting. Further, by the properties of $\pi_e$, each amalgamation
component projects down to a set of diameter zero. Hence \\

\begin{center}

 $d_{Gel}( \Pi_{\lambda , B} (p), \Pi_{\lambda , B} (q))
 \leq C_1 + 1$

\end{center}

Choosing $C$ as the maximum of these constants, we are through. $\Box$

\section{Construction of Quasiconvex Sets and Quasigeodesics}

\subsection{Construction of $B_\lambda$ and $\Pi_\lambda$}

Given a manifold $M$ of amalgamated geometry, we know that $M$ is
homeomorphic to $S \times J$ for $J = [0, \infty ) $ or $( - { \infty
  }, {\infty })$. By definition of amalgamated geometry, 
there exists a sequence $I_i$ of intervals and blocks $B_i$
where the metric on $S \times I_i$ coincides
  with that on some  building block $B_i$. Denote: \\
$\bullet$ $B_{\mu , B_i} = B_{i \mu }$ \\
$\bullet$ $\Pi_{\mu , B_i} = \Pi_{i \mu }$ \\

Now for a block $B = S \times I$ (thick or amalgamated),  a natural
map $\Phi_B$ may be defined taking
 $ \mu = B_{\mu , B} \cap \widetilde{S} \times \{ 0 \} $ to a
geodesic $B_{\mu , B} \cap \widetilde{S} \times \{ k \} = \Phi_B ( \mu
)$ where $k = 1$ or $3$ according as $B$ is thick or amalgamated. Let the map
$\Phi_{B_i}$ be denoted as $\Phi_i$ for $i \geq 0$.
For $i < 0$ we shall modify this by defining $\Phi_i$ to be the map
that takes 
 $ \mu = B_{\mu , B_i} \cap \widetilde{S} \times \{ k \} $ to a
geodesic $B_{\mu , B_i} \cap \widetilde{S} \times \{ 0 \} = \Phi_i ( \mu
)$ where $k = 1$ or $3$ according as $B$ is thick or amalgamated. 

We start with a reference block $B_0$ and a reference geodesic segment
$\lambda = \lambda_0$ on the `lower surface' $\widetilde{S} \times \{
0 \}$.
Now inductively define: \\
$\bullet$ $\lambda_{i+1}$ = $\Phi_i ( \lambda_i )$ for $i \geq 0$\\
$\bullet$ $\lambda_{i-1}$ = $\Phi_i ( \lambda_i )$ for $i \leq 0$ \\
$\bullet$ $B_{i \lambda }$ = $B_{\lambda_i } ( B_i )$ \\
$\bullet$ $\Pi_{i \lambda }$ = $\Pi_{\lambda_i , B_i}$ \\
$\bullet$ $B_\lambda = \bigcup_i B_{i \lambda }$ \\
$\bullet$ $\Pi_\lambda = \bigcup_i \Pi_{i \lambda }$ \\

Recall that 
each $\widetilde{S} \times i$ for $i = 0
\cdots m$ is called a {\bf horizontal sheet} of
$\widetilde{B}$, where $m = 1$ or $3$ according as $B$ is thick or
amalgamated. We will restrict our attention to the union of the horizontal
sheets $\widetilde{M_H}$ of $\widetilde{M}$ with the 
metric induced from the graph model. 

Clearly, 
$B_\lambda  \subset \widetilde{M_H} \subset \widetilde{M}$, and
$\Pi_\lambda$ is defined from  $\widetilde{M_H}$ to $B_\lambda$. Since
 $\widetilde{M_H}$ is a `coarse net' in $\widetilde{M}$ (equipped with
the {\it graph model metric}), we will be
able to get all the coarse information we need by restricting
ourselves to  $\widetilde{M_H}$.

By Lemmas \ref{retract-thick} and \ref{retract-thin}, we obtain the
fact that each $\Pi_{i \lambda }$ is a retract. Hence assembling all
these retracts together, we have the following basic theorem:

\begin{theorem}
There exists $C > 0$ such that 
for any geodesic $\lambda = \lambda_0 \subset \widetilde{S} \times \{
0 \} \subset \widetilde{B_0}$, the retraction $\Pi_\lambda :
\widetilde{M_H} \rightarrow B_\lambda $ satisfies: \\

Then $d_{Gel}( \Pi_{\lambda , B} (x), \Pi_{\lambda , B} (y)) \leq C
d_{Gel}(x,y) + 
C$.
\label{retract}
\end{theorem}

{\bf Note 1} For Theorem \ref{retract} above, note that all that we
really
require is that the universal cover $\widetilde{S}$ be a hyperbolic
metric space. There is {\em no restriction on $\widetilde{M_H}$.} In fact,
Theorem \ref{retract} would hold for general stacks of hyperbolic
metric spaces with  blocks of amalgamated geometry. 

{\bf Note 2:} $M_H$ has been  built up out of {\bf graph models
  of thick and amalgamated blocks} and have sheets that are
  electrocuted along geodesics.

We want to make {\em Note 1} above explicit.
We first modify the definition of amalgamation geometry as follows,
retaining only local quasiconvexity.

\noindent {\bf Definition:} A manifold $M$ homeormorphic to $S \times
  J$, where $ J = [0,
  {\infty })$ or $J = ( - \infty , \infty )$, is said to be a model of {\bf
weak  amalgamation geometry} if  \\

\begin{enumerate}

\item there is a fiber preserving homeomorphism from $M$ to
  $\widetilde{S} \times J$ 
that lifts to  a quasi-isometry of universal covers \\ 

\item there exists a sequence $I_i$ of intervals (with disjoint
 interiors)
 and blocks $B_i$
where the metric on $S \times I_i$ is the same as
 that on some  building block $B_i$. Each block is either thick or has
 amalgamation geometry. \\

\item $\bigcup_i I_i = J$ \\

\item There exists $ C_0 > 0$ such that for all amalgamated blocks $B_i$
  and geometric cores $C \subset B_i$, all amalgamation components of
  $\widetilde{C}$ are $C_0$-quasiconvex in ${\widetilde{B}}_i$ \\

\end{enumerate}

Then as a consequence of {\em the proof} of Theorem \ref{retract}, we
have the following Corollary.

\begin{cor}
Let $M$ be a model manifold of {\bf weak amalgamation geometry}.
There exists $C > 0$ such that the following holds: \\
 Given
any geodesic $\lambda \subset \widetilde{S} \times \{ 0 \}$, let
 $B_\lambda$, $\Pi_\lambda$ be as before. Then  
for  $\lambda = \lambda_0 \subset \widetilde{S} \times \{
0 \} \subset \widetilde{B_0}$, the retraction $\Pi_\lambda :
\widetilde{M_H} \rightarrow B_\lambda $ satisfies: \\

Then $d_{Gel}( \Pi_{\lambda , B} (x), \Pi_{\lambda , B} (y)) \leq C
d_{Gel}(x,y) + 
C$.
\label{retract-cor}
\end{cor}

In fact, all that follows in this section may  just as well be done for
model manifolds of {\em weak amalgamation geometry}. We shall make
this explicit again at the end of this entire section.

Bur before we proceed, we would like to deduce one further Corollary
of Theorem \ref{retract}, which shall be useful towards the end of the
paper. Instead of constructing vertical hyperbolic ladders $B_\lambda$
for finite geodesic segments, first note that $\lambda$ might as well
be bi-infinite. Next, we would like to construct such a $B_\lambda$
{\bf equivariantly} under the action of $\Bbb{Z}$. That is to say, we
would like to construct a vertical annulus in the manifold $M$
homeomorphic to $S \times \Bbb{R}$.

To do this, we start with a simple closed geodesic $\sigma$
on $S \times \{ 0 \}$. Instead of performing the construction in the
universal cover, homotop $\sigma$ into $S \times \{ i \}$ for each
level $i$. Let $\sigma_i$ denote the shortest electro-ambient geodesic
in the {\em free homotopy class} of $\sigma \times \{ i \}$ in the
path pseudometric on $S \times \{ i \}$. Now let $B_\sigma$ denote the
set $B_\sigma = \bigcup_i {\tilde{\sigma}}_i$. Then the proof of
Theorem \ref{retract} ensures the quasiconvexity of $B_\sigma$ in the
{\em graph-metric}. Finally, since $B_\sigma$ has been constructed to
be equivariant under the action of the surface group, its quotient in
$M$ is an embedded `quasi-annulus' $A_{P\sigma}$ 
which partitions the manifold
locally. We use the term `quasi-annulus' because $A_{P\sigma}$ is a
collection of disjoint circles at different levels. We finally
conclude:

\begin{cor}
Let $M$ be a model manifold of {\bf weak amalgamation geometry}.
There exists $C > 0$ such that the following holds: \\
 Given
any simple closed 
geodesic $\sigma \subset {S} \times \{ 0 \}$, let
 $B_\sigma$ be as above. Then its quotient, the embedded quasi-annulus
 $A_{P\sigma}$ above is $C$-quasiconvex in $M$ with the graph metric.
\label{qcannulus}
\end{cor}

Another Corollary will be used later. Suppose $\Sigma = \Sigma \times
\{ 0 \}$ be a subsurface
of $S \times \{ 0 \}$ 
with geodesic boundary components $\sigma^1 \cdots
\sigma^k$. Let $\Sigma_i$ be the subsurface of $S \times \{ i \}$ that
is bounded by $\sigma^1_i \cdots \sigma^k_i$. Let $B_\Sigma =
\bigcup_i \Sigma_i$.

\begin{cor}
Let $M$ be a model manifold of {\bf weak amalgamation geometry}.
There exists $C > 0$ such that the following holds: \\
 Given
any  
subsurface $\Sigma \subset {S} \times \{ 0 \}$ with geodesic boundary
components, let
 $B_\Sigma$ be as above. Then $B_\Sigma$
 is $C$-quasiconvex in $M$ with the graph metric.
\label{qcsubsurf}
\end{cor}

\subsection{Heights of Blocks}

Recall that each geometric core $C \subset B$ is identified with $S \times I$
where each fibre $\{ x \} \times I$ has length $\leq l_C$ for some
$l_C$, called the {\it thickness} of the block $B$. If $C \subset B_i$
for one of the above blocks $B_i$, we shall denote $l_C$ as $l_i$. 

Instead of considering all the horizontal sheets, we would now like
to consider only the {\bf boundary horizontal sheets}, i.e. for a thick
block we consider $\widetilde{S} \times \{ 0, 1\}$ and for a thin
block we consider $\widetilde{S} \times \{ 0,3 \}$. The union of all 
boundary horizontal sheets will be denoted by $M_{BH}$. 

\smallskip

{\bf Observation 1:} $\widetilde{M_{BH}}$ is a `coarse net' in $\widetilde{M}$
in the {\bf graph model}, but not in the {\bf model of amalgamated geometry}.

In the graph model, any point can be connected by a vertical segment
of length $\leq 2$ to one of the boundary horizontal sheets.

However, in the model of amalgamated geometry, there are points within
amalgamation components which are at a distance of the order of $l_i$ from the boundary horizontal sheets. Since $l_i$ is arbitrary,
$\widetilde{M_{BH}}$ is no longer a `coarse net' in $\widetilde{M}$
equipped with 
the model of {\it amalgamated geometry}.

{\bf Observation 2:} $\widetilde{M_H}$ is defined only in the {\bf
  graph model}, but not in the model of amalgamated geometry.

{\bf Observation 3:} The electric metric on the model of amalgamated
geometry on $\widetilde{M}$ obtained by
electrocuting amalgamation components is quasi-isometric to the
graph model of $\widetilde{M}$.

\smallskip

{\bf Bounded Height of Thick Block}

\smallskip

Let $\mu \subset \widetilde{S} \times \{ 0 \} \widetilde{B_i}$ be a
geodesic in a (thick or amalgamated) block. 
Then  there exists a $(K_i, \epsilon_i )$- quasi-isometry $\psi_i$ ( =
$\phi_i$ for thick blocks)
from $\widetilde{S} \times \{ 0 \}$ to $\widetilde{S} \times \{ 1 \}$
and $\Psi_i$ is the induced map on geodesics. Hence, for any $x \in
\mu$, $\psi_i (x)$ lies within some bounded distance $C_i$ of $\Psi_i
( \mu )$. But $x$ is connected to $\psi_i (x)$ by \\

\smallskip

\noindent {\it Case 1 - Thick Blocks:} a vertical segment of length $1$  \\
{\it Case 2 - Amalgamated Blocks:}  the union of 

\begin{enumerate}

\item
two vertical segments of length $1$ between $\widetilde{S}
\times \{ i \}$ and $\widetilde{S} \times \{ i + 1 \}$ for $i = 0, 2$
\\

\item a horizontal segment of length bounded by (some uniform)
  $C^{\prime}$   (cf. Lemma \ref{ea})
connecting $(x,1)$ to a point on the
  electro-ambient geodesic 
  $B_\lambda (B) \cap \widetilde{S} \times \{ 1 \}$ \\

\item a vertical segment of electric length zero in the {\bf graph model}
 connecting $(x,1)$ to $(x, 2)$. Such a path has to travel
 {\em through an amalgamated block} in the model of {\bf amalgamated
 geometry} and has length less than
$l_i$, where $l_i$ is the thickness of the $i$th block $B_i$. \\

\item a horizontal segment of length less than $C^{\prime}$ (Lemma \ref{ea})
connecting $(\phi_i (x),3)$ to a point on the
  hyperbolic geodesic 
  $B_\lambda (B) \cap \widetilde{S} \times \{ 3 \}$ \\

\end{enumerate}

\smallskip

Thus $x$ can be connected to a point on $x^{\prime} \in \Psi_i ( \mu
)$ by a path of length less than $g(i) = 2 + 2C^{\prime} + l_i$.
 Recall that $\lambda_i$ is the geodesic on the lower
horizontal surface of the block $\widetilde{B_i}$. The same can be
done for blocks $\widetilde{B_{i-1}}$ and {\em going down} from
$\lambda_i$ to $\lambda_{i-1}$.
What we have thus shown is:

\begin{lemma}
There exists a function $g: \Bbb{Z} \rightarrow \Bbb{N}$ such that
for any block $B_i$ (resp. $B_{i-1}$), and $x \in \lambda_i$, there exists $x^{\prime} 
\in \lambda_{i+1}$ (resp. $\lambda_{i-1}$) for $i \geq 0$ (resp. $i \leq
0$), satisfying:

\begin{center}

$d(x, x^{\prime}) \leq g(i)$

\end{center}
\label{bddheight}
\end{lemma}

\smallskip

\subsection{Admissible Paths}

We want to define a collection of {\bf $B_\lambda$-elementary admissible paths}
 lying in a bounded
neighborhood of  $B_\lambda$. $B_\lambda$ is not connected. Hence, it
 does not make much sense to speak of the path-metric on
 $B_\lambda$. To remedy this we introduce a `thickening'
 (cf. \cite{gromov-ai}) of $B_\lambda$ which is path-connected and
 where the
 paths are controlled. A {\bf  $B_\lambda$-admissible path} will be a composition
 of  $B_\lambda$-elementary admissible paths.

Recall
that admissible paths in the graph model of bounded geometry 
consist of the following :

\begin{enumerate}

\item Horizontal segments along some $\widetilde{S} \times \{ i \}$
  for $ i = \{ 0,1,2,3 \}$ (amalgamated blocks) or $i = \{ 0, 1 \}$ (thick
  blocks).

\item Vertical segments $x \times [0,1]$ or $x \times [2,3]$ for amalgamated
  blocks, where $x \in \widetilde{S}$.

\item Hyperbolic geodesic segments of length 
 $\leq l_B$ in $K \subset B$ joining $x  \times \{ 1 \}$ to
$x  \times \{ 2 \}$ for amalgamated blocks.

\item Vertical segments of length $1$ joining $x  \times \{ 0 \}$ to
$x \times \{ 1 \}$ for thick blocks.

\end{enumerate}

We shall choose a subclass of these admissible paths to define
$B_\lambda$-elementary admissible paths.

\smallskip

{\bf  $B_\lambda$-elementary admissible paths in the thick block}

Let $B = S \times [i,i+1]$ be a thick block, where each $(x,i)$ is
connected by a vertical segment of length $1$ to $( x,i+1)$. Let
$\phi$ be the map that takes $(x,i)$ to $(x,i+1)$.
 Also  $\Phi$ is the map on geodesics induced by
$\phi$. Let $B_\lambda \cap \widetilde{B} = \lambda_i \cup
\lambda_{i+1}$ where $\lambda_i$ lies on $\widetilde{S} \times \{ i
\}$ and $\lambda_{i+1}$ lies on  $\widetilde{S} \times \{ i+1
\}$.  $\pi_j$, for $j = i, i+1$  denote nearest-point projections of
$\widetilde{S} \times \{ j \}$ onto $\lambda_j$. Next, since $\phi$
is a quasi-isometry, there exists $C > 0$ such that for all $(x,i) \in
\lambda_i$, $(x,i+1)$ lies in a $C$-neighborhood of $\Phi
(\lambda_i ) = \lambda_{i+1}$. The same holds for $\phi^{-1}$ and
points in $\lambda_{i+1}$, where $\phi^{-1}$ denotes the {\em
  quasi-isometric inverse} of $\phi$ from
$\widetilde{S} \times \{ i + 1 \}$ to $\widetilde{S} \times \{ i \}
$. The {\bf
  $B_\lambda$-elementary admissible paths} in $\widetilde{B}$ consist
of the following:

\begin{enumerate}

\item Horizontal geodesic subsegments of $\lambda_j$,  $j = \{ i, i+ 1 \}$.

\item Vertical segments of length $1$ joining $x  \times \{ 0 \}$ to
$x  \times \{ 1 \}$.

\item Horizontal geodesic segments lying in a $C$-neighborhood of
  $\lambda_j$, $j = i, i+1$.

\end{enumerate}

{\bf  $B_\lambda$-elementary admissible paths in the amalgamated block}

Let $B = S \times [i,i+3]$ be an amalgamated block, where each $(x,i+1)$ is
connected by a geodesic segment of zero electric length and hyperbolic
length $\leq C(B)$ (due to bounded thickness of $B$) to $( \phi
(x),i+2)$ (Here $\phi$ can be thought of as the  map from
$\widetilde{S} \times \{ i+1 \}$ to .$\widetilde{S} \times \{ i+2 \}$
that is the identity on the first component.
 Also  $\Phi$ is the map on canonical representatives of electric geodesics induced by
$\phi$. Let $B_\lambda \cap \widetilde{B} = \bigcup_{j=i \cdots i+3}
 \lambda_j$ where $\lambda_j$ lies on $\widetilde{S} \times \{ j 
\}$.  $\pi_j$  denotes nearest-point projection of
$\widetilde{S} \times \{ j \}$ onto $\lambda_j$ (in the appropriate
sense - hyperbolic for $j = i, i+3$ and electric for $j = i+1,
i+2$). Next, since $\phi$ 
is an electric isometry, but a hyperbolic quasi-isometry, there exists
$C > 0$ (uniform constant) and $K=K(B)$  such that for all $(x,i) \in
\lambda_i$, $(\phi(x),i+1)$ lies in an (electric) $C$-neighborhood and
a hyperbolic $K$-neighborhood of $\Phi
(\lambda_{i+1} ) = \lambda_{i+2}$. The same holds for $\phi^{-1}$ and
points in $\lambda_{i+2}$, where $\phi^{-1}$ denotes the {\em
  quasi-isometric inverse} of $\phi$ from
$\widetilde{S} \times \{ i + 2 \}$ to $\widetilde{S} \times \{ i + 1\}
$. 

Again, since $\lambda_{i+1}$ and $\lambda_{i+2}$ are electro-ambient
quasigeodesics, we further note that there exists $C > 0$ (assuming
the same $C$ for convenience) such that for all $(x,i) \in \lambda_i$,
$(x,i+1)$ lies in a (hyperbolic) $C$-neighborhood of
$\lambda_{i+1}$. Similarly for all $(x,i+2) \in \lambda_{i+2}$,
$(x,i+3)$ lies in a (hyperbolic) $C$-neighborhood of
$\lambda_{i+3}$. The same holds if we go `down' from $\lambda_{i+1}$
to $\lambda_i$ or from $\lambda_{i+3}$ to $\lambda_{i+2}$.
The {\bf
  $B_\lambda$-elementary admissible paths} in $\widetilde{B}$ consist
of the following:

\begin{enumerate}

\item Horizontal  subsegments of $\lambda_j$,  $j = \{ i,\cdots i+ 3 \}$.

\item Vertical segments of length $1$ joining $x  \times \{ j  \}$ to
$x  \times \{ j+1 \}$, for $j = i, i+2$.

\item Horizontal geodesic segments lying in a {\em hyperbolic}
 $C$-neighborhood of
  $\lambda_j$, $j = i,\cdots i+3$.

\item Horizontal hyperbolic segments of {\em electric length $\leq C$}
  and {\em hyperbolic length $\leq K(B)$} joining points of the form
  $( \phi (x), i+2)$ to a point on $\lambda_{i+2}$ for $(x, i+1) \in
  \lambda_{i+1}$. 

\item Horizontal hyperbolic segments of {\em electric length $\leq C$}
  and {\em hyperbolic length $\leq K(B)$} joining points of the form
  $( \phi^{-1} (x), i+1)$ to a point on $\lambda_{i+1}$ for $(x, i+2) \in
  \lambda_{i+2}$. 

\end{enumerate}

{\bf Definition:} A   $B_\lambda$-admissible path is a union of
$B_\lambda$-elementary admissible paths.

The following lemma follows from the above definition and Lemma
\ref{bddheight}.

\begin{lemma}
There exists a function $g: \Bbb{Z} \rightarrow \Bbb{N}$ such that
for any block $B_i$, and $x$ lying on a $B_\lambda$-admissible
path in $\widetilde{B_i}$, there exist $y
\in \lambda_{j}$ and $z \in \lambda_k$ where $\lambda_j \subset B_\lambda$ and
$\lambda_k \subset B_\lambda$ lie on the two boundary horizontal sheets, satisfying:

\begin{center}

$d(x, y) \leq g(i)$ \\
$d(x, z) \leq g(i)$ \\

\end{center}
\label{bddheight-adm}
\end{lemma}

Let $h(i) = \Sigma_{j = 0 \cdots i} g(j)$ be the sum of the values of
$g(j)$ as $j$ ranges from $0$ to $i$ (with the assumption that
increments are by $+1$ for $i \geq 0$ and by $-1$ for $i \leq
0$). Then we have from Lemma \ref{bddheight-adm} above,

\begin{cor}
There exists a function $h: \Bbb{Z} \rightarrow \Bbb{N}$ such that
for any block $B_i$, and $x$ lying on a $B_\lambda$-admissible
path in $\widetilde{B_i}$, there exist $y
\in \lambda_0 = \lambda$ such that:

\begin{center}

$d(x, y) \leq h(i)$ 

\end{center}
\label{bddht-final}
\end{cor}

{\bf Important Note:} In the above Lemma \ref{bddheight-adm} and
Corollary \ref{bddht-final}, it is important to note
that the distance $d$ is {\bf hyperbolic}, not electric. 
This is because the number $l_i$ occurring in elementary paths of
type $5$ and $6$ is a hyperbolic length depending only on $i$ (in $B_i$).

\smallskip

Next suppose that $\lambda$ lies outside $B_N(p)$, the $N$-ball about
a fixed reference point $p$ on the boundary horizontal surface
$\widetilde{S} \times \{ 0 \} \subset \widetilde{B_0}$. Then by
Corollary \ref{bddht-final}, any $x$ lying on a $B_\lambda$-admissible
path in $\widetilde{B_i}$ satisfies

\begin{center}

$d(x, p) \geq N - h(i)$ 

\end{center}

Also, since the electric, and hence hyperbolic `thickness' (the
shortest distance between its boundary horizontal sheets) is $\geq 1$,
we get,

\begin{center}

$d(x, p) \geq |i|$ 

\end{center}

Assume for convenience that $i \geq 0$ (a similar argument works,
reversing signs for $i < 0$). Then,

\begin{center}

$d(x, p) \geq min \{ i,  N - h(i) \}$ 

\end{center}

Let $h_1 (i) = h(i) + i$.  Then $h_1$ is a
monotonically increasing function on the integers.
If
 $h_1^{-1} (N)$ denote the largest positive
integer $n$ such that $h(n) \leq m$, then clearly, .
 $h_1^{-1} (N) \rightarrow \infty$
 as $N \rightarrow \infty$. We have thus shown:

\begin{lemma}
There exists a function $M(N) : \Bbb{N} \rightarrow \Bbb{N}$
such that  $M (N) \rightarrow \infty$
 as $N \rightarrow \infty$ for which the following holds:\\
For any geodesic $\lambda \subset \widetilde{S} \times \{ 0 \} \subset
\widetilde{B_0}$, a fixed reference point 
 $p \in \widetilde{S} \times \{ 0 \} \subset
\widetilde{B_0}$ and any $x$ on a $B_\lambda$-admissible path, 

\begin{center}

$d(\lambda , p) \geq N \Rightarrow d(x,p) \geq M(N)$.

\end{center}

\label{far-nopunct}
\end{lemma}

As pointed out before, the discussion and Lemmas of the previous two
subsections go through just as well in the context of {\em weak
  amalgamation geometry} manifolds. We make this explicit in the case
of Lemma \ref{far-nopunct} above.

\begin{cor}
Let $M$ be a model manifold of {\bf weak amalgamation geometry}.
Then  
there exists a function $M(N) : \Bbb{N} \rightarrow \Bbb{N}$
such that  $M (N) \rightarrow \infty$
 as $N \rightarrow \infty$ for which the following holds:\\
 Given
any geodesic $\lambda \subset \widetilde{S} \times \{ 0 \}$, let
 $B_\lambda$ be as before. 
For  $\lambda \subset \widetilde{S} \times \{ 0 \} \subset
\widetilde{B_0}$, a fixed reference point 
 $p \in \widetilde{S} \times \{ 0 \} \subset
\widetilde{B_0}$ and any $x$ on a $B_\lambda$-admissible path, 

\begin{center}

$d(\lambda , p) \geq N \Rightarrow d(x,p) \geq M(N)$.

\end{center}

\label{far-nopunct-cor}
\end{cor}

\subsection{Joining the Dots}

Recall that {\bf admissible paths} in a model manifold of bounded
geometry consist of:

\begin{enumerate}

\item Horizontal segments along some $\widetilde{S} \times \{ i \}$
  for $ i = \{ 0,1,2,3 \}$ (thin blocks) or $i = \{ 0, 1 \}$ (thick
  blocks).\\
\item Vertical segments $x \times [0,1]$ or $x \times [2,3]$ for amalgamated
  blocks.\\
\item Vertical segments of length $\leq l_i$ joining $x  \times \{ 1 \}$ to
$x  \times \{ 2 \}$ for amalgamated blocks.\\
\item Vertical segments of length $1$ joining $x  \times \{ 0 \}$ to
$x  \times \{ 1 \}$ for thick blocks.\\

\end{enumerate}

Our strategy in this subsection is: \\
$\bullet 1$ Start with an electric geodesic $\beta_e$ in
$\widetilde{M_{Gel}}$ joining the end-points of $\lambda$.
\\
$\bullet 2$ Replace it by an {\em admissible quasigeodesic}, i.e. an
admissible path that is a quasigeodesic. 
\\
$\bullet 3$ Project the intersection of the
 admissible quasigeodesic with the horizontal sheets onto
 $B_\lambda$. \\
$\bullet 4$ The result of step 3 above is disconnected. {\em Join the
   dots} using $B_\lambda$-admissible paths. \\

The end product is an electric quasigeodesic built up of $B_\lambda$
admissible paths.

Now for the first two steps: \\
$\bullet$ Since $\widetilde{B}$ (for a thick block
$B$)  has thickness $1$, any path lying in a thick block can be
pertubed to an admissible path lying in $\widetilde{B}$, changing the
length by at most a bounded multiplicative factor. \\
$\bullet$ For $B$ amalgamated, we
decompose paths into horizontal paths lying in some $\widetilde{S}
\times \{ j \}$, for $j = 0, \cdots 3$ and 
vertical paths of types (2) or (3) above. This can be done 
without altering electric length within $\widetilde{S} \times
[1,2]$. To see this, project any path $\overline{ab}$ beginning and
ending on $\widetilde{S} \times \{ 1, 2 \}$ onto $\widetilde{S} \times
\{ 1 \}$ along the fibres. To connect this to the starting and ending
points $a, b$, we have to at most adjoin vertical segments through $a,
b$. 
Note that this does not increase the electric length of
$\overline{ab}$, as the electric length is determined by the number of
amalgamation blocks that $\overline{ab}$ traverses. \\
$\bullet$ For paths lying in $\widetilde{S} \times [0,1]$ or
$\widetilde{S} \times [2,3]$, we can modify the path into an
admissible path,
changing lengths by a bounded multiplicative
constant. The result is therefore an electric quasigeodesic. \\
$\bullet$ Without
loss of generality, we can assume that the electric quasigeodesic is
one without back-tracking (as this can be done without increasing the
length of the geodesic - see \cite{farb-relhyp} or \cite{klarreich}
for instance). \\
$\bullet$ Abusing notation slightly, assume therefore that 
$\beta_e$ is an admissible electric quasigeodesic without backtracking
joining the end-points of $\lambda$. \\
This completes Steps $\bullet 1$ and $\bullet 2$. 

\smallskip

\noindent 
$\bullet$ Now act on $\beta_e \cap \widetilde{M_H}$ by $\Pi_\lambda$. From Theorem \ref{retract}, we
conclude, by restricting $\Pi_\lambda$ to the horizontal sheets of
$\widetilde{M_{Gel}}$ that the image $\Pi_\lambda ( \beta_e )$  is a
`dotted electric quasigeodesic' lying entirely on $B_\lambda$. This
completes step 3. \\
$\bullet$ Note that since $\beta_e$ consists of admissible
segments, we can arrange so that two nearest points on $\beta_e \cap
\widetilde{M_H}$ which are not 
connected to each 
other  form the end-points of
a vertical segment of type (2), (3) or (4). Let $\Pi_\lambda ( \beta_e
) \cap B_\lambda = \beta_d$, be the dotted quasigedoesic lying on
$B_\lambda$. We want to join the dots in $\beta_d$ converting it into
a {\bf connected} electric quasigeodesic built up of {\bf
  $B_\lambda$-admissible paths}.  \\
$\bullet$ For vertical segments of type (4) joining $p, q$ (say), $\Pi_\lambda
(p), \Pi_\lambda (q)$ are a bounded hyperbolic distance apart. Hence, by
the proof of Lemma \ref{retract-thick}, we can join 
 $\Pi_\lambda
(p), \Pi_\lambda (q)$ by a $B_\lambda$-admissible path of length
bounded by some $C_0$ (independent of $B$, $\lambda$).\\
$\bullet$
For vertical segments of type (2) joining $p, q$, we note that $\Pi_\lambda
(p), \Pi_\lambda (q)$ are a bounded hyperbolic distance apart. Hence, by
the proof of Lemma \ref{retract-thin}, we can join 
 $\Pi_\lambda
(p), \Pi_\lambda (q)$ by a $B_\lambda$-admissible path of length
bounded by some $C_1$ (independent of $B$, $\lambda$). \\
$\bullet$ 
This leaves us to deal with case (3). Such a  segment consists of
a segment lying within a lift of an
amalgamation block.
Such a piece has electric length one in the graph model. Its image,
too, has electric length one (See for instance, Case (3) of 
 the proof of Lemma \ref{retract-thin}, where we
noted that the projection of any amalgamation component lies within
an amalgamation component). 

After joining the dots, we can assume further that the quasigeodesic
thus obtained does not backtrack (cf \cite{farb-relhyp} and
\cite{klarreich}).

Putting all this together, we conclude:

\begin{lemma}
There exists a function $M(N) : \Bbb{N} \rightarrow \Bbb{N}$
such that  $M (N) \rightarrow \infty$
 as $N \rightarrow \infty$ for which the following holds:\\
For any geodesic $\lambda \subset \widetilde{S} \times \{ 0 \} \subset
\widetilde{B_0}$, and  a fixed reference point 
 $p \in \widetilde{S} \times \{ 0 \} \subset
\widetilde{B_0}$,
there exists a connected electric quasigeodesic $\beta_{adm}$  without 
backtracking, such that \\
$\bullet$ $\beta_{adm}$ is built up of $B_\lambda$-admissible
paths. \\
$\bullet$ $\beta_{adm}$ joins
the end-points of $\lambda$. \\
$\bullet$ 
$d(\lambda , p) \geq N \Rightarrow d(\beta_{adm},p) \geq M(N)$. \\
\label{adm-qgeod-props}
\end{lemma}

{\bf Proof:} The first two criteria follow from the discussion
preceding this lemma. The last follows from Lemma \ref{far-nopunct}
since the discussion above gives a quasigeodesic built up out of
admissible paths.
$\Box$

As in the previous subsections, Lemma \ref{adm-qgeod-props} goes
through for {\bf weak amalgamation geometry}. We state this below:

\begin{cor}
Suppose that $M$ is a manifold of {\bf weak amalgamation geometry}.
There exists a function $M(N) : \Bbb{N} \rightarrow \Bbb{N}$
such that  $M (N) \rightarrow \infty$
 as $N \rightarrow \infty$ for which the following holds:\\
For any geodesic $\lambda \subset \widetilde{S} \times \{ 0 \} \subset
\widetilde{B_0}$, and  a fixed reference point 
 $p \in \widetilde{S} \times \{ 0 \} \subset
\widetilde{B_0}$,
there exists a connected electric quasigeodesic $\beta_{adm}$  without 
backtracking, such that \\
$\bullet$ $\beta_{adm}$ is built up of $B_\lambda$-admissible
paths. \\
$\bullet$ $\beta_{adm}$ joins
the end-points of $\lambda$. \\
$\bullet$ 
$d(\lambda , p) \geq N \Rightarrow d(\beta_{adm},p) \geq M(N)$. \\
\label{adm-qgeod-props-cor}
\end{cor}

\subsection{Admissible Quasigeodesics and Electro-ambient
  Quasigeodesics}

{\bf Definition:}
We next define (as before)
 a $(k, \epsilon)$  electro-ambient quasigeodesic $\gamma$
in $\tilde{M}$ relative to the amalgamation components $\tilde{K}$
to be a $(k, \epsilon )$ quasigeodesic in the graph model of
$\tilde{M}$ such that in an ordering (from the left)
of the amalgamation components
that $\gamma$ meets, each $\gamma \cap \tilde{K}$ is a $(k, \epsilon)$ -
quasigeodesic in the induced path-metric on $\tilde{K}$.

\smallskip

This subsection is devoted to extracting an electro-ambient
quasigeodesic $\beta_{ea}$
from a $B_\lambda$-admissible quasigeodesic $\beta_{adm}$. 
 $\beta_{ea}$ shall satisfy  the property indicated by Lemma
\ref{adm-qgeod-props} above. We shall prove this Lemma under the
assumption of (strong) amalgamation geometry. However, a weaker
assumption (which we shall discuss later, while weakening {\em
  amalgamation geometry} to {\bf graph amalgamation geometry}) is
enough for the main Lemma of this subsection to go through.

\begin{lemma}
There exist $\kappa , \epsilon$
and  a function $M^{\prime}(N) : \Bbb{N} \rightarrow \Bbb{N}$
such that  $M^{\prime}(N) \rightarrow \infty$
 as $N \rightarrow \infty$ for which the following holds:\\
For any geodesic $\lambda \subset \widetilde{S} \times \{ 0 \} \subset
\widetilde{B_0}$, and  a fixed reference point 
 $p \in \widetilde{S} \times \{ 0 \} \subset
\widetilde{B_0}$,
there exists a $(\kappa , \epsilon )$
  electro-ambient quasigeodesic $\beta_{ea}$  without 
backtracking, such that \\
$\bullet$ $\beta_{ea}$ joins
the end-points of $\lambda$. \\
$\bullet$ 
$d(\lambda , p) \geq N \Rightarrow d(\beta_{ea},p) \geq M^{\prime}(N)$. \\
\label{adm-ea-props}
\end{lemma}

\noindent {\bf Proof:} From Lemma \ref{adm-qgeod-props}, we have a
$B_\lambda$ - admissible quasigeodesic $\beta_{adm}$ and a function $M(N)$
without backtracking satisfying
the conclusions of the Lemma.
Since $\beta_{adm}$ does not backtrack, we can decompose it as a union
of non-overlapping segments $\beta_1, \cdots \beta_k$, such that each
$\beta_i$ is either an admissible (hyperbolic) 
quasigeodesic lying outside 
amalgamation components, or 
a $B_\lambda$-admissible quasigeodesic lying entirely
within some amalgamation component ${\widetilde{K}}_i$. Further, since
$\beta_{adm}$ does not backtrack, we can assume that all $K_i$'s are
distinct.

We  modify $\beta_{adm}$ to an electro-ambient quasigeodesic
$\beta_{ea}$ as follows:

\begin{enumerate}
\item $\beta_{ea}$ coincides with $\beta_{adm}$ outside amalgamation
  components. \\
\item There exist $\kappa, \epsilon$
such that if some $\beta_i$ lies within an amalgamation component
  ${\widetilde{K}}_i$ then, by
  uniform quasiconvexity of the $K_i$'s, it may be
replaced by a $(\kappa , \epsilon )$
  (hyperbolic) quasigeodesic $\beta_{i}^{ea}$
joining the end-points of $\beta_i$
and lying within   ${\widetilde{K}}_i$.  
\end{enumerate}

The resultant path $\beta_{ea}$ is clearly an  electro-ambient
quasigeodesic without backtracking. Next, each component
$\beta_i^{ea}$ lies in a $C_i$ neighborhood of $\beta_i$, where $C_i$
depends only on
  the thickness $l_i$ of the amalgamation component $K_i$.

We let $C(n)$ denote the maximum of the values of $C_i$ for $K_i
\subset B_n$.
Then, as in the proof of Lemma \ref{far-nopunct}, we have for any $z
\in \beta_{ea} \cap B_n$, 

\begin{center}
$ d(z, p) \geq $ max $(n, M(N) - C(n))$
\end{center}

Again, as in Lemma \ref{far-nopunct}, this gives us a (new) function 
$M^{\prime}(N) : \Bbb{N} \rightarrow \Bbb{N}$
such that  $M^{\prime} (N) \rightarrow \infty$
 as $N \rightarrow \infty$ for which\\
$\bullet$ 
$d(\lambda , p) \geq N \Rightarrow d(\beta_{ea},p) \geq M^{\prime}(N)$. \\

This prove the Lemma. $\Box$

\smallskip

{\bf Note:} We have essentially used the following two properties of
amalgamation components in concluding Lemma \ref{adm-ea-props}:

\begin{enumerate}
\item any path lying inside an amalgamation component $\tilde{K}$
 may be replaced
  by a (uniform) hyperbolic quasigeodesic joining its end-points and
 lying within the same $\tilde{K}$ \\
\item Each electro-ambient quasigeodesic joining the end-points of an
  admissible quasigodesic  in ${\widetilde{K}} \subset
  {\widetilde{B}}_n$ lies in a (hyperbolic) $C(n)$-neighborhood of the
  latter.
\end{enumerate}

We shall have occasion to use this when we discuss {\bf
  graph-quasiconvexity}.

\section{Cannon-Thurston Maps for Surfaces Without Punctures}

It is now time to introduce hyperbolicity of $\tilde{M}$, global
quasiconvexity  of amalgamation components, (and hence) model
manifolds of (strong) amalgamation geometry. We shall
assume till the end of this section that 

\begin{enumerate}

\item there exists a hyperbolic
manifold $M$ and a homeomorphism from $\widetilde{M}$ to $\widetilde{S} \times
\Bbb{R}$. We identify $\widetilde{M}$ with $\widetilde{S} \times
\Bbb{R}$ via this homeomorphism.

\item  $\widetilde{S} \times \Bbb{R}$ admits a 
quasi-isometry $g$ to a model manifold of {\em amalgamated geometry}

\item
$g$  preserves the fibers over $\Bbb{Z} \subset \Bbb{R}$

\end{enumerate}

We shall henceforth ignore the quasi-isometry $g$ and think of
$\widetilde{M}$ itself as the universal cover of a model manifold of
{\em amalgamated geometry}.

\subsection{ Electric Geometry Revisited}

We note the following properties of the pair $(X, \mathcal{H})$
where $X$ is the graph model of $\widetilde{M}$ and $\mathcal{H}$
consists of the amalgamation components.
 There exist $C, D, \Delta$ such that

\begin{enumerate}

\item Each amalgamation component  is $C$-quasiconvex.

\item Any two amalgamation components are $1$-separated.

\item $\widetilde{M_{Gel}} = X_{Gel}$ is $\Delta$-hyperbolic, (where
  $\widetilde{M_{Gel}} = X_{Gel}$ is the electric metric on
  $\widetilde{M} = X$
  obtained by electrocuting all amalgamation components, i.e. all
  members of $\mathcal{H}$).

\item Given $K, \epsilon$, there exists $D_0$ such that if $\gamma$ 
be  a $(K, \epsilon)$ hyperbolic quasigeodesic joining $a, b$ and if
$\beta$ be a $(K, \epsilon)$ electro-ambient quasigeodesic joining $a,
b$, then $\gamma$ lies in a $D_0$ neighborhood of $\beta$.

\end{enumerate}

The first property follows from the definition of a manifold of
amalgamation geometry.

The second follows from the construction of the graph model.

The third follows from Lemma \ref{farb1A}.

The fourth follows from Lemma \ref{ea-strong}.

\subsection{Proof of Theorem}

We shall now assemble the proof
of the main Theorem.

\begin{theorem}
Let $M$ be a 3 manifold homeomorphic to $S \times J$ (for $J = [0,
  \infty ) $ or $( - \infty , \infty )$). Further suppose that $M$ has
  {\em amalgamated geometry}, where $S_0 \subset B_0$ is the lower
  horizontal surface of the building block $B_0$. Then the inclusion
  $i : \widetilde{S} \rightarrow \widetilde{M}$ extends continuously
  to a map 
  $\hat{i} : \widehat{S} \rightarrow \widehat{M}$. Hence the limit set
  of $\widetilde{S}$ is locally connected.
\label{crucial}
\end{theorem}

{\bf Proof:}
Suppose $\lambda \subset \widetilde{S}$ lies outside
a large $N$-ball about $p$. By Lemma 
\ref{adm-ea-props} we obtain an electro-ambient
 quasigeodesic without backtracking
$\beta_{ea}$
 lying outside an $M(N)$-ball about $p$ (where $M(N) \rightarrow
 \infty $ as $N \rightarrow \infty $). 

Suppose that $\beta_{ea}$ is
 a $(\kappa, \epsilon)$ electro-ambient quasigeodesic. Note that $\kappa, \epsilon$ depend on
 `the Lipschitz constant' of $\Pi_\lambda$ and hence only on
 $\widetilde{S}$ and $\widetilde{M}$.

From Property (4) above, (or Lemma \ref{ea-strong}) we find that
if $\beta^{h}$ denote the hyperbolic geodesic in
$\widetilde{M}$ joining the end-points of $\lambda$, then $\beta^h$
lies in a (uniform) $C^{\prime}$ neighborhood of $\beta_{ea}$. 

Let $M_1(N) = M(N) - C^{\prime}$.
Then $M_1(N) \rightarrow
\infty$ as $N \rightarrow \infty$. Further, the hyperbolic geodesic 
 $\beta^h$ 
 lies outside an $M_1(N)$-ball around $p$. Hence, by Lemma
 \ref{contlemma}, 
the inclusion
  $i : \widetilde{S} \rightarrow \widetilde{M}$ extends continuously
  to a map 
  $\hat{i} : \widehat{S} \rightarrow \widehat{M}$. 

Since the
  continuous image of a compact locally connected set is locally
  connected (see \cite{hock-young} )
and the (intrinsic) boundary of $\widetilde{S}$ is a circle, we
  conclude that the limit set
  of $\widetilde{S}$ is locally connected.

This proves the theorem. $\Box$

\section{Modifications for Surfaces with Punctures}

In this section, we shall describe the modifications necessary for
Theorem \ref{crucial} to go through for surfaces with punctures.

\subsection{ Partial Electrocution}

In this subsection, we indicate a modification of Farb's \cite{farb-relhyp}
 notion of
{\em strong relative hyperbolicity}
 and construction of an electric metric, described earlier in this paper.
Though much of this works in the context of relative hyperbolicity
with {\em Bounded Penetration Property} \cite{farb-relhyp}
or, equivalently, strong relative hyperbolicity \cite{bowditch-relhyp},
we shall focus on the case we need, viz. convex hyperbolic 3-manifolds
with boundary of the form $\sigma \times P$, where $P$ is either
an interval or a circle, and $\sigma$ is a horocycle of some fixed 
length $e_0$. In the universal cover, if we excise (open) horoballs, we
are left with a manifold whose boundaries are flat horospheres of the
form $\widetilde{\sigma} \times \tilde{P}$. Note that $\tilde{P} = P$ if
$P$ is an interval, and $\Bbb{R}$ if $P$ is a circle (the case for
a $(Z + Z)$-cusp ).

\smallskip

\noindent 
{\bf Partial Electrocution} of a horosphere $H$ will be defined as putting the
 zero metric in the $\tilde{\sigma}$ direction, and retaining the 
usual Euclidean metric in the $\tilde{P}$ direction.

The construction of {\it partially electrocuted}
horospheres  is half way between the spirit of
Farb's construction (in Lemmas \ref{farb1A}, \ref{farb2A},
where
the entire horosphere is coned off), and McMullen's Theorem \ref{ctm}
(where nothing is coned off, and properties of
 {\it ambient quasigeodesics} are investigated).

In the partially electrocuted case, instead of coning all of a
horosphere down to a point 
we cone only horocyclic leaves of a foliation of the horosphere.
Effectively, therefore, we have a cone-line rather a cone-point. 

 We explicitly describe  below {\it partial electrocution}
for convex hyperbolic 3-manifolds.

\smallskip

\noindent {\bf Partial Electrocution of Horospheres} \\
Let $Y$ be a convex simpy connected hyperbolic 3-manifold.
Let $\mathcal{B}$ denote a collection of horoballs. Let $X$ denote
$Y$ minus the interior of the horoballs in $\mathcal{B}$. Let 
$\mathcal{H}$ denote the collection of boundary horospheres.Then each
$H \in \mathcal{H}$ with the induced metric is isometric to a Euclidean
product $E^{n-2} \times L$ for an interval $L\subset \mathbb{R}$. 
Partially electrocute  each 
$H$ by giving it the product of the zero metric with the Euclidean metric,
i.e. on $E^{n-2}$ give the zero metric and on $L$ give the Euclidean
metric. The resulting space is exactly what one would get by gluing
to each $H$ the mapping cylinder of the projection of $H$ onto the $L$-factor.

\smallskip

Much
of what follows would go through in the following more general
setting:

\begin{enumerate}

\item $X$ is (strongly) hyperblic relative to a collection of subsets
$H_\alpha$, thought of as horospheres (and {\em not horoballs}). \\
\item For each $H_\alpha$ there is a uniform large-scale
retraction $g_alpha : H_\alpha \rightarrow L_\alpha$ to some
(uniformly) $\delta$-hyperbolic metric space $L_\alpha$, i.e. there
exist $\delta , K, \epsilon > 0$ such that for all $H_\alpha$ there exists
a $\delta$-hyperbolic $L_\alpha$ and a map 
$g_\alpha : H_\alpha \rightarrow L_\alpha$ with
$d_{L_\alpha} (g_\alpha (x), g_alpha (y)) \leq Kd_{H_\alpha}(x,y)
+ \epsilon $ for all $x, y \in H_\alpha$. \\
\item The coned off space corresponding to $H_\alpha$ is the (metric)
mapping cylinder for the map  $g_alpha : H_\alpha \rightarrow L_\alpha$. \\

\end{enumerate}

In Farb's construction $L_\alpha$ is just a single point. However,
the notions and arguments of \cite{farb-relhyp} or Klarreich
 \cite{klarreich} or
the proof of quasiconvexity of a {\it hyperbolic geodesic union
  horoballs it meets} in McMullen
\cite{ctm-locconn}
go
through even in this setting. The metric, and geodesics and quasigeodesics
in the partially electrocuted space will be referred to as the 
partially electrocuted metric $d_{pel}$, and partially
electrocuted geodesics and quasigeodesics respectively. In this
situation, we conclude as in Lemmma \ref{farb1A}:

\begin{lemma}
$(X,d_{pel})$ is a hyperbolic metric space and the sets $L_\alpha$
are uniformly quasiconvex.
\label{pel}
\end{lemma}

\noindent {\bf Note 1:} When $K_\alpha$ is a point, the last statement is a 
triviality.

\noindent {\bf Note 2:} $(X, d_{pel})$ is strongly hyperbolic relative to
the sets $\{ L_\alpha \}$. In fact the space obtained by electrocuting the
sets $L_\alpha$ in $(X,d_{pel})$ is just the space $(X,d_e)$ obtained by 
electrocuting the sets $\{ H_\alpha \}$ in $X$.

\noindent {\bf Note 3:} The proof of Lemma \ref{pel} and other such
results below follow Farb's \cite{farb-relhyp} constructions. For
instance, consider a hyperbolic geodesic $\eta$ in a convex complete
simply connected hyperbolic 3-manifold $X$. Let $H_i$, $i = 1\cdots k$
be the partially electrocuted horoballs  it
meets. Let $N(\eta )$ denote the union of $\eta$ and $H_i$'s. Let $Y$
denote $X$ minus the interiors of the $H_i$'s. The
first step is to show that $N(\eta ) \cap Y$ is quasiconvex in $(Y,
d_{pel})$. To do this one takes a hyperbolic $R$-neighborhood of
$N(\eta )$ and projects
$(Y, d_{pel})$ onto it, using the hyperbolic projection. It was shown
by Farb in \cite{farb-relhyp} that the projections of all horoballs
 are uniformly
bounded in hyperbolic diameter. (This is essentially mutual
coboundedness). Hence, given $K$, choosing
 $R$  large enough, any path that goes out of
 an $R$-neighborhood of $N( \eta )$ cannot be a $K$-partially
 electrocuted
quasigeodesic. This is the one crucial step that allows the results of
\cite{farb-relhyp}, in particular, Lemma \ref{pel}
 to go through in the context of partially
electrocuted spaces.

\smallskip

As in Lemma \ref{farb2A}, partially electrocuted quasigeodesics
and geodesics without backtracking have the same intersection patterns
with {\em horospheres and boundaries of lifts of tubes} as electric geodesics
without backtracking.  
Further, since 
electric geodesics and hyperbolic quasigeodesics have similar intersection
 patterns with {\em horoballs and lifts of tubes} it follows that
 partially electrocuted 
quasigeodesics and hyperbolic quasigeodesics have similar intersection
patterns with {\em horospheres and boundaries of lifts of tubes}. We
state this formally below: 

\begin{lemma}
Given $K, \epsilon \geq 0$, there exists $C > 0$ such that the following
holds: \\
Let $\gamma_{pel}$ and $\gamma$ denote respectively a $(K, \epsilon )$
partially electrocuted quasigeodesic in $(X,d_{pel})$ and a hyperbolic
$(K, \epsilon )$-quasigeodesic in $(Y,d)$ joining $a, b$. Then $\gamma \cap X$
lies in a (hyperbolic) $C$-neighborhood of (any representative of) 
$\gamma_{pel}$. Further, outside of  a $C$-neighborhood of the horoballs
that $\gamma$ meets, $\gamma$ and $\gamma_{pel}$ track each other.
\label{pel-track}
\end{lemma}

Next, we note that partial electrocution preserves quasiconvexity. Suppose
that $A \subset Y$ as also $A \cap H$ for all $H \in \mathcal{H}$
are $C$-quasiconvex. Then given $a, b \in A \cap X$,
the hyperbolic geodesic $\lambda$ in $X$ joining $a, b$ lies 
in a $C$-neighborhood
of $A$. Since horoballs are convex, $\lambda$ cannot backtrack. 
Let $\lambda_{pel}$
be the partially electrocuted geodesic joining $a, b \in (X, d_{pel})$.
Then by Lemma ref{pel-track} above, we conclude that 
for all $H \in \mathcal{H}$
that $\lambda$ intersects, there exist points of $\lambda_{pel}$ (hyperbolically)
near the entry and exit points of $\lambda$ with respect to $H$. Since these 
points lie near 
 $A \cap H$, and since the corresponding $L$ is quasiconvex in
$(X, d_{pel})$, we conclude that $\lambda_{pel}$ lies within a bounded
distance from $A$ near horoballs. For the rest of $\lambda_{pel}$ the conclusion 
follows from Lemma \ref{pel-track}. We conclude:

\begin{lemma} 
Given $C_0$ there exists $C_1$ such that if
 $A \subset Y$ and $A \cap H$ are $C_0$-quasiconvex for all
$H \in \mathcal{H}$, then $(A,d_{pel})$ is $C_1$-quasiconvex 
in $(X,  d_{pel})$.
\label{pel-qc}
\end{lemma}

\subsection{Amalgamated Geometry for Surfaces with Punctures}

\noindent {\bf Step 1:}
For a hyperbolic surface $S^h$ (possibly) with punctures, we fix a (small) 
$e_0$, and excise the cusps leaving horocyclic boundary components of
(ordinary or Euclidean) length 
$e_0$. We then take the induced {\em path metric} on $S^h$ minus cusps and call
the resulting surface $S$. This induced path metric will still be referred to
as the hyperbolic metric on $S$ (with the understanding that now $S$ possibly
has boundary). 

\noindent {\bf Step 2:}
The definitions and constructions of {\bf amalgamated building blocks} and 
{\bf amalgamation components} now go through with appropriate changes. The
only difference is that $S$ now might have boundary curves of length $e_0$.
For thick blocks, we assume (as in \cite{brahma-ibdd} ) that a thick block
is the universal curve over a Teichmuller geodesic (of length less than $D$
for some uniform $D$) minus {\em cusps $\times I$}. 

 There is one subtle point about global quasiconvexity (in
$\tilde{M}$) of
amalgamation components. This does not hold in the metric obtained by
merely excising the cusps and equipping the resulting horospheres with
the Euclidean metric. What we demand is that each amalgamation
component along with the parts of the horoballs that  meet the
boundary 
(horocycle times closed interval)'s be 
quasiconvex in $\tilde{M}$. When we partially electrocute horospheres
below, and consider quasiconvexity in the resulting partially
electrocuted space, amalgamation components in this sense remain
quasiconvex by Lemma \ref{pel-qc}. 

\noindent {\bf Step 3:}
Next, we modify the metric on $S$ 
by electrocuting its boundary components so that the metric
on the boundary components of each block $S \times I$ is the product
of the zero metric on the horocycles of fixed (Euclidean) length $e_0$
and the Euclidean metric on the $I$-factor.
The resulting blocks will be called {\bf partially electrocuted
blocks}. 
 We demand that in the model $M_{pel}$ obtained by gluing
together partially electrocuted blocks, the amalgamation components
are uniformly quasiconvex.  
By Lemma \ref{pel-qc}, this follows from quasiconvexity of
amalgamation components in the sense of the note above. Note that
$M_{pel}$ may also be constructed directly from $M$ by excising a
neighborhood of the cusps and partially electrocuting the resulting
horospheres. By Lemma \ref{pel} ${\tilde{M}}_{pel}$ is a hyperbolic
metric space. 

\noindent {\bf Step 4:}
 Again,
the definitions and constructions of {\bf amalgamated building blocks} and 
{\bf amalgamation components} now go through {\em mutatis mutandis} for 
partially electrocuted blocks. 

\noindent {\bf Step 5:} Next, let $\lambda^h$ be a hyperbolic
geodesic in $\tilde{S^h}$. We replace pieces of
$\lambda^h$ that lie within horodisks by shortest
horocyclic segments joining its entry
and exit points (into the corresponding horodisk). Such a path
is called a  horo-ambient quasigeodesic in \cite{brahma-pared}. See
Figure below:

\smallskip

\begin{center}

\includegraphics[height=4cm]{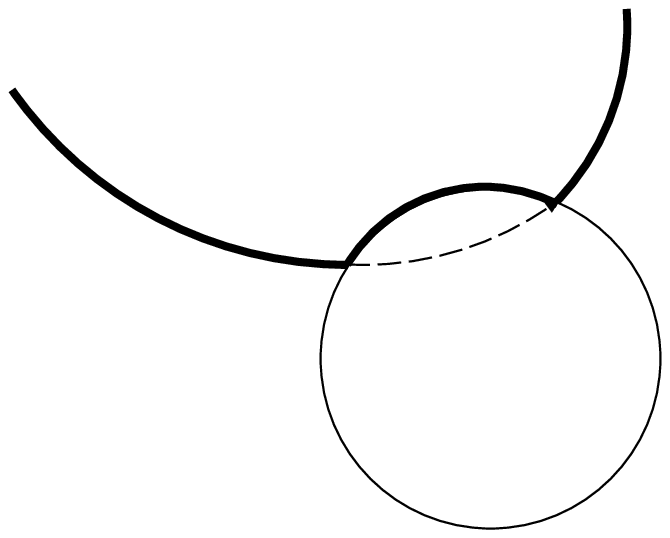}

\smallskip

\underline{Figure 4:{\it Horo-ambient quasigeodesic} }

\end{center}

\smallskip

A small modification might be introduced if we electrocute horocycles.
Geodesics and
quasigeodesics without backtracking
then travel for free along the zero metric horocycles.
This does not change matters much as the geodesics and
quasigeodesics in the two resulting constructions track each other by
Lemma \ref{farb2A}.

\noindent {\bf Step 6:}
Thus, our starting point for the construction of the hyperbolic ladder
$B_\lambda$ is not a hyperbolic geodesic $\lambda^h$ but a horoambient
quasigeodesic $\lambda$. We construct the {\bf graph model} as before.
{\it By Lemma \ref{pel-qc} quasiconvexity of amalgamation components as
well as lifts of Margulis tubes is preserved by partial
electrocution.}

\noindent {\bf Step 7:}
The construction of $B_\lambda , \Pi_\lambda$ and their properties
go through {\it mutatis mutandis}
 and we conclude that $B_\lambda$ is quasiconvex
in the {\em graph model } of the
partially electrocuted space ${\tilde{M}}_{pel}$. As before,
${\widetilde{M_H}}_{pel}$ will denote the collection of horizontal sheets.
The modification of Theorem \ref{retract} is given below:

\begin{theorem}
There exists $C > 0$ such that 
for any horo-ambient geodesic
 $\lambda = \lambda_0 \subset \widetilde{S} \times \{
0 \} \subset \widetilde{B_0}$, the retraction $\Pi_\lambda :
{\widetilde{M_H}}_{pel} \rightarrow B_\lambda $ satisfies: \\

 $d_{pel}( \Pi_{\lambda , B} (x), \Pi_{\lambda , B} (y)) \leq C d(x,y) +
C$.
\label{retract-punct}
\end{theorem}

\noindent {\bf Step 8:}
From this step on, the modifications for punctured
surfaces follow \cite{brahma-pared}
As in \cite{brahma-pared}, we decompose $\lambda$ into 
portions $\lambda^c$ and $\lambda^b$  that lie
along horocycles and those that do not.  Accordingly, we decompose
$B_\lambda$ into two
parts $B^c_\lambda$ and $B^b_\lambda$ consisiting of parts that lie
along horocycles and those that do not.  Dotted geodesics and
admissible
paths are constructed as before.
As in Lemma \ref{far-nopunct}, we get

\begin{lemma}
There exists a function $M(N) : \Bbb{N} \rightarrow \Bbb{N}$
such that  $M (N) \rightarrow \infty$
 as $N \rightarrow \infty$ for which the following holds:\\
For any horo-ambient quasigeodesic $\lambda \subset \widetilde{S}
\times \{ 0 \} \subset 
\widetilde{B_0}$, a fixed reference point 
 $p \in \widetilde{S} \times \{ 0 \} \subset
\widetilde{B_0}$ and any $x$ on  $B^b_\lambda$, 

\begin{center}

$d(\lambda^b , p) \geq N \Rightarrow d(x,p) \geq M(N)$.

\end{center}

\label{far-punct}
\end{lemma}

\noindent {\bf Step 9:}
 Construct a `dotted' ambient electric
quasigeodesic lying on $B_\lambda$ by projecting some(any) ambient
electric quasigeodesic onto $B_\lambda$ by $\Pi_\lambda$.
Join the dots using admissible paths to get
 a
connected ambient electric quasigeodesic $\beta_{amb}$.

\noindent {\bf Step 10} Construct from $\beta_{amb} \subset
 \widetilde{M}$ an electric 
quasigeodesic $\gamma$ in ${\widetilde{M}}_{pel}$ as in the previous section
and note that parts of $\gamma$ not lying along horocycles lie close
 to $B^b_\lambda$.

\noindent {\bf Step 11} Conclude that if $\lambda^h$ lies outside large
balls in $S^h$ then each point of  $\gamma$  lying outside partially 
electrocuted horospheres also lies outside
large balls. 
\\
\noindent {\bf Step 12} Let $\gamma^h$ denote 
the hyperbolic geodesic in ${\widetilde{M}}^h$ joining the end-points
of $\gamma$. By Lemma \ref{pel-track}  $\gamma$ and $\gamma^h$
track each other off  a bounded (hyperbolic) neighborhood of the
electrocuted horoballs. Recall that $X$ denotes ${\widetilde{M}}^h$
minus interiors of horoballs. Then, every point of
$\gamma^h \cap X$ 
 must lie close to
 some point of $\gamma$ lying outside partially 
electrocuted horospheres. Hence from Step (11), 
if $\lambda^h$ lies outside large
balls about $p$ in $S^h$ then   $\gamma^h \cap X$   also lies outside
large balls about $p$ in $X$. In particular,  $\gamma^h$ enters and
leaves horoballs at large distances from $p$. From this we conclude
that
$\gamma^h$ lies  outside large balls. Hence
by Lemma \ref{contlemma} there exists a Cannon-Thurston map and the
limit set is locally connected.
\\

We state the conclusion below:

\begin{theorem}
Let $M^h$ be a 3 manifold homeomorphic to $S^h \times J$ (for $J = [0,
  \infty ) $ or $( - \infty , \infty )$). Further suppose that $M^h$ has
  {\em amalgamated geometry}, where $S^h_0 \subset B_0$ is the lower
  horizontal surface of the building block $B_0$. Then the inclusion
  $i : {\widetilde{S}}^h \rightarrow {\widetilde{M}}^h$ extends continuously
  to a map 
  $\hat{i} : {\widehat{S}}^h \rightarrow {\widehat{M}}^h$. Hence the limit set
  of ${\widetilde{S}}^h$ is locally connected.
\label{crucial-punct}
\end{theorem}

\section{Weakening the Hypothesis I: Graph Quasiconvexity and Graph
  Amalgamation Geometry}

We now proceed to weaken the hypothesis
of amalgamation geometry in the hope of capturing all Kleinian surface
groups.  Recall that in the definition of
amalgamation geometry,  
 two criteria  were used - {\em local and global quasiconvexity of
   amalgamation components}. We shall retain local quasiconvexity, and
 replace global quasiconvexity by a {\bf weaker condition}
 which we shall term {\bf graph quasiconvexity}. The rationale behind
 this terminology shall be made clear later.
We first modify the definition of amalgamation geometry as follows,
retaining only local quasiconvexity.
We first recall the definition of {\it weak amalgamation geometry}.

\noindent A manifold $M$ homeormorphic to $S \times
  J$, where $ J = [0,
  {\infty })$ or $J = ( - \infty , \infty )$, is said to be a model of {\bf
weak  amalgamation geometry} if  \\

\begin{enumerate}

\item there is a fiber preserving homeomorphism from $M$ to
  $\widetilde{S} \times J$ 
that lifts to  a quasi-isometry of universal covers \\ 

\item there exists a sequence $I_i$ of intervals (with disjoint
 interiors)
 and blocks $B_i$
where the metric on $S \times I_i$ is the same as
 that on some  building block $B_i$. Each block is either thick or has
 amalgamation geometry. \\

\item $\bigcup_i I_i = J$ \\

\item There exists $ C > 0$ such that for all amalgamated blocks $B_i$
  and geometric cores $K \subset B_i$, all amalgamation components of
  $\widetilde{K}$ are $C$-quasiconvex in ${\widetilde{B}}_i$ \\

\end{enumerate}

\noindent {\bf Definition:} An amalgamation component $K \subset B_n$
is said to be
 (m. $\kappa$  )  {\bf graph - quasiconvex} if there exists a
$\kappa$-quasiconvex (in the hyperbolic metric)
 subset  $CH(K)$ containing $K$ such that 

\begin{enumerate}
\item $CH(K) \subset N_m^G(K)$ 
where $N_m^G(K)$ denotes the $m$ neighborhood of $K$ in the graph
model of $M$. \\
\item For each $K$ there exists $C_K$ such that
$K$ is $C_K$-quasiconvex in $CH(K)$.
\end{enumerate} 

\smallskip

Since the quasiconvex sets (thought of as convex hulls of $K$) lie
within a bounded distance from $K$ in the {\em graph model} we have
used the term {\em graph-quasiconvex}.

\noindent {\bf Definition:} A manifold $M$ of weak amalgamation
 geometry
is said to be a model of  {\bf
graph  amalgamation geometry} if there exist $m, \kappa$ such that
 each amalgamation geometry component is $(m, \kappa )$ -graph -
quasiconvex.

A manifold $N$ is said to have {\bf graph amalgamation geometry} if there
 is a level-preserving homeomorphism
 from $N$ to a model manifold of {\it graph amalgamation geometry}
that lifts to a quasi-isometry at the level of universal covers.

\smallskip

\noindent {\bf Note:} As before, we proceed with the assumption that
for surfaces with punctures, $S$ corresponds to a complete hyperbolic
surface $S^h$ minus a neighborhood of the cusps with horocycles
electrocuted. Further, $M$ corresponds to $M^h$ minus a neighborhood
of the cusps with resultant horospheres {\em partially electrocuted}.

\smallskip

Now, let us indicate the modifications necessary to carry out the
proof of the Cannon-Thurston Property for manifolds of graph
amalgamation geometry (suppressing  the quasi-isometry to a
model manifold). 
As in Theorem \ref{crucial}, the proof consists of two steps:

\begin{enumerate}
\item Constructing a quasiconvex set $B_\lambda$
in an auxiliary electric space
  (the {\bf graph model} ), and from this an admissible electric
  quasigeodesic $\beta$. \\
\item Recovering from $\beta$ and its intersection pattern,
  information about the hyperbolic geodesic joining its end-points.\\
\end{enumerate}

The first step is the same as that for models of {\em amalgamation
  geometry} as it goes through for {\em weak amalgamation geometry}. 
Then from Corollary \ref{retract-cor}  we have:

\noindent {\bf Step 1A:} Given $\lambda \subset {\widetilde{S}} \times
\{ 0 \}$, construct $B_\lambda$, $\Pi_\lambda$ as before.
There exists $C > 0$ such that the the retraction $\Pi_\lambda :
\widetilde{M_H} \rightarrow B_\lambda $ satisfies: \\
 $d_{Gel}( \Pi_{\lambda } (x), \Pi_{\lambda} (y)) \leq C d_{Gel}(x,y) +
C$, where $d_{Gel}$ denotes the metric in the graph model.

\medskip

Again, from  Corollary
\ref{adm-qgeod-props-cor} we have:

\noindent {\bf Step 1B:} \\
There exists a function $M(N) : \Bbb{N} \rightarrow \Bbb{N}$
such that  $M (N) \rightarrow \infty$
 as $N \rightarrow \infty$ for which the following holds:\\
For any geodesic $\lambda \subset \widetilde{S} \times \{ 0 \} \subset
\widetilde{B_0}$, and  a fixed reference point 
 $p \in \widetilde{S} \times \{ 0 \} \subset
\widetilde{B_0}$,
there exists a connected $B_\lambda$-admissible
 quasigeodesic $\beta_{adm}$  without 
backtracking, such that \\
$\bullet$ $\beta_{adm}$ is built up of $B_\lambda$-admissible
paths. \\
$\bullet$ $\beta_{adm}$ joins
the end-points of $\lambda$. \\
$\bullet$ 
$d(\lambda , p) \geq N \Rightarrow d(\beta_{adm},p) \geq M(N)$. ($d$
is the ordinary, non-electric metric.) \\

\medskip

\noindent {\bf Summary of Step 2:} \\
Now we come to the second step: {\bf recovering a hyperbolic geodesic
  from an electric geodesic}. 

\smallskip

This step can be further subdivided into two parts. In the
  first part  we construct
  a second auxiliary space $M_2$ 
by electrocuting
  the elements $CH(K)$. We show that the spaces ${\tilde{M}}_1$ and
${\tilde{M}}_2$ are quasi-isometric. In fact we show that the identity
  map on the underlying subset is a quasi-isometry. This step requires
  only the first condition in the definition of {\it graph
  quasiconvexity}. The second stage extracts information
  about an electro-ambient quasi-geodesic in ${\tilde{M}}_2$ from an
  admissible path in ${\tilde{M}}_1$. It is at this second stage that
  we require the second condition: (not necessarily uniform)
quasi-convexity of amalgamation components. 

We now furnish the details.

\medskip

\noindent {\bf Step 2A:} \\
Let $M_1$ denote $M$ with the graph metric
obtained by electrocuting amalgamation components. Next, let $M_2$
denote $M$ with an electric metric obtained by electrocuting the
family of sets $CH(K)$ (for amalgamation components $K$) appearing in
the definition of {\bf graph amalgamation geometry}. 
\begin{lemma}
 The identity map on the underlying set $M$ from
$M_1$ to $M_2$ induces a quasi-isometry of universal covers
${\widetilde{M}}_1$ and ${\widetilde{M}}_2$. 
\label{qi12}
\end{lemma}

\noindent {\bf Proof:} Let $d_1$, $d_2$ denote the electric
metrics on ${\widetilde{M}}_1$ and ${\widetilde{M}}_2$. Since $K
\subset CH(K)$ for every amalgamation component, we have right off \\
\begin{center}
$d_1 (x,y) \leq d_2 (x,y)$ for all $x, y \in \tilde{M}$
\end{center}

To prove a reverse inequality with appropriate constants, it is enough
to show that each set $CH(K)$ (of diameter one in $M_2$) has uniformly
bounded diameter in $M_1$. To see this, note that by definition of
graph-quasiconvexity, there exists $n$ such that for all $K$
and each point $a$ in $CH(K)$, there exists a point $b \in K$
with $d_1(x,y) \leq n$. Hence by the triangle inequality, \\
\begin{center}
$d_2 (x,y) \leq 2n+1$ for all $x, y \in \widetilde{CH(K)}$
\end{center}
  Therefore, \\
\begin{center}
$d_2 (x,y) \leq (2n+1)d_1 (x,y)$ for all $x, y \in \tilde{M}$
\end{center}

This proves the Lemma. $\Box$

\smallskip

\noindent {\bf Step 2B:} \\
Now let $\beta_{adm}$ denote an admissible $B_\lambda$
 quasigeodesic in
${\tilde{M}}_1 $, which does not backtrack relative to the
amalgamation components. By Lemma \ref{qi12} above, $\beta_{adm}$ is a
quasigeodesic in ${\tilde{M}}_2$. As in Lemma \ref{adm-ea-props},
using the Note following it, we conclude: \\

There exists a  $\kappa , \epsilon$-electro-ambient quasigeodesic
$\beta_{ea}$
in
${\widetilde{M}}_2$ (as opposed to ${\widetilde{M}}_1$, which is what
we needed in the {\em amalgamation geometry} case). (See Lemma
\ref{adm-ea-props}. ) Note that in ${\widetilde{M}}_2$, we electrocute
 the lifts of the sets $CH(K)$ rather than $\tilde{K}$'s. 

We thus obtain, as in Lemma \ref{adm-ea-props}
 a function $M^{\prime}(N) : \Bbb{N} \rightarrow \Bbb{N}$
such that  $M^{\prime}(N) \rightarrow \infty$
 as $N \rightarrow \infty$ for which the following holds:\\
For any geodesic $\lambda \subset \widetilde{S} \times \{ 0 \} \subset
\widetilde{B_0}$, and  a fixed reference point 
 $p \in \widetilde{S} \times \{ 0 \} \subset
\widetilde{B_0}$,
there exists a $(\kappa , \epsilon )$
  electro-ambient quasigeodesic $\beta_{ea}$  without 
backtracking, such that \\
$\bullet$ $\beta_{ea}$ joins
the end-points of $\lambda$. \\
$\bullet$ 
$d(\lambda , p) \geq N \Rightarrow d(\beta_{ea},p) \geq M^{\prime}(N)$. \\

\medskip

Finally, as in the proof of Theorem \ref{crucial}, we use Lemma
\ref{ea-strong} to conclude that the hyperbolic geodesic in
$\tilde{M}$
joining the
end-points of $\lambda$ lies in a uniform hyperbolic neighborhood of
$\beta_{ea}$. This gives us Theorem \ref{crucial} with {\bf graph
  amalgamation geometry} replacing {\em amalgamation geometry}.

\begin{theorem}
Let $M$ be a 3 manifold homeomorphic to $S \times J$ (for $J = [0,
  \infty ) $ or $( - \infty , \infty )$). Further suppose that $M$ has
  {\em graph amalgamation geometry}, where $S_0 \subset B_0$ is the lower
  horizontal surface of the building block $B_0$. Then the inclusion
  $i : \widetilde{S} \rightarrow \widetilde{M}$ extends continuously
  to a map 
  $\hat{i} : \widehat{S} \rightarrow \widehat{M}$. Hence the limit set
  of $\widetilde{S}$ is locally connected.
\label{crucial-graph}
\end{theorem}

The modifications for the case with punctures are as before (See
Theorem \ref{crucial-punct}. Thus, we
conclude:

\begin{theorem}
Let $M^h$ be a 3 manifold homeomorphic to $S^h \times J$ (for $J = [0,
  \infty ) $ or $( - \infty , \infty )$). Further suppose that $M^h$ has
  {\em amalgamated geometry}, where $S^h_0 \subset B_0$ is the lower
  horizontal surface of the building block $B_0$. Then the inclusion
  $i : {\widetilde{S}}^h \rightarrow {\widetilde{M}}^h$ extends continuously
  to a map 
  $\hat{i} : {\widehat{S}}^h \rightarrow {\widehat{M}}^h$. Hence the limit set
  of ${\widetilde{S}}^h$ is locally connected.
\label{crucial-punct-graph}
\end{theorem}

\section{Weakening the Hypothesis II: Split  Geometry}

In this section, we shall weaken the hypothesis of {\em graph
  amalgamation geometry} further to include the possibility of
  Margulis tubes cutting across the blocks $B_i$. 
But before we do this, let us indicate a straightforward
  generalisation
of {\em amalgamation geometry} or {\em graph amalgamation geometry}

\subsection{More Margulis Tubes in a Block}

A straightforward generalisation of Theorem \ref{crucial} (or 
Theorem \ref{crucial-graph}) is to
the case where more than one Margulis tube is allowed per block $B$, and
 each of these tubes  splits the block $B$ locally. On the
 surface $S$, this corresponds to a number of disjoint
(uniformly) bounded length curves.
 As before we require that each
amalgamation component be uniformly quasiconvex (or graph quasiconvex)
in $\tilde{M}$ for the
proof of Theorem \ref{crucial} (or Theorem \ref{crucial-graph})
to go through. See the figure below for
a schematic rendering of the model block of amalgamation geometry.

\smallskip

\begin{center}

\includegraphics[height=4cm]{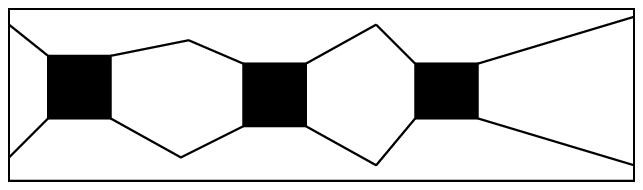}

\smallskip

\underline{Figure 5:{\it Building Block for
 Generalised Amalgamation Geometry} }

\end{center}

\smallskip

\subsection{Motivation for Split Geometry}

So far, we have assumed that the boundaries of {\em amalgamated
  geometry blocks} or
 {\em graph amalgamated
  geometry blocks} are all of bounded geometry. This assumption needs
  to be relaxed to accomodate general surface Kleinian groups. 
Before we define the objects of interest, we shall first informally
  analyse what went into the construction of the hyperbolic ladder
  $B_\lambda$. We require: \\

\begin{enumerate}
\item Horizontal surfaces $S_i$, all abstractly homeomorphic to each
  other \\
\item A block decomposition $M = \cup B_i$, where $B_{i-1} \cap B_{i}
  = S_i$ \\
\item Given a geodesic $\lambda_i \subset {\widetilde{S}}_i$, we require
 a (uniformly) large-scale retract $\pi_i$ of 
${\widetilde{S}}_i$ onto $\lambda_i$ and a prescription to construct 
$\lambda_{i+1} \subset {\widetilde{S}}_{i+1}$. Thus, starting with 
$\lambda_0 \subset {\widetilde{S}}_0$, we first construct $\pi_0$ and
 then inductively construct the pairs $(\lambda_i , \pi_i )$. \\
\item Each block $B_i$ has an auxiliary metric or pseudometric which
  induces the given path metrics on $S_{i-1}, S_i$. 
\end{enumerate}

We want to relax the assumption that $S_i$'s have bounded geometry,
while retaining the essential properties of bounded geometry. As
elsewhere in this paper 
we invoke the following 
(uncomfortably dictatorial) policy that we have adopted:

\smallskip

\noindent {\bf Policy:} {\it Electrocute anything that gives trouble.}

\smallskip

{\em What this policy means is that whenever some 
  construction possibly gives rise to non-uniformity of some
  parameter(s), locate the source of non-uniformity and electrocute
  it. Then, at the end of the game, re-instate the original geometry
  by using comparison properties between ordinary hyperbolic geometry
  and electric geometry}.

\smallskip

Thus, each $S_i$ is now allowed to have a pseudometric where a finite
number of disjoint, bounded length (uniformly, independent of $i$)
collection of simple closed geodesics are electrocuted. Then, instead
of geodesics $\lambda_i \subset {\widetilde{S}}_i$, we shall require
the $\lambda_i$ to be only electro-ambient geodesics. This will allow
us to go ahead with the construction of $B_\lambda$.

One further comment as to how this solves the problem. Let us fix a
small
(less than Margulis constant) $\epsilon_0$. Given any hyperbolic
surface $S^h$, we can simply electrocute {\it thin parts},
i.e. tubular neighborhoods of short (less than $\epsilon_0$) geodesics with
boundaries of length $\epsilon_0$. Alternately, we can first
cut out the
interiors
of these
thin parts. Next, corresponding to 
each {\em Margulis annulus} that has been cut out, 
glue the corresponding boundary components of length $\epsilon_0$
 together, and then electrocute the resulting closed curves.

This construction
 is adapted to the construction of {\em split  level surfaces} in
 Minsky \cite{minsky-elc1}, and Brock-Canary-Minsky
 \cite{minsky-elc2}. 

\subsection{Definitions}
  
Toplogically, a {\bf split subsurface} $S^s$ of a surface $S$ is a
(possibly disconnected, proper)
subsurface with boundary such that $S - S^s$ consists of a non-empty
family of
non-homotopic annulii, which in turn are not homotopic into the
boundary of $S^s$. \\

\smallskip

Geometrically, we assume that $S$ is given some finite volume
hyperbolic
structure. A split subsurface $S^s$ of $S$ has bounded geometry, i.e.\\
\begin{enumerate}
\item each boundary component of $S^s$ is of length $\epsilon_0$, and
  is in fact a component of the boundary of $N_k(\gamma )$, where
  $\gamma$ is a hyperbolic geodesic on $S$, and $N_k(\gamma )$ denotes
  its $k$-neighborhood. \\
\item For any  closed geodesic $\beta$ on $S$, either $\beta \subset S
  - S^s$, or, the length of any component of $\beta \cap (S - S^s )$
  is greater than $\epsilon_0$. \\
\end{enumerate}

\smallskip

Topologically, a {\bf split block} $B^s \subset B = S \times I$
is a topological product $S^s \times I$ for some {\em connected
  $S^s$}. However, its upper and lower boundaries need not be $S^s
\times 1$ and $S^s \times 0$. We only require that the upper and
lower boundaries
be split subsurfaces of $S^s$. This is to allow for Margulis tubes
starting (or ending) within the split block. Such tubes would split
one of the horizontal boundaries but not both. We shall call such
tubes {\bf hanging tubes}. See figure below: \\


\begin{center}

\includegraphics[height=4cm]{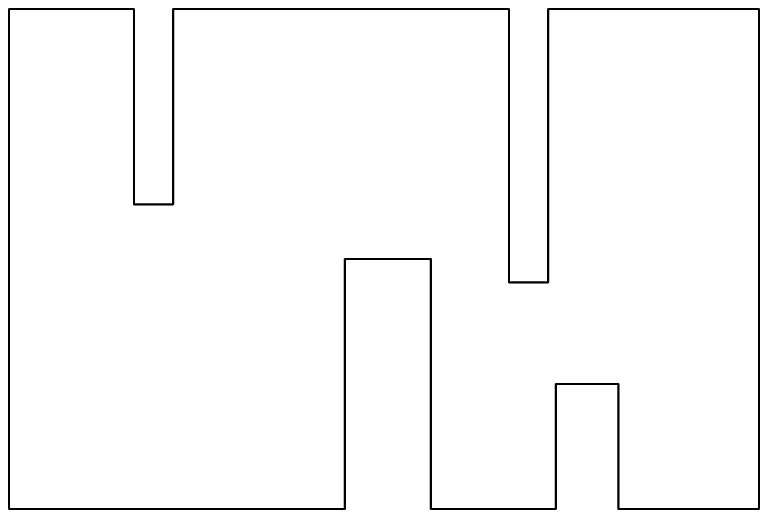}

\underline{Figure 6:  {\it Split Block with hanging tubes} }

\end{center}

\smallskip

Geometrically,
we require that the metric on  a split block
induces a path metric on its {\em upper and lower
horizontal} boundary components, which are subsurfaces of
$S^s \times \partial I$, such that each horizontal
boundary component is a
(geometric) split surface. Further, the metric on $B^s$
induces on each {\em vertical} boundary
component of a Margulis tube $\partial S^s \times I$
   the product metric. Each
boundary component for
Margulis tubes that `travel all the way from the lower to the upper
boundary' is an annulus of height
 equal to length of $I$. We demand further that {\em hanging tubes}
 have {\em  length uniformly bounded
  below by $\eta_0 > 0$}. Further, each such annulus
has cross section a round circle of length $\epsilon_0$.
This leaves us to decide the metric on
lower and upper boundaries of hanging tubes. Such
 boundaries are declared to have a metric equal to
that on $S^1 \times [-\eta, \eta]$, where $S^1$ is a round circle of
length $\epsilon_0$ and $\eta$ is a sufficiently small number.

\smallskip

{\bf Note:} In the above definition, we {\em do not} require that the
upper (or lower) horizontal boundary of a split block $B^s$ be
connected for a connected $B^s$. This happens due to the presence of
{\em hanging tubes}.

\smallskip

We further require that the distance between horizontal boundary
components
is at least $1$, i.e. for a component $R$ of $S^s$ $d(R \times 0, R
\times 1) \geq 1$. We define the {\bf thickness} of a split block to be the
supremum of the lengths of $x \times I$ for $x \in S^s$ and demand
that it be finite (which holds under all reasonable conditions,
e.g. a smooth metric; however, since we shall have occasion to deal
with possibly discontinuous pseudometrics, we make this explicit). We
shall denote the thickness of a split block $B^s$ by $l_B$.

Each component of a split block shall be called a {\bf split
  component}. We
further require that the `vertical boundaries' (corresponding to
  Euclidean annulii)  of split components be
  uniformly (independent of choice of a block and a split component)
quasiconvex in the corresponding split component.

Note that the boundary of each split block has an intrinsic metric
that is flat and corresponds to a Euclidean torus.

A lift of a split block to the universal cover of the block $B = S
\times I$ shall be termed a {\bf split component} of $\tilde{B}$.

\noindent {\bf Remark:} The notion of {\it split components} we
deal with here  is closely related to
the notion of {\bf bands} described by Bowditch in
\cite{bowditch-endinv}, \cite{bowditch-model} and also to the
notion of {\bf scaffolds} introduced by Brock, Canary and Minsky
in \cite{minsky-elc2}.

\smallskip

We define a {\bf welded split block} to be a split block with
identifications as follows: Components
of $\partial S^s \times 0$ are glued together if and only if
they correspond to the same geodesic in $S - S^s$. The same is done
for components
of $\partial S^s \times 1$. A simple closed curve that results from
such an identification shall be called a {\bf weld curve}. For hanging
tubes, we also weld  the boundary circles of their
{\em lower or upper boundaries} by simply collapsing $S^1 \times
[-\eta, \eta]$ to $S^1 \times \{ 0 \}$.

 This may be done topologically or
geometrically while retaining Dehn twist information about the
curves. To record information about the Dehn twists, we have to define
(topologically) a map that takes the lower boundary of a welded split
block to the upper boundary. We define a map that takes $x \times 0$
to $x \times 1$ for every point in $S^s$. This clearly induces a map
from the lower boundary of a welded split block to its upper
boundary. However, this is not enough to give a well-defined map on
paths. To do this, we have to record {\em twist information} about
{\em weld curves}.
The way to do this is to define a map on transversals to weld
curves. The map is defined on transversals by recording the number of
times a transversal to  a weld curve $\gamma \times 0$ twists around
$\gamma \times 1$ on  the upper boundary of the welded split block.
(A related context in which such transversal information is important
is that of markings described in Minsky \cite{minsky-elc1}.)

\smallskip

Let the metric product
$S^1 \times [0,1]$ be called the {\bf standard annulus}
 if each horizontal
$S^1$ has length $\epsilon_0$. For hanging tubes the standard annulus
 will be taken to be $S^1 \times [0,1/2]$.

\smallskip

Next, we require another pseudometric on $B$ which we shall term the
{\bf tube-electrocuted metric}. We first
 define a map  from
each boundary
annulus $S^1 \times I$ (or $S^1 \times [0,1/2]$ for hanging annulii)
to the corresponding standard annulus that is affine on the
second factor and an isometry on the first. Now glue the mapping
cylinder of this map to the boundary component. The resulting `split
block'
has a number of standard annulii as its boundary components. Call the
split block $B^s$ with the above mapping cylinders attached, the {\it
  stabilized split block} $B^{st}$.

Glue boundary components of $B^{st}$
corresponding to the same geodesic together
to get the {\bf tube electrocuted metric} on $B$ as follows.
Suppose that two boundary components of $B^{st}$ correspond to the
same geodesic $\gamma$. In this case, these boundary components are
both of the form $S^1 \times I$ or $S^1 \times [0, \frac{1}{2} ]$
where there is a  projection onto the horizontal $S^1$ factor
corresponding to $\gamma$. Let $S^1_l \times J$ and
$S^1_r \times J$ denote these two boundary components (where $J$
denotes $I$ or $[0, \frac{1}{2} ]$). Then each $S^1 \times \{ x \}$ has length
$\epsilon_0$.
Glue $S^1_l \times J$ to $S^1_r \times J$ by the natural `identity
map'.   Finally, on each resulting $S^1 \times \{ x \}$ put the zero
metric. Thus the annulus $S^1 \times J$ obtained via this
identification has the zero metric in the {\it horizontal direction}
$S^1 \times \{ x \}$ and the Euclidean metric in the {\it vertical
  direction} $J$. The resulting block will be called the {\bf
  tube-electrocuted block} $B_{tel}$ and the pseudometric on it will
be denoted as $d_{tel}$.   Note that $B_{tel}$ is
homeomorphic to $S \times I$. The operation of obtaining a {\em tube
  electrocuted block and metric}  $(B_{tel}, d_{tel})$
from a split block $B^s$ shall be called
{\em tube electrocution}.

\smallskip

Next, fix a hyperbolic structure on a Riemann surface $S$ and construct the
metric product $S \times \Bbb{R}$. Fix further a positive real
number $l_0$.

\begin{definition} An annulus $A$ will be said to be {\bf vertical}
if it is of the form $\sigma \times J$ for $\sigma$ a geodesic of length
less than $l_0$ on $S$ and
$J = [a,b]$
a closed sub-interval of $\Bbb{R}$. $J$ will be called
the {\bf vertical interval} for the vertical annulus $A$. \\
A disjoint  collection of annulii is said to be a
{\bf vertical system} of annulii if
each annulus in the collection is vertical. \\
\end{definition}

The above definition is based on a definition due to Bowditch
\cite{bowditch-endinv},\cite{bowditch-model}.

Suppose now that
$S \times \Bbb{R}$ is equipped with a vertical system $\mathcal{A}$
 of annulii.
We shall call  $z \in  \Bbb{R}$ a

\begin{enumerate}

\item a {\bf beginning level} if $z$ is the lower bound of a vertical
 interval
for some  annulus $A \in \mathcal{A}$. \\
\item an {\bf ending level}  if  $z$ is the lower bound of a vertical
 interval
for some  annulus $A \in \mathcal{A}$. \\
\item an {\bf intermediate level}   if  $z$ is an interior point
 of a vertical
 interval
for some  annulus $A \in \mathcal{A}$. \\

\end{enumerate}

In the figure below (where for convenience, all appropriate levels are
marked with integers), $2, 5, 11$ and $14$ are {\it beginning levels}, $4, 7,
13$ and $16$ are {\it ending levels},
$3, 6, 9, 12$ and $15$ are {\it intermediate
levels}.  We shall also allow Dehn twists to occur while
going
along the annulus. \\

\smallskip

\begin{center}

\includegraphics[height=5cm]{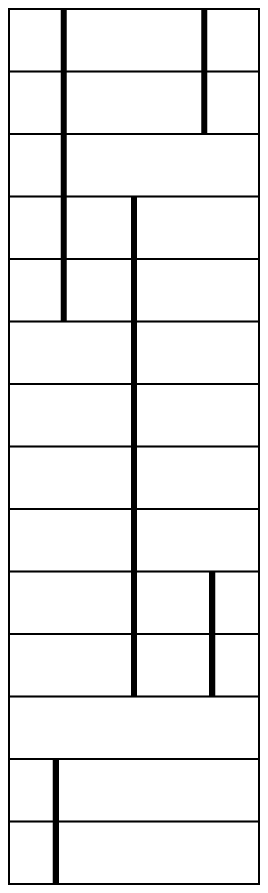}

\smallskip

\underline{Figure 7:  {\it Vertical Annulus Structure} }

\end{center}

\smallskip

A slight modification of the vertical annulus structure will sometimes
be useful.

Replacing each geodesic $\gamma$ on $S$ by a neighborhood $N_\epsilon
(\gamma )$ for sufficiently small $\epsilon$,
we obtain a {\bf vertical Margulis tube} structure after taking
products with vertical intervals. The family of Margulis tubes shall
be denoted by $\mathcal{T}$ and the union of their  interiors as $Int
\mathcal{T}$. The union of $Int \mathcal{T}$ and its  horizontal
boundaries (corresponding to neighborhoods of geodesics $\gamma
\subset S$ ) shall be denoted as $Int^+ \mathcal{T}$.

\smallskip

\noindent {\bf Thick Block}

Fix constants $D, \epsilon$ and let $\mu = [p,q]$ be  an $\epsilon$-thick
Teichmuller geodesic of length less than $D$. $\mu$ is
$\epsilon$-thick means that for any $x \in \mu$ and any closed
geodesic $\eta$ in the hyperbolic
surface $S_x$ over $x$, the length of $\eta$ is greater than
$\epsilon$.
Now let $B$ denote the universal curve over $\mu$ reparametrized
such that the length of $\mu$ is covered in unit time.
 Thus $B = S \times [0,1]$
topologically.

 $B$ is given the path
metric and is
called a {\bf thick building block}.

Note that after acting by an element of the mapping class group, we
might as well assume that $\mu$ lies in some given compact region of
Teichmuller space. This is because the marking on $S \times \{ 0 \}$
is not important, but rather its position relative to $S \times \{ 1 \}$
Further, since we shall be constructing models only upto
quasi-isometry, we might as well assume that $S \times \{ 0 \}$ and $S
\times \{ 1 \}$ {\em lie in the orbit} under the mapping class group
of some fixed base surface. Hence $\mu$ can be further simplified to
be a Teichmuller geodesic joining a pair $(p, q)$ amongst a finite set of
points in the orbit of a fixed hyperbolic surface $S$.

\smallskip

\noindent {\bf Weak Split Geometry}

A manifold $S \times \Bbb{R}$
equipped with a vertical Margulis tube
structure is said to be a model of
 {\bf weak split geometry}, if it is equipped with
a new metric satisfying the following conditions: \\
\begin{enumerate}
\item  $S \times [m,m+1] \cap Int \mathcal{T} = \emptyset$ (for $m \in
  \Bbb{Z} \subset \Bbb{R} ) $ implies that $S \times [m, m+1]$ is a
  thick block. \\
\item  $S \times [m,m+1] \cap Int \mathcal{T} \neq \emptyset$ (for $m \in
  \Bbb{Z} \subset \Bbb{R} ) $ implies that $S \times [m, m+1] - Int^+
  \mathcal{T} $ is (geometrically) a split block.  \\
\item There exists a uniform upper bound on the lengths of vertical
  intervals for vertical Margulis tubes \\
\item The metric on each component Margulis tube $T$ of $\mathcal{T}$
  is hyperbolic\\
\end{enumerate}

\noindent
{\bf Note 1:} Dehn twist information can still be implicitly recorded in a
model of {\em weak split geometry} by the Dehn filling information
corresponding to tubes $T$. \\
{\bf Note 2:} The metric on a model of {\em weak split geometry} is
possibly discontinuous along the boundary torii of Margulis tubes. If
necessary, one could smooth this out. But we would like to carry on
with
the above metric. \\

\smallskip

Removing the interiors of Margulis tubes and tube electrocuting each
block, we obtain a new pseudo-metric on ${M}$ called the {\bf
  tube electrocuted metric } $d_{tel}$ on ${M}$. The pseudometric
$d_{tel}$ may also be lifted to $\tilde{M}$.

The induced pseudometric on ${\tilde{S}}_i$'s shall be referred to as
{\bf split electric metrics}. The notions of {\em electro-ambient
  metrics, geodesics and quasigeodesics} 
go through in this context.

\smallskip

Next, we shall describe a {\bf graph metric} on $\tilde{M}$ which is
almost (but not quite)
 the metric on the nerve of the covering of $\tilde{M}$ by split
components (where each edge is assigned length $1$). This is not
strictly true as thick blocks are retained with their usual geometry
in the graph metric. However the analogy with the nerve is exact if
all blocks have {\em weak split geometry}.

\smallskip

For each split component $\tilde{K}$ assign a single vertex $v_K$ and
construct a cone of height $1/2$ with base $\tilde{K}$ and vertex
$v_K$. The metric on
the resulting space (coned-off or electric space in the sense
of Farb \cite{farb-relhyp}) shall be called the {\bf graph
metric} on $\tilde{M}$.

\smallskip

The union of a  split component of $\tilde{B}$ and the lifts of
Margulis tubes (to $\tilde{M}$) that intersect its boundary shall be
called a {\bf split amalgamation component} in $\tilde{M}$.

\noindent {\bf Definition:} A split amalgamation component $K$
is said to be
{\bf (m. $\kappa$  ) -  graph  quasiconvex} if there exists a
$\kappa$-quasiconvex (in the hyperbolic metric)
 subset  $CH(K)$ containing $K$ such that

\begin{enumerate}
\item $CH(K) \subset N_m^G(K)$
where $N_m^G(K)$ denotes the $m$ neighborhood of $K$ in the graph
metric on  $M$. \\
\item For each $K$ there exists $C_K$ such that
$K$ is $C_K$-quasiconvex in $CH(K)$.
\end{enumerate}

\noindent {\bf Definition:} A model manifold $M$ of weak split
 geometry
is said to be a model of  {\bf
split geometry} if there exist $m, \kappa$ such that
 each split amalgamation  component is $(m, \kappa )$ - graph
quasiconvex.

\subsection{The Cannon-Thurston Property for Manifolds of Split
 Geometry}

We shall first extract information about geodesics in the {\em tube
  electrocuted} model. As with Theorem \ref{crucial} and Theorem
  \ref{crucial-graph}, the proof splits into two parts: \\

\begin{enumerate}
\item Construction of $B_\lambda$ and its quasiconvexity in an
  auxiliary graph metric. The end-product of this step is an
  electro-ambient quasigeodesic in the graph model \\
\item Extraction of information about a hyperbolic geodesic and its
  intersection pattern with blocks from the electro-ambient
  quasigeodesic constructed in Step 1 above. 
\end{enumerate}

\noindent {\bf Details of Step 1:}

\smallskip

\noindent {\bf Step 1A: Construction of $B_\lambda$} \\
It is at this stage that the construction differs somewhat from 
the constrcution of $B_\lambda$ for manifolds of {\em graph
  amalgamated geometry}.

We start with the (tube-electrocuted) metric $d_{tel}$ on the model
manifold of {\em split geometry}. Then there exists a sequence of
split surfaces $S_i$ exiting the end(s). 

{\em Recall} that in the construction of $B_\lambda$ (for all
preceding cases) we are {\em not} interested in the metric on each
${\tilde{S}}_i$ per se, but in geodesics on ${\tilde{S}}_i$.

The metric $d_{tel}$ on the model manifold induces the {\bf split
  electric metric} on each $S_i$ obtained by electrocuting the {\bf
  weld curves}. The natural geodesics to consider on ${\tilde{S}}_i$
  are therefore the {\em electro-ambient } quasigedoescis where the
  electrocuted subsets correspond to geodesics representing the weld
  curves. 

Thus we start off with a hyperbolic geodesic $\lambda$ in
${\tilde{S}}_0$ joining $a, b$ say. We let $\lambda_0$ denote the
electro-ambient quasigeodesic joining $a, b$ in the split electric
metric on ${\tilde{S}}_0$. Now construct $B_\lambda$ inductively as
follows:

\smallskip

\noindent $\bullet$ Each split block $B_i$ and hence ${\tilde{B}}_i$
comes equipped with a (topological) product structure. Thus there is a
canonical map $\Phi_i : {\tilde{S}}_i \rightarrow {\tilde{S}}_{i+1}$
which maps each $(x,i)$ to a point $(x,i+1)$ by lifting the map from
$S_i$ to $S_{i+1}$ ($i \geq 0$ corresponding to the product
structure).  \\
\noindent $\bullet$ Next, if $\lambda_i$ 
is an electro-ambient quasi-geodesic in the split electric metric on
${\tilde{S}}_i$ joining $(a,i)$ and $(b,i)$ we let $\lambda_{i+1}$
denote the electro-ambient quasigeodesic in the split-electric metric
on ${\tilde{S}}_{i+1}$ joining $(a,i+1)$ and $(b,i+1)$. This gives us
a prescription for constructing $\lambda_{i+1}$ from $\lambda_i$ for
$i \geq 0$. Similarly, for $i \leq 0$ (in the totally degenerate case)
we can construct $\lambda_{i-1}$ from $\lambda_i$. Then as before,
define \\
\begin{center}
$B_\lambda = \bigcup_i \lambda_i$
\end{center}

\noindent $\bullet$ Again, $\pi_i : {\tilde{S}}_i \rightarrow
\lambda_i$ is defined as the retarction that minimises the ordered
pair of distances in the split electric metric and the hyperbolic
metric (without electrocuting weld curves). $\Pi_\lambda$ is obtained
in the graph metric by defining it on the horizontal sheets
${\tilde{S}}_i$ as \\
\begin{center}
$\Pi_\lambda (x) = \pi_i (x)$ for $x \in {\tilde{S}}_i$. 
\end{center}

\noindent $\bullet$ 
Then as before we conclude that in the graph model for $\tilde{M}$,
with the metric $d_{Gel}$, $\Pi_\lambda$ does not stretch distances
much, i.e. there exists a uniform $C \geq 0$ such that \\
\begin{center}
$d_{Gel} ( \Pi_\lambda (x), \Pi_\lambda (y) ) \leq C d_{Gel} (x,y) +
  C$
\end{center}

\smallskip

\noindent {\bf Step 1B: Construction of admissible quasigedoesic} \\
The above construction of $\Pi_\lambda$ may be used to construct a
  $B_\lambda$-
  admissible quasigeodeic $\beta_{adm}$ in the tube-electrocuted
  model. As before we have: \\
There exists a function $M(N) : \Bbb{N} \rightarrow \Bbb{N}$
such that  $M (N) \rightarrow \infty$
 as $N \rightarrow \infty$ for which the following holds:\\
For any geodesic $\lambda \subset \widetilde{S} \times \{ 0 \} \subset
\widetilde{B_0}$, and  a fixed reference point 
 $p \in \widetilde{S} \times \{ 0 \} \subset
\widetilde{B_0}$,
there exists a connected $B_\lambda$-admissible
 quasigeodesic $\beta_{adm}$  without 
backtracking, such that \\
$\bullet$ $\beta_{adm}$ is built up of $B_\lambda$-admissible
paths. \\
$\bullet$ $\beta_{adm}$ joins
the end-points of $\lambda$. \\
$\bullet$ 
If $d(\lambda , p) \geq N$ then for any $x \in \beta_{adm} - Int
  \mathcal{T}$, $ d(x,p) \geq M(N)$. ($d$
is the ordinary, hyperbolic, or non-electric metric.) \\

\smallskip

\noindent {\bf Step 2: Recovering a  quasigeodesic in the tube electrocuted
  model from an
  admissible quasigeodesic} \\ 
We now follow the proof of Theorem \ref{crucial-graph}. 

\smallskip

\noindent {\bf Step 2A:}
As in {\it Step 2A} in the proof of Theorem \ref{crucial-graph}
we construct
  a second auxiliary space $M_2$ 
by electrocuting
  the elements $CH(K)$ for split components $K$. 
The spaces ${\tilde{M}}_1$ and
${\tilde{M}}_2$ are quasi-isometric by uniform
{\it graph quasiconvexity} of split components. 
In fact  the identity
  map on the underlying subset is a quasi-isometry as in Lemma
  \ref{qi12}. 

\smallskip

\noindent {\bf Step 2B}
Next, as in {\it Step 2B} in the proof of Theorem \ref{crucial-graph},
we extract information
  about an electro-ambient quasi-geodesic in ${\tilde{M}}_2$ from an
  admissible path in ${\tilde{M}}_1$. It is at this second stage that
  we require the  condition that split components are (not necessarily
  uniformly) 
quasi-convex in the hyperbolic metric, and hence by Lemma \ref{pel-qc}
in the tube electrocuted metric $d_{tel}$. 

\smallskip

We may assume that $\beta_{adm}$  does not backtrack relative to the
split components. From Step 2A above, $\beta_{adm}$ is a
quasigeodesic in ${\tilde{M}}_2$. Then we conclude: \\

There exists a  $\kappa , \epsilon$-electro-ambient quasigeodesic
$\beta_{tea}$ 
in
${\widetilde{M}}_2$  (Note that in ${\widetilde{M}}_2$, we electrocute
 the lifts of the sets $CH(K)$ rather than $\tilde{K}$'s). 

We finally obtain
 a function $M^{\prime}(N) : \Bbb{N} \rightarrow \Bbb{N}$
such that  $M^{\prime}(N) \rightarrow \infty$
 as $N \rightarrow \infty$ for which the following holds:\\
For any geodesic $\lambda \subset \widetilde{S} \times \{ 0 \} \subset
\widetilde{B_0}$, and  a fixed reference point 
 $p \in \widetilde{S} \times \{ 0 \} \subset
\widetilde{B_0}$,
there exists a $(\kappa , \epsilon )$
  electro-ambient quasigeodesic $\beta_{tea}$ (in the {\bf tube
 electrocuted metric})  without 
backtracking, such that \\
$\bullet$ $\beta_{tea}$ joins
the end-points of $\lambda$. \\
$\bullet$ 
If $\lambda$ lies outside a large ball about a fixed reference point
 $p \in {\tilde{S}}_0$, then
each point of 
$\beta_{tea} \cap ({\tilde{M}} - Int \mathcal{T} )$ also lies outside a
 large ball about $p$. \\

\noindent {\bf Step 3: Recovering a  hyperbolic geodesic from
the tube electrocuted quasigeodesic $\beta_{tea}$} \\
This is a new step that comes from the extra phenomenon of tube
electrocution which makes the metric $d_{tel}$ an `intermediate'
metric between the hyperbolic metric $d$ and the graph metric $d_{Gel}$.

\smallskip

Observe that lifts of  Margulis tubes to $( \tilde{M} , d_{Gel} )$
 have uniformly bounded diameter in the metric $d_{Gel}$ and consequently in
 the metric $d_{tel}$ by uniform boundedness of vertical
 intervals of vertical Margulis tubes. Hence the {\it tube electrocuted
 metric} $d_{tel}$
 on $\tilde{M}$ is quasi-isometric to the metric $d_{fe}$
where lifts of
 Margulis tubes are electrocuted (i.e. fully electrocuted rather than
 just tube electrocuted, and hence
 each tube has diameter $1$). Let ${\tilde{M}}_{fe}$
 denote $\tilde{M}$ equipped with this new metric. Then
 geodesics without backtracking in the tube electrocuted metric become
(uniform)  quasi-geodesics without backtracking in ${\tilde{M}}_{fe}$.

\smallskip

\noindent {\bf Note:} It is at this (rather late)
stage that we need to assume that $\tilde{M}$ is a hyperbolic metric space.

\smallskip

Let $\gamma^h$ denote a hyperbolic geodesic joining the end-points of
$\beta_{tea}$ and hence $\lambda$. By Lemma \ref{farb2A}, $\gamma^h$
and $\beta_{tea}$ track each other off Margulis tubes. Hence
$\gamma^h  \cap ({\tilde{M}} - Int \mathcal{T} )$ lies outside a large
ball about $p$.
In particular, this is true for 
entry
and exit points of $\gamma^h$ with respect to Margulis tubes.
This implies (See for instance Lemma 7.3 of \cite{brahma-ibdd} )
that the parts of $\lambda^h$ lying within Margulis tubes
also  lie outside large balls about $p$. As before, by Lemma
\ref{contlemma} we infer the Cannon-Thurston property for manifolds of
    {\em split geometry}.

\begin{theorem}
Let $M$ be a 3 manifold homeomorphic to $S \times J$ (for $J = [0,
  \infty ) $ or $( - \infty , \infty )$). Further suppose that $M$ has
  {\em split geometry}, where $S_0 \subset B_0$ is the lower
  horizontal surface of the building block $B_0$. Then the inclusion
  $i : \widetilde{S} \rightarrow \widetilde{M}$ extends continuously
  to a map 
  $\hat{i} : \widehat{S} \rightarrow \widehat{M}$. Hence the limit set
  of $\widetilde{S}$ is locally connected.
\label{crucial-split}
\end{theorem}

There is a bit of ineffective ambiguity in the above theorem. In split
geometry, $S_0$ is only a split surface. We can extend this to any
surface $S_0$ so long as the annulii that we glue on to construct the
full surface lie entirely within Margulis tubes.
The modifications for the case with punctures are as before:
conclude:

\begin{theorem}
Let $M^h$ be a 3 manifold homeomorphic to $S^h \times J$ (for $J = [0,
  \infty ) $ or $( - \infty , \infty )$). Further suppose that $M^h$ has
  {\em split geometry}, where $S^h_0 \subset B_0$ is the lower
  horizontal surface of the building block $B_0$. Then the inclusion
  $i : {\widetilde{S}}^h \rightarrow {\widetilde{M}}^h$ extends continuously
  to a map 
  $\hat{i} : {\widehat{S}}^h \rightarrow {\widehat{M}}^h$. Hence the limit set
  of ${\widetilde{S}}^h$ is locally connected.
\label{crucial-punct-split}
\end{theorem}

\section{Generalisation: Incompressible away from Cusps}

The aim of this section is to sketch the proof of the following more
general theorem:

\begin{theorem}
Suppose that $N^h \in H(M,P)$ is a hyperbolic structure of {\em
split geometry} 
on a pared manifold $(M,P)$ with incompressible boundary $\partial_0 M$. Let
$M_{gf}$ denote a geometrically finite hyperbolic structure adapted
to $(M,P)$. Then the map  $i: \widetilde{M_{gf}}
\rightarrow \widetilde{N^h}$ extends continuously to the boundary
$\hat{i}: \widehat{M_{gf}}
\rightarrow \widehat{N^h}$. If $\Lambda$ denotes the limit set of
$\widetilde{M}$, then $\Lambda$ is locally connected.
\label{main3}
\end{theorem}

See \cite{brahma-pared} for definition of pared manifold with
incompressible boundary (this coincides with the notion of manifolds
whose boundary is incompressible away from cusps). 
Theorem \ref{crucial-punct} and its proof takes the place of Theorem
4.15 of \cite{brahma-pared}. Since nothing else is new, given these
constituents, we content ourselves with giving an outline of the
proof.

\medskip

{\bf Outline of Proof of Theorem \ref{main3} }

\smallskip

\noindent {\bf Step 1} Construct $B_\lambda$ in $\widetilde{M}$ ( =
$\widetilde{M^h}$ - cusps) as in Section 4.1 of
\cite{brahma-pared}. The only difference is that for an end $E$ of
{\it split geometry}, $\widetilde{E}$ is given the graph metric
corresponding to the graph model.
\\
\noindent {\bf Step 2} As in Sections 4.2, 4.3 of \cite{brahma-pared}
we obtain a retract $\Pi_\lambda$ onto $B_\lambda$.
\\
\noindent {\bf Step 3} Construct a `dotted' ambient electric
quasigeodesic lying on $B_\lambda$ by projecting some(any) (graph)
geodesic onto $B_\lambda$ by $\Pi_\lambda$.
\\
\noindent {\bf Step 4} Join the dots using admissible paths.
 This results in a
connected ambient electric quasigeodesic $\beta_{amb}$.
\\
\noindent {\bf Step 4A} Construct from $\beta_{amb}$ an
electro-ambient
quasigeodesic $\beta_{ea}$ replacing bits that lie within blocks by
hyperbolic geodesics (which lie within a bounded distance from it in
the graph metric, by graph-quasiconvexity). \\
\noindent {\bf Step 5} Conclude that the segments of $\beta_{ea}$ that
lie off {\em partially electrocuted horospheres} in fact lie outside a
large ball about a fixed reference point if $\lambda^h$ (the
hyperbolic geodesic joining the end-points of $\lambda$ in
${\widetilde{S}}^h$ does so.
\noindent {\bf Step 6} Construct from $\beta_{ea} \subset
 \widetilde{M}$ an electro-ambient 
quasigeodesic $\gamma$ in $\widetilde{M^h}$ by
 replacing pieces of $\beta_{ea}$ that lie along partially
 electrocuted
 horoballs (if any) by 
hyperbolic quasigeodesics that lie (apart from bounded length segments
at the beginning and end) within horoballs. 
\\
\noindent {\bf Step 7} Conclude that if $\lambda^h$ lies outside large
balls in $S^h$ then each point of the path $\gamma$ also lies outside
large balls. 
\\
\noindent {\bf Step 8} Let $\gamma^h$ be the hyperbolic geodesic
joining the end-points of $\gamma$.
Since the underlying set of $\gamma^h$ lies in a neighborhood of
$\gamma$,
 by Lemma \ref{ea-strong}, it must  lie  outside large balls. Hence
by Lemma \ref{contlemma} there exists a Cannon-Thurston map and the
limit set is locally connected.
\\
\noindent {\bf Step 8}
As in \cite{brahma-pared}, the Steps 1-7 above are carried out first
for manifolds of {\it p-incompressible boundary}. Then in the last
step (as in Section 5.4  of \cite{brahma-pared}) the hypothesis is
relaxed and the result proven for pared manifolds with incompressible
boundary. (Recall from \cite{brahma-pared} that p-incompressibility
roughly means the absence of accidental parabolics in any hyperbolic
structure.) Note also that the definition of pared manifolds with
incompressible boundary coincides with the notion of
`incompressibility away from cusps' introduced by Brock, Canary and
Minsky in \cite{minsky-elc2}.

\section{The Minsky Model and Split Geometry: A Sketch}

The aim of this section is to sketch a 
proof of the following theorem:

\smallskip

\noindent {\bf Theorem \cite{mahan-split}:}
Let $M$ be a hyperbolic manifold corresponding to a totally degenerate
surface group. Then $M$ has split geometry.

\smallskip

We shall use a model manifold that was built by Minsky
in \cite{minsky-elc1} to prove the Ending Lamination Conjecture. It
was shown by Brock, Canary and Minsky in \cite{minsky-elc2} that the
model is bi-Lipschitz equivalent to a hyperbolic manifold with the
same ending laminations.

We refer the reader to Minsky \cite{minsky-elc1} for the definitions
of the relevant terms, particularly {\bf hierarchy,
resolution} and other related notions. Fix a hyperbolic surface $S$.

\smallskip

\noindent {\bf Step 1: {\underline Constructing a sequence of split 
surfaces}} \\

\medskip

We  require the following:

\begin{enumerate}

\item {\bf Resolution sweep} - Lemma 5.8 of \cite{minsky-elc1} \\
\item {\bf J(v) is an interval} - Lemma 5.16 of
  \cite{minsky-elc1}. What this means is the following: \\
Given a vertex $v$ (corresponding to a simple closed curve on the
  surface) occurring in the hierarchy $H$ obtained from the ending
  laminations, fix a resolution $\{ \tau_i \}_{i \in {\mathcal{I}}}$
  of $H$ with ${\mathcal{I}}$ a subinterval of $\Bbb{Z}$. In the
  doubly (resp. simply) degenerate case $\mathcal{I}$ can be thought
  of
as $\Bbb{Z}$ (resp. $\Bbb{N}$ ). Let \\
\begin{center}
$J(v) = \{ i \in \mathcal{I} :$ $ v \in $ {\bf base} $\mu_{\tau_i}$ $\}$
\end{center}
where {\bf base} $\mu_{\tau_i}$ denotes the pants decomposition
induced by the marking  $\mu_{\tau_i}$. We might as well assume that
there are no repetitions in $J(v)$ (see the proof of Theorem 8.1 in
\cite{minsky-elc1}). Then $J(v) $ is an interval. \\
\item Again from the proof of Theorem 8.1 of \cite{minsky-elc1}
  we obtain a flat orientation preserving embedding of the Minsky model
  minus Margulis tubes (denoted as $M_{\nu} (0)$)
into $S \times \Bbb{R}$.\\
\item To each $\tau_i$ Minsky associates a {\it split-level surface}
  $F_i$.
{\bf This is the point at which the notions we have introduced in this
  paper and its predecessor \cite{brahma-ibdd} converge with those in
  Minsky's construction of his model in \cite{minsky-elc1}. }

In fact, the term {\it split geometry} was chosen with this in view.
In what follows, we shall construct {\it split surfaces} (as per our 
definitions) from the {\it split level surfaces} of Minsky.\\

\end{enumerate}

From the Minsky model we shall construct:\\

\begin{enumerate}
\item A sequence of split surfaces $S^s_i$ exiting the
end(s) of $M$. These will determine the levels for the split blocks
and split geometry. There is a lower bound on the distance between
$S^s_i$ and $S^s_{i+1}$ \\
\item A collection of Margulis tubes $\mathcal{T}$. \\
\item For each  complementary annulus of $S^s_i$ with core $\sigma$,
  there is a Margulis tube $T$ whose core is freely homotopic to $\sigma$
  and  such that $T$ intersects the level $i$. (What this roughly
  means is that there is a $T$ that contains the complementary
  annulus.) \\
\item For all $i$, either there exists a Margulis tube splitting both $S^s_i$
  and $S^s_{i+1}$ and hence $B^s_i$, or else $B_i$ is a thick block. \\
\item  $T \cap S^s_i$ is either empty or consists of a pair of
  boundary components of $S^s_i$ that are parallel in $S_i$. \\
\item There is a uniform upper bound $n$ on the number of surfaces that
 $T$ splits. \\
\end{enumerate}

We  define $S_{i}^{s}$  to be the first
 split level surface in which $v_i$ occurs. The region between $S^s_i$
 and $S^s_{i+1}$ is temporarily deignated $B^s_i$. We shall describe
 in \cite{mahan-split} a
 procedure for interpolating auniformly bounded number of
 split surfaces between $S^s_i$ and
 $S^s_{i+1}$. .

It will be shown in \cite{mahan-split} that 
there exists $n$ such that each thin Margulis tube
splits at most $n$ split surfaces in the above sequence.

This allows us to conclude that
the Minsky model has weak aplit geometry.

\medskip

\noindent {\bf Step 2: {\underline Graph quasiconvexity of Split Components}}

\smallskip

In order to prove that the Minsky model enjoys the property of {\bf
  split geometry}, we need to show further that any of the split components is
  (not necessarily uniformly) quasiconvex in the hyperbolic metric,
  and uniformly quasiconvex in the graph metric, i.e. we require to
  show {\em hyperbolic quasiconvexity} and {\em uniform graph
  quasiconvexity} of split components. 

\smallskip

\noindent {\bf Step 2A:}  Hyperbolic 
quasiconvexity is easy to prove and follows from the Thurston-Canary
covering Theorem \cite{Thurstonnotes} \cite{canary-cover}. 
\smallskip

\noindent {\bf Step 2B:}
Next, we need to prove that each split component of $S^s_i$
corresponding to some subsurface $\Sigma$ of $S$ is
uniformly graph quasiconvex. First off, 
any simple closed curve  in $\Sigma$
 must be realised within a
uniformly bounded distance in the graph metric. To prove this, we show
in \cite{mahan-split} 
that any pleated surface which contains at least one boundary geodesic
of $\Sigma$ in its pleating locus is realised within a 
uniformly bounded distance of $S^s_i$ in the graph metric.

Next. any split component is bounded by Margulis tubes. We drill out
these tubes and appeal to the Drilling Theorem
\cite{brock-bromberg-density} to conclude that the drilled manifold
and the complement of the Margulis tube in the original manifold are
both uniformly bi-Lipschitz to the corresponding hyperbolic manifolds.
Now in the drilled manifold the subsurface $\Sigma$ gives us a genuine
quasifuchsian group, whose convex hull boundary is pleated and hence
within a uniform distance in the graph metric from the split
component. 

But the convex hull $CH_\Sigma$  of a lift $\tilde{\Sigma }$
in the drilled hyperbolic manifold may also be regarded as a
quasiconvex set in the hyperbolic manifold corresponding to the
surface group. (This requires some additional argument which is
supplied in \cite{mahan-split}.)

Since (using this identification) $CH_\Sigma$ is uniformly graph
quasiconvex in the drilled manifold, it is also uniformly graph
quasiconvex  in the split geometry model for
$\tilde{M}$. 

This shows that the Minsky model is of split geometry. Combining this
fact with Theorems \ref{crucial-split} and \ref{crucial-punct-split}
we shall obtain:

\smallskip

\noindent {\bf Theorem: \cite{mahan-split}}
Let $\rho$ be a representation of a surface group $H$ (corresponding
to the surface $S$) into
$PSl_2(C)$ without accidental parabolics. Let $M$ denote the (convex
core of) ${\Bbb{H}}^3 / \rho 
(H)$.  Further suppose that $i: S \rightarrow M$, taking
 parabolic to parabolics, induces a homotopy
equivalence.
  Then the inclusion
  $\tilde{i} : \widetilde{S} \rightarrow \widetilde{M}$ extends continuously
  to a map 
  $\hat{i} : \widehat{S} \rightarrow \widehat{M}$. Hence the limit set
  of $\widetilde{S}$ is locally connected.

\section{Extending the Sullivan-McMullen Dictionary }

A celebrated theorem of Yoccoz in Complex Dynamics (see Hubbard
\cite{hubbard-yoccoz}, or Milnor
\cite{milnor-yoccoz}) proves the local connectivity of certain Julia
sets using a technique called `puzzle pieces'. We shall not describe
this in any detail. What we shall simply say is that it consists of a
decomposition of a complex domain into pieces each of which under
iteration by a quadratic map converges to a single point. The
dynamical system can then be regarded as a semigroup ${\Bbb{Z}}_+$
of transformations acting on a complex domain. 

In the case of split (or amalgamation) geometry each of the split (or
amalgamation) components can be regarded as a 3-dimensional analogue
of puzzle pieces. Let us try to justify this analogy.
Suppose there is a group $G$ acting on the manifold
$\tilde{M}$. Let $H \subset G$ denote the fundamental group of a split
component. Let $G/H$ denote the coset space. 
Then what we require first is that if one takes a sequence of elements
$g_i$ going to infinity in the
coset space, the iterates of the split component converge to a point
in the limit sphere. However, this does not give all the information
as $G$ does not act co-compactly on $\tilde{M}$. In the cases we are
interested in $G/H$ correspond to normal directions to the split
component lying within the block containing the split component. This
does not help. 
To compensate, we look at the graph
model. Here, there is no group in sight. However, normal directions can
be salvaged from the {\bf graph metric}. Thus, instead of going to
infinity by iteration, we go to infinity in the graph metric. Further,
the analogue of the  requirement that iterates go to infinity, is that
the visual diameter goes to zero as we move to infinity in the graph
metric. This is easily ensured by hyperbolic quasiconvexity, and also
follows easily from {\bf graph quasiconvexity}. Note that {\bf graph
quasiconvexity} is a statement that gives uniform shrinking of visual
diameter to zero as one goes to infinity.

Thus we extend the Sullivan-McMullen dictionary (see
\cite{sullivan-dict}, \cite{ctm-renorm})between Kleinian
groups and complex dynamics by suggesting the following analogy:

\begin{enumerate}

\item {\em Puzzle pieces} are analogous to {\bf split components} \\
\item {\em Convergence  to a point under iteration} is analogous to
  {\bf graph quasiconvexity} \\
\end{enumerate}

One issue that gets clarified by the above analogy is a point raised
by McMullen in \cite{ctm-locconn}. McMullen indicates that though the
Julia set $J(P_\theta )$, where \\

\begin{center}

$P_\theta (z) = e^{2\pi i \theta}z + z^2$

\end{center}

need not be locally connected in general by a result of Sullivan
\cite{sullivan-lc}, the limit set of the punctured torus groups are
nevertheless locally connected. By extending the analogy of puzzle
pieces, this issue is to an extent clarified. 

\smallskip

An analogue of the ${\Bbb{Z}}_+$ dynamical system may also be
extracted from the split geometry model. Note that each block
corresponds to a splitting of the surface group, and hence an action
on a tree. As $i \rightarrow \infty$, the split blocks $B^s_i$ and
hence the induced splittings also go to infinity, converging to a {\bf
free
action of the surface group on an ${\Bbb{R}}$-tree dual to the ending
lamination}. Thus iteration of the quadratic function correponds to
taking a sequence of splittings of the surface group converging to a
(particular)
action on an ${\Bbb{R}}$-tree. 

\smallskip

{\bf Problem:} The building of the Minsky model and its bi-Lipschitz
equivalence to a hyperbolic manifold \cite{minsky-elc1}
\cite{minsky-elc2} gives rise to a speculation that there should be a
purely combinatorial way of doing much of the work. Bowditch's
rendering \cite{bowditch-endinv}, \cite{bowditch-model}
of the Minsky, Brock-Canary-Minsky results is a step in this
direction. This paper brings out the possibility that the whole thing
should be do-able purely in terms of actions on trees. Of course there
is an action of the surface group on a tree dual to a pants
decomposition. So we do have a starting point. However, one ought to
be able to give a purely combinatorial description, {\em ab initio},
in terms of a sequence of actions of surface groups on trees
converging to an action on an ${\Bbb{R}}$-tree. This would open up the
possibility of extending these results (including those of this paper)
to other hyperbolic groups with
infinite automorphism groups, notably free groups.

\bibliography{amalgeo}
\bibliographystyle{plain}

\end{document}